\documentstyle[epsf,amssymb,12pt]{amsart}

\newcommand{\rrkm}{\stackrel{\mbox{\tiny $r_{k,m}$}}%
{\relbar\joinrel\relbar\joinrel\longrightarrow}}
\newcommand{\rkm}{\stackrel{r_{k,m}}%
{\relbar\joinrel\longrightarrow}}

\newcommand{\Ch}{\mbox{ch}}

\newcommand{\DOT}{\setlength{\unitlength}{1pt}\begin{picture}(2.5,2)(1,1)
\put(1,2){\circle*{2}}\end{picture}}

\newcommand{\Span}[1]{\langle #1 \rangle}
\newcommand{\SPan}[1]{\langle\langle #1 \rangle\rangle}

\newcommand{\Calfdot}{{{\mathcal F}\!_{\DOT}}}
\newcommand{\Gr}{\mbox{\it Grass}}

\newcommand{\Edot}{{E_{\DOT}}}
\newcommand{\Epdot}{{E_{\DOT}}\!'}

\newcommand{\Fdot}{{F\!_{\DOT}}}
\newcommand{\Fpdot}{{{F\!_{\DOT}}'}}
\newcommand{\Fppdot}{{{F\!_{\DOT}}''}}

\newcommand{\Gdot}{{G_{\DOT}}}
\newcommand{\Gpdot}{{{G_{\DOT}}'}}

\newcommand{\llra}{\relbar\joinrel\longrightarrow}

\newcommand{\llllra}{\relbar\joinrel\relbar\joinrel\llra}
\newcommand{\Lra}{\relbar\joinrel\relbar\joinrel%
\relbar\joinrel\relbar\joinrel\llra}

\newcommand{\precdot}{{\prec\!\!\!\cdot\,}}

\newcommand{\subsetdot}{{\subset\!\!\!\!\cdot\,\,}}

\newcommand{\Pdot}{{P\!_{\DOT}}}

\newcommand{\rpp}{{\stackrel{r_{p,n+1-p}}%
{\Lra}}}
\newcommand{\cpp}{{\stackrel{c_p}%
{\llra}}}

\newcommand{\QED}{
\setlength{\unitlength}{1.0pt}%
\begin{picture}(7.5,10)
\put(0,-5){\rule{2.5pt}{5pt}}
\put(0,0){\rule{5pt}{2.5pt}}
\put(0,2.5){\rule{7.5pt}{2.5pt}}
\end{picture}\vspace{10pt}}

\newcommand{\QEDnoskip}{
\setlength{\unitlength}{1.0pt}%
\begin{picture}(7.5,10)
\put(0,-5){\rule{2.5pt}{5pt}}
\put(0,0){\rule{5pt}{2.5pt}}
\put(0,2.5){\rule{7.5pt}{2.5pt}}
\end{picture}}

\newtheorem{lem}{Lemma}[subsection]
\newtheorem{thm}[lem]{Theorem}
\newtheorem{prop}[lem]{Proposition}
\newtheorem{cor}[lem]{Corollary}
\newtheorem{ex}[lem]{Example}

\newtheorem{alg}[lem]{Algorithm}
\newtheorem{rem}[lem]{Remark}
\newtheorem{defn}[lem]{Definition}

\begin{document}

\title[Schubert Polynomials and the Bruhat order]{Schubert polynomials,
the Bruhat order, and the geometry of flag manifolds}

\author{Nantel Bergeron \and Frank Sottile}

\address{Department of Mathematics and Statistics\\
        York University\\
        North York, Ontario M3J 1P3\\
	CANADA}
\email[Nantel Bergeron]{bergeron@@mathstat.yorku.ca}
\address{MSRI\\
	1000 Centennial Drive\\
	Berkeley CA, 94720\\
	USA}
\address{On leave from: Department of Mathematics\\
        University of Toronto\\
        100 St.~George Street\\
	Toronto, Ontario  M5S 3G3\\
	CANADA}
\email[Frank Sottile]{sottile@@msri.org\qquad sottile@@math.toronto.edu}
\date{27 February 1997}
\thanks{First author supported in part by an NSERC grant}
\thanks{Second author supported in part by NSERC grant  OGP0170279 \ 
and NSF grant DMS-9022140}
\thanks{revised version of MSRI preprint \# 1996 - 083}

\begin{abstract}
We illuminate the relation between the Bruhat order and structure
constants for the polynomial ring in terms of its basis of Schubert
polynomials. 
We use  combinatorial, algebraic, and geometric methods, notably a study
of intersections of Schubert varieties and maps between flag manifolds. 
We establish a number of new identities among these
structure constants.
This leads to formulas for some of these 
constants and new results on the enumeration of chains in
the Bruhat order.
A new graded partial order on the symmetric group which contains 
Young's lattice arises from these investigations. 
We also derive formulas for certain specializations of Schubert
polynomials.
\end{abstract}

\maketitle

\mbox{ } 

\begin{center}
To the memory of Marcel Paul  Sch{\"u}tzenberger
\end{center}

\section*{Introduction}
Extending work of Demazure~\cite{Demazure} and of Bernstein, Gelfand, and
Gelfand~\cite{BGG},  in 1982 Lascoux and
Sch\"utzenberger~\cite{Lascoux_Schutzenberger_polynomes_schubert}
defined remarkable polynomial representatives for Schubert classes in
the cohomology of a flag manifold, which they called 
Schubert polynomials.
For each permutation $w$ in ${\cal S}_\infty$, there is
a Schubert polynomial ${\frak S}_w\in{\Bbb Z}[x_1,x_2,\ldots]$.
The collection of all Schubert polynomials forms an additive homogeneous
basis for this polynomial ring.
Thus the identity
	\begin{equation}\label{eq:structure}
	{\frak S}_u \cdot {\frak S}_v \quad 
	=\quad \sum_w c^w_{u\, v} {\frak S}_w
	\end{equation}
defines integral {\em structure constants} $c^w_{u\, v}$ for the ring of
polynomials with respect to its Schubert basis.
Littlewood-Richardson coefficients are
a special case of the $c^w_{u\, v}$ as 
every Schur symmetric polynomial  is a Schubert polynomial.
The  $c^w_{u\, v}$  are positive integers:
Evaluating a Schubert polynomial at certain Chern classes
gives a Schubert class in the cohomology of the flag manifold.
Hence, $c^w_{u\, v}$ enumerates the flags in a
suitable triple intersection of Schubert varieties.
This evaluation exhibits the cohomology of the flag 
manifold as the quotient:
$$
	{\Bbb Z}[x_1,x_2,\ldots]/\Span{{\frak S}_w\,|\, 
	w\not\in {\cal S}_n}.
$$

These constants,  $c^w_{u\, v}$, are readily computed:
The MAPLE libraries ACE~\cite{ACE} include routines for  manipulating
Schubert polynomials, Gr{\"o}bner basis methods are applied in~\cite
{Fomin_Gelfand_Postnikov}, and a new approach, using orbit values of
Kostant polynomials, is developed in~\cite{Billey}.  
However, it remains an open problem to understand these constants 
combinatorially.
By this we mean a combinatorial interpretation or a bijective
formula for these constants.
We expect such a formula will have the form
	\begin{equation}\label{eq:expected_rule}
	c^w_{u\, v}\quad =\quad\# \left\{
	\mbox{\begin{minipage}[c]{3.6in}
	(saturated) chains in the Bruhat order on ${\cal S}_\infty$ from
	\mbox{\ $u$ to $\,w$}  \,satisfying 
	{\em \,some} \,condition \,imposed \,by $\,v$ \end{minipage}}
	\right\}.
	\end{equation}

The Littlewood-Richardson rule~\cite{Littlewood_Richardson}, which this
would generalize, may be expressed in this form
({\em cf.}~\S\ref{sec:chain_description}), as standard Young tableaux are
chains in Young's Lattice, a suborder of the Bruhat order. 
A relation between chains in the Bruhat order and the multiplication of
Schubert polynomials has previously been
noted~\cite{Hiller_intersections}. 
A new proof of the classical Pieri's formula for
Grassmannians~\cite{sottile_explicit_pieri} 
suggests a geometric rationale for such `chain-theoretic' formulas.
Lastly, known formulas for these constants, particularly Monk's
formula~\cite{Monk}, Pieri-type
formulas (first stated by Lascoux and
Sch\"utzenberger~\cite{Lascoux_Schutzenberger_polynomes_schubert}
but only recently given proofs using geometry~\cite{sottile_pieri_schubert}
and algebra~\cite{Winkel_multiplication}), and other
formulas of~\cite{sottile_pieri_schubert}, are all of this form. 
Recently, Ciocan-Fontanine~\cite{Ciocan_partial} has 
generalized these chain-theoretic Pieri-type formulas to the quantum
cohomology rings of manifolds of partial flags.
This has also been announced~\cite{Kirillov_Maeno} in the special case
of quantum Schubert polynomials~\cite{Fomin_Gelfand_Postnikov}.
In \S\ref{sec:chains_and_orders}, we give a
refinement of~(\ref{eq:expected_rule}).

We establish several new identities for the $c^w_{u\,v}$, including
Theorems~\ref{thm:B} ({\em ii}) and~\ref{thm:C} ({\em ii}). 
One, Theorem~\ref{thm:A} ({\em i})({\em b}), 
gives a recursion for $c^w_{u\,v}$
when one of the permutations $wu^{-1}, wv^{-1}$, or $w_0uv^{-1}$ has a fixed
point and a condition on its inversions holds. 
When  ${\frak S}_v$ is a Schur polynomial, 
we give a chain-theoretic interpretation
(Theorem~\ref{thm:skew_shape}) for some $c^w_{u\,v}$, 
determine many more (Theorem~\ref{thm:skew_permutation}) in terms of the
classical Littlewood-Richardson coefficients,  and 
show how a map that takes certain chains in the Bruhat order to
standard Young tableaux and satisfies some additional properties would
give a combinatorial interpretation for these  $c^w_{u\,v}$
(Theorem~\ref{thm:combinatorial}). 

Most of these identities have an order-theoretic companion
which could imply them,
were a description such as (\ref{eq:expected_rule}) known.
The one identity (Theorem~\ref{thm:D}) lacking such a companion yields a
new result about the enumeration of chains in
the Bruhat order (Corollary~\ref{cor:equal_chains}).
In \S\ref{sec:orders}, we study a suborder called the
$k$-Bruhat order, which is relevant in (\ref{eq:expected_rule}) when
${\frak S}_v$ is  Schur symmetric 
polynomial in $x_1,\ldots,x_k$.
This leads to a new graded partial order on ${\cal S}_\infty$
containing every interval in
Young's lattice as an induced suborder for which many group
homomorphisms are order preserving (Theorem~\ref{thm:new_order}).

Some of these identities require the computation of maps on the
cohomology of flag manifolds induced by certain embeddings, including
Theorem~\ref{thm:projection}, Lemma~\ref{lemma:fixed_pts}, 
and Theorem~\ref{thm:many_identities}.  
We use these to determine the effect of some homomorphisms of 
${\Bbb Z}[x_1,x_2,\ldots]$ on its Schubert basis.
For example, let $P\subset {\Bbb N}$ and list 
the elements of $P$ and ${\Bbb N}-P$ in order:
$$
P\ :\ p_1<p_2<\cdots\qquad \qquad
{\Bbb N}-P\ :\ p^c_1<p^c_2<\cdots
$$
Define $\Psi_P:{\Bbb Z}[x_1,x_2,\ldots]\rightarrow
{\Bbb Z}[y_1,y_2,\ldots,z_1,z_2,\ldots]$ by:
$$
\Psi_P(x_{p_j}) \ =\ y_j\qquad\mbox{and}\qquad
\Psi_P(x_{p^c_j}) \ =\ z_j.
$$
Then there exist integers $d^{u\,v}_w$ such that
$$	%\begin{equation}\label{eq:hopf}
\Psi_P({\frak S}_w(x))\ =\ 
\sum_{u,v} d^{u\,v}_w{\frak S}_u(y){\frak S}_v(z).
$$	%\end{equation}
We show (Theorem~\ref{thm:substitution}) there exist 
$\pi\in{\cal S}_\infty$ (depending upon $w$ and $P$) 
such that $d^{u\,v}_w = c^{(u\times v)\pi}_{\pi\;w}$. 
In particular, the coefficients $d^{u\,v}_w$ are nonnegative.
This generalizes 1.5 of~\cite{Lascoux_Schutzenberger_structure_de_Hopf},
where it is shown that the $d^{u\,v}_w$ are non-negative when
$P=\{1,2,\ldots,n\}$.

Algebraic structures in the cohomology of a flag manifold also yield
identities among the  $c^w_{u\,v}$ such as 
$c^w_{u\,v}=c^w_{v\,u}$ (imposed by commutativity) or
$c^w_{u\,v}=c^{w_0u}_{w_0w\:v}=
c^{\overline{w}}_{\overline{u}\:\overline{v}}$,
where $\overline{w}:= w_0 w w_0$, (imposed by Poincar{\'e} duality among
the Schubert classes).
Such `algebraic' identities for the classical Littlewood-Richardson
coefficients were 
studied combinatorially in~\cite{Zelevinsky,%
Berenstein_Zelevinsky,Hanlon_Sundaram,%
benkart_sottile_stroomer_switching,Fulton_tableaux}.
We expect the identities established here will similarly lead to some
beautiful combinatorics, once a combinatorial interpretation for the  
$c^w_{u\,v}$ is known.
These identities impose stringent conditions on the form of
any combinatorial interpretation 
and should be useful in guiding the search for such an
interpretation.

\tableofcontents

\section{Description of results}
\subsection{Suborders of the Bruhat order and the
$c^w_{u\,v}$}\label{sec:chains_and_orders} 
Suppose the Schubert polynomial ${\frak S}_v$ in (\ref{eq:structure}) is
replaced by the Schur polynomial 
$S_\lambda(x_1,\ldots,x_k)$. 
The resulting identity
$$ %	\begin{equation}\label{eq:two}
	{\frak S}_u \cdot S_\lambda(x_1,\ldots,x_k)
	 \quad =\quad \sum_w c^w_{u\, v(\lambda,k)}\, {\frak S}_w
\eqno(1.1.1)
$$%	\end{equation}
defines integer constants $c^w_{u\, v(\lambda,k)}$, which we call {\em
Littlewood-Richardson coefficients for Schubert polynomials}, as 
we show they share many properties
with the classical Littlewood-Richardson coefficients.
They are related to chains
in a suborder of the Bruhat order called the 
{\em $k$-Bruhat order}, $\leq_k$. 
Its covers  coincide with the index
of summation in Monk's formula~\cite{Monk}:
$$		
{\frak S}_u \cdot {\frak S}_{(k,\,k{+}1)} \quad =\quad 	
{\frak S}_u \cdot (x_1+\cdots+x_k)
 \quad =\quad \sum  {\frak S}_{u(a, b)},
$$	
the sum over all $a\leq k<b$ where $\ell(u(a,b))=\ell(u)+1$.
The set of permutations comparable to the identity
in the $k$-Bruhat order is isomorphic to Young's lattice of
partitions  with at most $k$ parts.
This is the set of Grassmannian permutations with descent 
$k$, those permutations whose Schubert polynomials are Schur
symmetric polynomials in $x_1,\ldots,x_k$.
If $f^\lambda$ is the number of standard Young tableaux of shape $\lambda$,
then~\cite[ I.5, Example 2]{Macdonald_symmetric},
$$
 (x_1+\cdots+x_k)^m
 \quad =\quad \sum_\lambda f^\lambda S_\lambda(x_1,\ldots,x_k),
$$
the sum over all $\lambda$ which partition the integer $m$.
Considering the coefficient of ${\frak S}_w$ in the product 
${\frak S}_u\cdot (x_1+\cdots+x_k)^m$ and the 
definition (1.1.1) %\ref{eq:two})
of  $c^w_{u\,v(\lambda,k)}$, we obtain:

\begin{prop}\label{prop:chains}
The number of chains in the $k$-Bruhat order from $u$ to $w$ is
$$
\sum_\lambda f^\lambda c^w_{u\,v(\lambda,k)}.
$$
\end{prop}

In particular, $c^w_{u\,v(\lambda,k)}=0$ unless $u\leq_k w$.
A chain-theoretic description of the constants
$c^w_{u\,v(\lambda,k)}$ should provide a bijective proof of
Proposition~\ref{prop:chains}. 
By this we mean a function $\tau$ from the set of chains in $[u,w]_k$ to
the set of standard Young tableaux $T$ whose shape is a partition of
$\ell(w)-\ell(u)$ with the further condition that whenever $T$ has shape 
$\lambda$, then $\#\tau^{-1}(T)=c^w_{u\,v(\lambda,k)}$.
For the classical Littlewood-Richardson coefficients, 
Schensted insertion~\cite{Schensted} furnishes a
proof~\cite{Thomas_schensted_construction} 
({\em cf.}~\S\ref{sec:chain_description}), as does 
Sch{\"u}tzenberger's  {\em jeu de
taquin}~\cite{Schutzenberger_jeu_de_taquin}.
In~\S\ref{sec:further} we show (Theorem~\ref{thm:combinatorial}) that if
$\tau$ is a function where $\#\tau^{-1}(T)$ depends only upon the shape
of $T$ and satisfies a 
further condition, then $\#\tau^{-1}(T)=c^w_{u\,v(\lambda,k)}$.
Such a function $\tau$ would be a generalization of Schensted insertion
to this setting.

The $k$-Bruhat order has a more intrinsic formulation, which we
establish in \S\ref{sec:k-bruhat}:

\begin{thm}\label{thm:k-length}
Let $u,w\in {\cal S}_\infty$. 
Then $u\leq _k w$ if and only if
\begin{enumerate}
\item[I.] $a\leq k < b$ implies $u(a)\leq w(a)$ and $u(b)\geq w(b)$.
\item[II.] If $\/a<b$, $u(a)<u(b)$, and $w(a)>w(b)$, then 
$a\leq k< b$.
\end{enumerate}
\end{thm}

The $k$-Bruhat order and its connection to the Littlewood-Richardson
coefficients $c^w_{u\,v(\lambda,k)}$ may be generalized, which leads to
a refinement of (\ref{eq:expected_rule}).
A {\em parabolic subgroup} $P$ of 
${\cal S}_\infty$~\cite{Bourbaki_Groupes_IV} is a subgroup generated by
some  adjacent transpositions, $(i,i{+}1)$. 
Given a parabolic subgroup $P$ of ${\cal S}_\infty$, define the 
{\em P-Bruhat order}  by its covers.
A cover $u\lessdot_P w$ in the $P$-Bruhat order is a cover in the Bruhat
order where $u^{-1}w \not\in P$.
When $P$ is generated by all adjacent
transpositions except $(k,k{+}1)$, this is the $k$-Bruhat order.

Let $I\subset \{1,2,\ldots,n{-}1\}$ index the adjacent transpositions
{\em not} in $P$. 
A {\em coloured chain} in the $P$-Bruhat order is a chain together with
an element of $I\bigcap \{a,a{+}1,\ldots,b{-}1\}$ for each cover 
$u\lessdot_P u(a,b)$ in the chain.
This notion of colouring the Bruhat order was introduced
in~\cite{Lascoux_Schutzenberger_symmetry}.
Iterating Monk's rule, we obtain:
$$
%\begin{equation}\label{eq:parabolic_chains}
	\left(\sum_{i\in I}\:{\frak S}_{(i,i{+}1)}\right)^m
	\quad =\quad 
	\sum_v\, f^v_e(P)\; {\frak S}_v,
\eqno(1.1.2)	%\end{equation}
$$
where $f^v_e(P)$ counts the coloured chains in the $P$-Bruhat order
from $e$ to $v$.
This number, $f^v_e(P)$, is nonzero only for those $v$ which are minimal
in their coset $vP$.
More generally, let $f^w_u(P)$ count the coloured chains in the
$P$-Bruhat order from $u$ to $w$.
Multiplying (1.1.2) by ${\frak S}_u$ and equating
coefficients of ${\frak S}_w$, gives a generalization of
Proposition~\ref{prop:chains}:

\begin{thm}\label{thm:chains}
Let $u,w\in{\cal S}_\infty$ and $P$ be any parabolic subgroup of 
${\cal S}_\infty$.
Then
$$
f^w_u(P)\quad =\quad 
\sum_v\, c^w_{u\,v}\, f^v_e(P).
$$
\end{thm}

This also shows $c^w_{u\,v}=0$ unless $u\leq_P w$, whenever $v$ is
minimal in $vP$.
Theorem~\ref{thm:chains} suggests a refinement of
(\ref{eq:expected_rule}):
Let $u,v,w\in{\cal S}_\infty$, and let $P$ be any parabolic subgroup such
that $v$ is minimal in $vP$.
(There always is such a $P$.)
Then,  for every coloured chain $\gamma$ in the $P$-Bruhat order from $e$
to $v$, we expect that
$$ %	\begin{equation}\label{eq:refinement}
	c^w_{u\, v}\quad =\quad\# \left\{
	\mbox{\begin{minipage}[c]{3.6in}
	coloured chains in the $P$-Bruhat order on ${\cal S}_\infty$ from
	$u$ to $w$  which satisfy 
	{\em some} condition imposed by $\gamma$ \end{minipage}}
	\right\}.
\eqno(1.1.3)%	\end{equation}
$$
Moreover, this rule should give a bijective proof of
Theorem~\ref{thm:chains}.

This $P$-Bruhat order may be defined for every parabolic subgroup of every
Coxeter group. 
Likewise, the problem of finding the structure constants for a Schubert
basis also generalizes.
For Weyl groups, the basis is the Schubert classes in the cohomology of 
a generalized flag manifold $G/B$ or the analogues of Schubert
polynomials in this
case~\cite{Billey_Haiman_Schubert,Fomin_Kirillov_Bn,%
Fulton_orthogonal,Pragacz_Ratajski_formulas}.
For finite Coxeter groups, the basis is the `Schubert classes' of
Hiller~\cite{Hiller_schubert} in the coinvariant algebra.
Likewise, Theorem~\ref{thm:chains} and the expectation
(1.1.3) have analogues in this more general setting.
Of the known formulas in this
setting~\cite{Chevalley91,Hiller_Boe,Pragacz_S-Q,%
Stembridge_shifted,Pragacz_Ratajski_Pieri_Odd_I,%
Pragacz_Ratajski_Pieri_Odd_II,Pragacz_Ratajski_Pieri_Even} 
(see also the survey~\cite{Pragacz_divided}), 
few~\cite{Chevalley91,Hiller_Boe,Pragacz_S-Q,%
Stembridge_shifted} have been expressed in such a chain-theoretic
manner.

\subsection{Substitutions and the Schubert basis}
In \S\S\ref{sec:endomorphism} and \ref{sec:fixed_point_identities}, we
study the $c^w_{u\, v}$ when $w(p)=u(p)$ for some $p$.
This leads to a formula for the substitution of 0 for $x_p$
in terms of the Schubert basis, a recursion for some $c^w_{u\, v}$, and
new identities.
For $w\in{\cal S}_{n+1}$ and $1\leq p\leq n+1$, let 
$w/_p\in {\cal S}_n$ be defined by deleting the $p$th row and $w(p)$th
column from the permutation matrix of  $w$.
If $y\in{\cal  S}_n$ and $1\leq q\leq n+1$, then 
$\varepsilon_{p,q}(y)\in {\cal S}_{n+1}$ is the permutation such that
$\varepsilon_{p,q}(y)/_p = y$ and $\varepsilon_{p,q}(y)(p)=q$.
The index of summation in a particular case of the Pieri-type
formula~\cite{Lascoux_Schutzenberger_polynomes_schubert,%
sottile_pieri_schubert,Winkel_multiplication}, 
$$
{\frak S}_v \cdot (x_1\cdots x_{p-1}) \quad =\quad
\sum_{v \cpp w} {\frak S}_w,
$$
defines the relation $v\cpp w$,
which is described in more detail before Theorem~\ref{thm:projection}.
Define $\Psi_p:{\Bbb Z}[x_1,x_2,\ldots] \rightarrow
{\Bbb Z}[x_1,x_2,\ldots]$ by 
$$
\Psi_p(x_j)\ =\ \left\{\begin{array}{ll} x_j&\mbox{ if } j<p\\
0& \mbox{ if } j=p\\
x_{j-1}&\mbox{ if } j>p\end{array}\right..
$$

\begin{thm}\label{thm:A}
Let $u,w\in {\cal S}_\infty$ and $p\in{\Bbb N}$.
\begin{enumerate}
\item[({\em i})]
Suppose $w(p)=u(p)$ and 
$\ell(w)-\ell(u)=\ell(w/_p)-\ell(u/_p)$.
Then
\begin{enumerate}
\item[({\em a})]
$\varepsilon_{p,u(p)} :
[u/_p,w/_p]  \stackrel{\sim}{\longrightarrow}  [u,w]$.
\item[({\em b})] For every $v\in {\cal S}_\infty$,
$$
	c^w_{u\, v}\quad =\quad
	\sum_{\stackrel{\mbox{\scriptsize $y\in{\cal S}_\infty$}}%
	{v \cpp \varepsilon_{p,1}(y)}} c^{w/_p}_{u/_p\: y}.
$$
\end{enumerate}
\item[({\em ii})] 
For every $v\in {\cal S}_\infty$,
$$
 \Psi_p({\frak S}_v) =  
\sum_{\stackrel{\mbox{\scriptsize $y\in{\cal S}_\infty$}}%
	{v \cpp \varepsilon_{p,1}(y)}} {\frak S}_y.
$$
\end{enumerate}
\end{thm}

The first statement (Lemma~\ref{lem:expanding_bruhat} ({\em ii})) is
proven  using combinatorial
arguments, while the second (Theorem~\ref{thm:coeff_sum}) and third
(Theorem~\ref{thm:theorem_A_iii}) are 
proven by computing certain maps on cohomology. 
Since  $c^w_{u\,v}=c^w_{v\,u}=c^{w_0 u}_{v\:w_0w}$, 
Theorem~\ref{thm:A} ({\em i})({\em b}) gives a
recursion for $c^w_{u\, v}$ when one of $wu^{-1}, wv^{-1}$, or 
$w_0uv^{-1}$
has a fixed point and the condition on lengths is satisfied.

We also compute the effect of other substitutions of the variables in
terms of the Schubert basis: 
Let  $P\subset {\Bbb N}$ and list the
elements of $P$ and ${\Bbb N}-P$ in order:
$$
P\ :\ p_1<p_2<\cdots\qquad \qquad
{\Bbb N}-P\ :\ p^c_1<p^c_2<\cdots
$$
Define $\Psi_P:{\Bbb Z}[x_1,x_2,\ldots]\rightarrow
{\Bbb Z}[y_1,y_2,\ldots, z_1,z_2,\ldots]$ by:
$$
\Psi_P(x_{p_j}) \ =\ y_j\qquad\mbox{and}\qquad
\Psi_P(x_{p^c_j}) \ =\ z_j.
$$
In Remark~\ref{rem:I_P}, we define an infinite set $I_P$ of permutations
with the following property:

\begin{thm}\label{thm:substitution}
For every $w\in {\cal S}_\infty$, there exists an integer $N$
such that if $\pi\in I_P$ and $\pi\not\in{\cal S}_N$, then 
$$
\Psi_P({\frak S}_w)\ =\ 
\sum_{u,\, v} 
c^{(u\times v)\cdot \pi}_{\pi\; w}\;
{\frak S}_u(y)\;{\frak S}_v(z).
$$
\end{thm}

A precise version of Theorem~\ref{thm:substitution}
(Theorem~\ref{thm:substitution_constants}) is proven in
\S\ref{sec:substitution}.
Theorem~\ref{thm:substitution} gives infinitely many identities of
the form 
$c^{(u\times v)\cdot \pi}_{\pi\; w} =
c^{(u\times v)\cdot \sigma}_{\sigma\; w}$ for $\pi,\sigma\in I_P$.
Moreover, for these $u,v$, and $\pi$ with 
$c^{(u\times v)\pi}_{\pi\: w}\neq 0$, we have 
$[\pi,\; (u\times v)\cdot \pi]\simeq [e,u]\times[e,\,v]$, 
which is suggestive of a chain-theoretic basis for these identities.

Theorem~\ref{thm:substitution} extends 1.5
of~\cite{Lascoux_Schutzenberger_structure_de_Hopf}, where it is
shown that the $d^{u\,v}_w([n])$ are non-negative.
A combinatorial proof of the non-negativity of these coefficients
$d^{u\;v}_w(P)$ and of Theorem~\ref{thm:A} ({\em ii}) using, perhaps,
one of the combinatorial constructions of Schubert polynomials~\cite{%
Kohnert,Bergeron,BJS,Fomin_Kirillov_YB,%
% 1990   1992    1993           1993          
Fomin_Stanley,Bergeron_Billey,Winkel_kohnert_rule} 
%      1994         1994         1995
may provide insight into the problem of
determining the $c^w_{u\,v}$.

Theorems~\ref{thm:A} ({\em ii}) and ~\ref{thm:substitution} enable the
computation of rather general substitutions:
Let $\Pdot:= (P_0,P_1,\ldots)$ be any (finite or infinite) partition of 
${\Bbb N}$.
For $i>0$, let $\underline{x}^{(i)}:= x^{(i)}_1,x^{(i)}_2,\ldots$
be a set of variables in bijection with  $P_i$.
Define $\Psi_{\Pdot}:{\Bbb Z}[x_1,x_2,\ldots]\rightarrow
{\Bbb Z}[\underline{x}^{(1)},\underline{x}^{(2)},\ldots]$
by
$$
\Psi_{\Pdot}(x_j)\ =\ \left\{\begin{array}{ccl}
0&&\mbox{if } j\in P_0\\
x_l^{(i)}&&\mbox{if $j$ is the $l$th element of $P_i$}
\end{array}\right..
$$

\begin{cor}\label{cor:general_substitution}
For every partition $\Pdot$ of\/ ${\Bbb N}$ and $w\in {\cal S}_\infty$, 
$$
\Psi_{\Pdot}({\frak S}_w(x))\ =\ 
\sum_{u_1,u_2,\ldots} d^{u_1,u_2,\ldots}_w(\Pdot)\:
{\frak S}_{u_1}(\underline{x}^{(1)})
{\frak S}_{u_1}(\underline{x}^{(2)})\cdots, 
$$
where each $d^{u_1,u_2,\ldots}_w(\Pdot)$ is a(n  explicit) sum of
products of the 
$c^z_{v\,y}$, hence non-negative.
\end{cor}

A {\em ballot sequence}~\cite[\S 4.9]{Sagan} $A = (a_1,a_2,\ldots)$ is a
sequence of non-negative integers where, for each $i,j \geq 1$,
$$
\#\{k\leq j\;|\; a_k=i\}\quad \geq \quad
\#\{k\leq j\;|\; a_k=i+1\}.
$$
(Traditionally, the $a_i>0$.  
One should consider $a_i=0$ as a vote for `none of the
above'.) 
Given a ballot sequence $A$, define
$\Psi_A: {\Bbb Z}[x_1,x_2,\ldots]\rightarrow {\Bbb Z}[x_1,x_2,\ldots]$
by
$$
\Psi_A(x_i)\quad =\quad \left\{\begin{array}{ll}
0& a_i=0\\
x_{a_i}& a_i \neq 0\end{array}\right..
$$

\begin{cor}
For every ballot sequence $A$ and $w\in{\cal S}_n$, there exist
non-negative integers $d^u_w(A)$ for $u,w\in {\cal S}_\infty$
such that 
$$
\Psi_A({\frak S}_w(x))\quad=\quad
\sum_u d^u_w(A)\:{\frak S}_u(x).
$$
Moreover, each $d^u_w(A)$ is a(n explicit) sum of products of the
$c^z_{v\,y}$. 
\end{cor}

\noindent{\bf Proof. }
If $P_0:=\{i\;|\; a_i=0\}$ and for $j>0$
$$
P_j \ :=\ 
\{i\;|\;a_i\mbox{ is the $j$th occurrence of some integer in $A$}\},
$$
then $\Psi_A = \Delta \circ \Psi_{(P_0,P_1,\ldots)}$,
where $\Delta$ is the diagonal map, $\Delta(x^{(i)}_j)= x_j$.
\QEDnoskip

\subsection{Identities when ${\frak S}_v$ is a Schur
polynomial}\label{sec:Schur_identities} 

If $\lambda, \mu$, and $\nu$ are partitions with at most $k$ parts then 
the classical Littlewood-Richardson coefficients 
$c^\nu_{\mu\,\lambda}$ are defined by the identity
$$
S_\mu(x_1,\ldots,x_k)\cdot S_\lambda(x_1,\ldots,x_k) \quad = 
\quad \sum_\nu c^\nu_{\mu\,\lambda} S_\nu(x_1,\ldots,x_k).
$$
The $ c^\nu_{\mu\,\lambda}$ depend only upon $\lambda$ and
the skew partition $\nu/\mu$.
That is, if $\kappa$ and $\rho$ are partitions with at most $l$ parts,
and $\kappa/\rho = \nu/\mu$, then 
for all partitions $\lambda$,
$$
c^\nu_{\mu\,\lambda} \quad=\quad
c^\kappa_{\rho\,\lambda}.
$$
Moreover, $c^\kappa_{\rho\,\lambda}$ is the coefficient of 
$S_\kappa(x_1,\ldots,x_l)$ when 
$S_\rho(x_1,\ldots,x_l)\cdot S_\lambda(x_1,\ldots,x_l)$
is expressed as a sum of Schur polynomials.
The order type of the interval in Young's lattice
from $\mu$ to $\nu$ is determined by $\nu/\mu$.
These facts generalize to the Littlewood-Richardson coefficients
$c^w_{u\,v(\lambda,k)}$. 

If $u\leq_k w$, let $[u,w]_k$ be the interval between $u$ and $w$ in the
$k$-Bruhat order, a graded poset.
Permutations $\zeta$ and $\eta$ are {\em shape equivalent} if there
exist sets of integers $P=\{p_1<\cdots<p_n\}$ and $Q=\{q_1<\cdots<q_n\}$,
where $\zeta$ (respectively $\eta$) acts as the identity on 
${\Bbb N}- P$ (respectively ${\Bbb N}- Q$), and 
$$
\zeta(p_i)\ =\ p_j \quad \Longleftrightarrow\quad
\eta(q_i)\ =\ q_j.
$$

\begin{thm}\label{thm:B}
Suppose $u\leq_kw$ and $x\leq_l z$ where
$wu^{-1}$ is shape equivalent to $zx^{-1}$.
Then
\begin{enumerate}
\item[({\em i})] $[u,w]_k\simeq[x,z]_l$.
When $wu^{-1}=zx^{-1}$, this isomorphism is given by $v\mapsto vu^{-1}x$.
\item[({\em ii})] For all partitions $\lambda$,
$c^w_{u\,v(\lambda,k)} = c^z_{x\,v(\lambda,l)}$.
\end{enumerate}
\end{thm}

Theorem~\ref{thm:B} ({\em i}) is a consequence of
Theorems~\ref{thm:k-order} and~\ref{thm:new_order}, which are proven
using combinatorial arguments.
Theorem~\ref{thm:B} ({\em ii}) is proven in \S\ref{sec:proof_thm_B}
using geometric arguments.
By Theorem~\ref{thm:B}, we may define the constant $c^\zeta_\lambda$ for
$\zeta\in {\cal S}_\infty$ and $\lambda$ a partition by 
$c^\zeta_\lambda := c^{\zeta u}_{u\,v(\lambda,k)}$ 
and also define $|\zeta| := \ell(\zeta u) - \ell(u)$ for any 
$u \in {\cal S}_\infty$ with $u\leq_k \zeta u$.
In \S\ref{sec:new:order}, this analysis leads  to a graded 
partial order $\preceq$ on ${\cal S}_\infty$
with rank function  $|\zeta|$
which has the defining property:
Let $[e,\zeta]_\preceq$  be the interval in the $\preceq$-order from the
identity to $\zeta$.  
If $u\leq_k \zeta u$, then the map 
$[e,\zeta]_\preceq \rightarrow [u,\zeta u]_k$ defined by
$$
\eta \quad \longmapsto \quad \eta \, u
$$
is an order isomorphism.
Each interval in Young's lattice is an induced suborder  in 
$({\cal S}_\infty,\preceq)$ rooted at the
identity, as is the
lattice of partitions with at most $l$ parts
(these are embedded differently for different $l$).
Proposition~\ref{prop:chains} may be stated in terms of this order:
$\sum_\lambda f^\lambda c^\zeta_\lambda$ counts the 
chains in $[e,\zeta]_\preceq$.
This order is studied further in~\cite{bergeron_sottile_order},
where an upper bound is given for $c^\zeta_\lambda$.

For $\eta\in{\cal S}_n$, the map $\eta\mapsto\overline{\eta}$ induces an
order isomorphism $[e,\zeta]_\preceq \stackrel{\sim}{\rightarrow}
[e,\overline{\zeta}]_\preceq$.
The involution ${\frak S}_w \mapsto {\frak S}_{\overline{w}}$ shows
$c^\zeta_\lambda = c^{\overline{\zeta}}_{\lambda^t}$, where 
$\lambda^t$ is the conjugate or transpose of $\lambda$.
Thus $[e,\zeta]_\preceq \simeq [e,\eta]_\preceq$ is not sufficient to
guarantee $c^\zeta_\lambda = c^\eta_\lambda$.
We express some of the Littlewood-Richardson coefficients in
terms of chains in the Bruhat order.
If $u\lessdot_k u(a,b)$ is a cover in the $k$-Bruhat order, label that edge
of the Hasse diagram with the integer $u(b)$.
The {\em word} of a chain in the $k$-Bruhat order is the sequence of
labels of edges in the chain.

\begin{thm}\label{thm:skew_shape}
Suppose $u\leq_k w$ and $wu^{-1}$  is shape equivalent to
$v(\mu,l)\cdot v(\nu,l)^{-1}$, for some $l$ and partitions $\mu,\nu$.
Then, for all partitions $\lambda$ and standard Young tableaux $T$ of
shape $\lambda$, 
$$
c^w_{u\,v(\lambda,k)}\ =\ 
\#\left\{\begin{array}{cc}\mbox{Chains in $k$-Bruhat order from $u$ to
$w$ whose word} \\\mbox{has recording tableau $T$ for Schensted
insertion} \end{array}\right\}.
$$
\end{thm}

Theorem~\ref{thm:skew_shape} gives a combinatorial proof of
Proposition~\ref{prop:chains} for some $u,w$.
We prove this in \S\ref{sec:chain_description} using
Theorem~\ref{thm:B} and combinatorial arguments.

If a skew partition $\theta=\rho\coprod\sigma$ is the union
of incomparable skew partitions $\rho$ and $\sigma$, then
$$
\rho\coprod\sigma\ \simeq\ \rho\times\sigma,
$$
as graded posets.
The skew Schur function $S_\theta$ is
defined~\cite[I.5]{Macdonald_symmetric}  to be 
$\sum_\lambda c^\theta_\lambda\, S_\lambda$ and 
$S_{\rho\coprod \sigma}=S_\rho\cdot
S_\sigma$~\cite[I.5.7]{Macdonald_symmetric}.
Thus 
$$	%\begin{equation}\label{eq:disjoint_skew}
	c^{\rho\coprod\sigma}_\lambda\quad =\quad \sum_{\mu,\,\nu}
	c^\lambda_{\mu\,\nu}\, c^\rho_\mu\, c^\sigma_\nu.
\eqno(1.3.1)%	\end{equation}
$$
Permutations $\zeta$ and $\eta$  are {\em disjoint} if 
$\zeta$ and $\eta$ have disjoint supports and
$|\zeta\eta|=|\zeta|+|\eta|$. 

\begin{thm}\label{thm:C}
Let $\zeta$ and $\eta$ be disjoint permutations.
Then
\begin{enumerate}
\item[({\em i})]  The map 
$(\xi,\chi)\mapsto \xi\chi$ induces an isomorphism of 
graded posets 
$$
[e,\zeta]_\preceq\times[e,\eta]_\preceq
\stackrel{\sim}{\longrightarrow} 
[e,\zeta\eta]_\preceq.
$$
\item[({\em ii})] For every partition $\lambda$,\  
${\displaystyle
c^{\zeta\eta}_\lambda\ =\ \sum_{\mu,\,\nu}
c^\lambda_{\mu\,\nu}\, c^\zeta_\mu\,  c^\eta_\nu.
}$
\end{enumerate}
\end{thm}

The first statement is proven in \S\ref{sec:disjoint_permutations}
(Theorem~\ref{thm:disjoint_iso}), using a characterization of
disjointness related to non-crossing partitions~\cite{Kreweras}, and the
second in \S\ref{sec:proof_C} using geometry.

Our last identity has no analogy with the
classical Littlewood-Richardson coefficients.
Let $(1\,2\,\ldots\, n)$ be the permutation which cyclicly
permutes $[n]$.

\begin{thm}[Cyclic Shift]\label{thm:D}
Suppose $\zeta\in S_n$ and $\eta = \zeta^{(1\,2\,\ldots\,n)}$.
Then, for every partition $\lambda$, 
$c^\zeta_\lambda = c^\eta_\lambda$.
\end{thm}

This is proven in \S\ref{sec:thmd} using geometry.
Combined with Proposition~\ref{prop:chains},
we obtain:

\begin{cor}\label{cor:equal_chains}
If $u\leq_k w$ and $x\leq_k z$ with $wu^{-1},zx^{-1}\in{\cal S}_n$
and 
$(wu^{-1})^{(1\,2\,\ldots\,n)} = zx^{-1}$, then 
the two intervals $[u,w]_k$ and $[x,z]_k$ each have the same number of
chains.
\end{cor}

The two intervals
$[u,w]_k$ and $[x,z]_k$ of Corollary~\ref{cor:equal_chains} are
typically non-isomorphic: 
For example, in ${\cal S}_4$ let $u=1234$, $x=2134$, and $v=1324$.
If $\zeta=(1243)$, $\eta=(1423)= \zeta^{(1234)}$, and 
$\xi=(1342)=\eta^{(1234)}$, then 
$$
u \ \leq_2\ \zeta u,\quad x\ \leq_2\  \eta x,\quad\mbox{and}
\quad v\ \leq_2 \ \xi v.
$$
Here are the intervals $[u,\zeta u]_2$,  $[x,\eta x]_2$,
and  $[v,\xi v]_2$.
\begin{figure}[htb]
$$\epsfxsize=4in \epsfbox{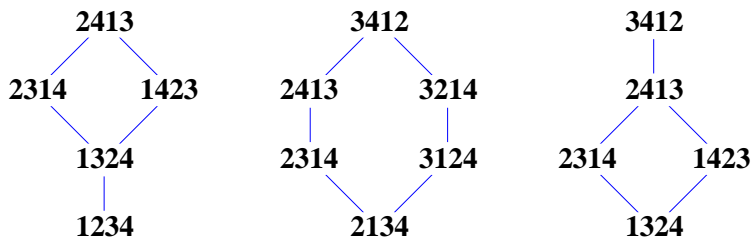}$$
\caption{Effect of cyclic shift on intervals\label{fig:fig1}}
\end{figure}

The Theorems of this section, together with the `algebraic' identities
$c^w_{u\,v}=c^{w_0w}_{w_0u\:v}=
c^{\overline{w}}_{\overline{u}\:\overline{v}}$,
greatly reduce the number of distinct Littlewood-Richardson coefficients
$c^w_{u\; v(\lambda,k)}$ which need to be determined.
We make this precise.
For $\lambda$ a partition, let $\lambda^t$ be its conjugate, or
transpose.
Note $\overline{v(\lambda,k)}= v(\lambda^t,n-k)$.

\begin{thm}[Symmetries of the $c^w_{u\; v(\lambda,k)}$]
Let $\lambda$ be a partition and $\zeta\in{\cal S}_n$.
Then
$$
c^\zeta_\lambda \ =\ c^{\zeta^{-1}}_{\lambda^t}\ =\ 
c^{\overline{\zeta}}_{\lambda^t}\ =\ 
c^{\zeta^{(1\,2\,\ldots\,n)}}_\lambda.
$$
\end{thm}

Let $D_n$ be the dihedral group with $2n$ elements.
These identities show that the action of 
${\Bbb Z}/2{\Bbb Z}\times D_n$ on these coefficients
leaves their values invariant.

\section{Preliminaries}

\subsection{Permutations}\label{sec:permutations}
Let ${\cal S}_n$ be the group of permutations of
$[n] := \{1,2,\ldots,n\}$. 
Let $(a,b)$ be the transposition interchanging $a<b$.
The {\em length\/} $\ell(w)$ of a permutation $w\in {\cal S}_n$
counts the {\em inversions\/}, $\{i<j\,|\, w(i)>w(j)\}$, of $w$.
The Bruhat order $\leq$ on ${\cal S}_n$ is the partial order
whose cover 
relation is $w \lessdot w(a,b)$ if $w(a)<w(b)$ and whenever $a<c<b$,
either $w(c)<w(a)$, or $w(b)<w(c)$.
Thus $\ell(w)+1=\ell(w(a,b))$, so the Bruhat order is
graded by length with minimal element  the identity, $e$. 
If $u\leq w$, let $[u,w]:=\{v\,|\, u\leq v\leq w\}$ be the interval
between $u$ and $w$ in ${\cal S}_n$, a poset graded by
$\ell(v)-\ell(u)$. 
Let $w_0^{(n)}\in {\cal S}_n$ (or simply $w_0$) be
defined by $w_0(j)=n+1-j$.

A permutation $w\in {\cal S}_n$ acts on $[n{+}1]$, fixing $n+1$.
Thus ${\cal S}_n\subset{\cal S}_{n+1}$.
Define ${\cal S}_\infty:= \bigcup_n {\cal S}_n$, the permutations of
the positive integers ${\Bbb N}$ fixing all but finitely many integers.
For $P=\{p_1<p_2<\cdots\}\subset {\Bbb N}$,
define $\phi_P: {\cal S}_{\#P} \rightarrow {\cal S}_\infty$
by requiring that $\phi_P$ act as
the identity on ${\Bbb N}-P$ and 
 $\phi_P(\zeta)(p_i)=p_{\zeta(i)}$.
This  injective homomorphism is a map of posets, but not of graded
posets, as $\phi_P$ typically does not preserve length. 
If $P=\{n+1,n+2,\ldots\}$, then $\phi_P$ does preserve length.
For this $P$, set $1^n\times w:= \phi_P(w)$.
If there exist permutations $\xi,\zeta$, and $\eta$ and sets of positive
integers $P,Q$ such that $\phi_P(\xi)=\zeta$ and $\phi_Q(\xi)=\eta$,
then $\zeta$ and $\eta$ are {\em shape equivalent}.

\subsection{Schubert polynomials}\label{sec:schubert}
Lascoux and Sch\"utzenberger invented and then developed the elementary
theory of Schubert  polynomials in a series of
papers~\cite{Lascoux_Schutzenberger_polynomes_schubert,%
Lascoux_Schutzenberger_structure_de_Hopf,%
Lascoux_Schutzenberger_symmetry,%
Lascoux_Schutzenberger_interpolation,%
Lascoux_Schutzenberger_schub_LR_rule,%
Lascoux_Schutzenberger_operators}. %,Lascoux_Schutzenberger_grothendieck
A self-contained exposition of some of this elegant theory is found
in~\cite{Macdonald_schubert}. 
For an interesting historical survey, see~\cite{Lascoux_historique}.

${\cal S}_n$ acts on polynomials in $x_1,\ldots,x_n$ 
by permuting the variables.
For a polynomial $f$, 
$f - (i,i{+}1) f$ is antisymmetric in $x_i$ and $x_{i+1}$, hence
divisible by $x_i - x_{i+1}$.
Define the {\em divided difference} operator 
$$
\partial_i \ :=\  (x_i-x_{i+1})^{-1} (e - (i,i{+}1)).
$$
If $w = (a_1, a_1{+}1)\cdots (a_p,a_p{+}1)$ is a factorization of $w$ 
into adjacent transpositions with minimal length ($p = \ell(w)$),
then $\partial_{a_1}\circ\cdots\circ\partial_{a_p}$
depends only upon $w$,
defining an operator $\partial_{w}$ for each $w\in {\cal S}_n$.
For $w \in {\cal S}_n$, Lascoux and
Sch\"utzenberger~\cite{Lascoux_Schutzenberger_polynomes_schubert}
defined the {\em Schubert polynomial}  ${\frak S}_{w}$ by
$$
{\frak S}_{w}\ :=\ \partial_{w^{-1}w_0}
\left( x_1^{n-1} x_2^{n-2}\cdots x_{n-1} \right).
$$
The degree of $\partial_i$ is $-1$, so ${\frak S}_{w}$ is homogeneous of 
degree ${n\choose 2} - \ell(w^{-1}w_0) = \ell(w)$.
Since $w_0^{(n)}= (n,n{+}1)\cdots(2,3)(1,2)w_0^{(n+1)}$
and $x_1^{n-1}\cdots x_{n-1}= 
\partial_n\circ\cdots\circ\partial_1
(x_1^n \cdots x_{n-1}^2 x_n)$,
${\frak S}_w$ is independent of $n$ (when $w\in {\cal S}_n$).
This defines polynomials ${\frak S}_w\in{\Bbb Z}[x_1,x_2,\ldots]$ for 
$w\in{\cal S}_\infty$.

A {\em partition} $\lambda$ is a decreasing sequence 
$\lambda_1\geq\lambda_2\geq\cdots\geq\lambda_k\geq 0$ of integers.
Each $\lambda_j$ is a {\em part} of $\lambda$.
For partitions $\lambda$ and $\mu$, write $\mu\subset \lambda$ if 
$\mu_i\leq \lambda_i$ for all $i$.
{\em Young's lattice} is the set of partitions ordered by $\subset$.
The partition with $n$ parts each equal to $m$ is written $m^n$.
For a partition $\lambda$ with $\lambda_{k+1}=0$, the Schur
polynomial $S_\lambda(x_1,\ldots,x_k)$ is
$$
S_\lambda(x_1,\ldots,x_k) \ :=\  
\frac{\det\left|x_j^{k-i+\lambda_i}\right|_{i,j=1}^k}%
{\det\left|x_j^{k-i}\right|_{i,j=1}^k}.
$$
$S_\lambda(x_1,\ldots,x_k)$ is symmetric in 
$x_1,\ldots,x_k$ and  homogeneous of degree
$|\lambda|:=\lambda_1+\cdots+\lambda_k$.

A permutation $w$ is {\em Grassmannian of descent $k$} if
$j\neq k \Rightarrow w(j)<w(j+1)$.
A Grassmannian permutation $w$ with descent $k$ defines, and is defined
by, a partition $\lambda$ with $\lambda_{k+1}=0$:
$$
\lambda_{k+1-j}\ =\ w(j) - j\qquad j=1,\ldots, k.
$$
(The condition $w(k{+}1)<w(k{+}2)<\cdots$ determines the remaining
values of $w$.)
In this case, write $w = v(\lambda,k)$.
The {\em raison d'etre} for this  definition is that
${\frak S}_{v(\lambda,k)} = S_\lambda(x_1,\ldots,x_k)$.
Thus the Schubert polynomials form a basis for 
${\Bbb Z}[x_1,x_2,\ldots]$ which contains all Schur symmetric
polynomials $S_\lambda(x_1,\ldots,x_k)$ for all $\lambda$ and $k$.

\subsection{The flag manifold}\label{sec:flag}

Let $V\simeq {\Bbb C}^n$.
A {\em flag} $\Fdot$ in $V$ is a sequence
$$
\{0\}\  =\ F_0 \subset F_1 \subset F_2\subset \cdots \subset F_{n-1}
\subset F_n\ =\ V,
$$
of subspaces with $\dim_{{\Bbb C}} F_i = i$.
Flags $\Fdot$ and $\Fpdot$ are {\em opposite} if
$F_{n-j}\cap F'_j = \{0\}$ for all $j$.
The set of all flags is an ${n\choose 2}$-dimensional complex 
manifold, ${\Bbb F}\ell V$ (or ${\Bbb F}\ell_n$), called the 
{\em flag manifold}.
There is a tautological flag $\Calfdot$ of 
bundles over ${\Bbb F}\ell V$ whose fibre at $\Fdot$ is $\Fdot$.
Let $-x_i$  be the first Chern class of the line bundle 
${\cal F}_i/{\cal F}_{i-1}$.
Borel~\cite{Borel_cohomology} showed the cohomology ring of 
${\Bbb F}\ell V$ to be
$$
{\Bbb Z}[x_1,\ldots,x_n]/\Span{e_i(x_1,\ldots,x_n)\,|\,i=1,\ldots,n},
$$
where $e_i(x_1,\ldots,x_n)$ is the $i$th elementary symmetric polynomial
in $x_1,\ldots, x_n$.  

Given a subset $S \subset V$, let $\Span{S}$ be its linear span.
For subspaces $W\subset U$, let $U- W$ be their 
set-theoretic difference.
An ordered 
basis $f_1,f_2,\ldots,f_n$ for $V$ determines a flag
$\Edot:=\SPan{f_1,\ldots,f_n}$, where 
$E_i = \Span{f_1,\ldots,f_i}$ for $1\leq i \leq n$.
A fixed flag $\Fdot$ gives a decomposition due to 
Ehresmann~\cite{Ehresmann} of ${\Bbb F}\ell V$ into 
affine cells indexed by permutations $w$ of $S_n$.
The cell determined by $w$ has two equivalent descriptions:
$$%	\begin{equation}	%   remove for the published version
\label{eq:Schubert_cell_definition}
X^{\circ}_w \Fdot \quad  := \quad 
\left\{\begin{array}{ll}
\{ \Edot\in {\Bbb F}\ell V\,|\,
\dim E_i\bigcap F_j = \#\{p\leq i\,|\, w(p)>n-j\}\}, \\ 
\{ \Edot=\SPan{f_1,\ldots,f_n}\,|\,
f_i \in F_{n+1-w(i)}- F_{n-w(i)}, \,1\leq i\leq n\}.\rule{0pt}{15pt}
\end{array}\right.
$$ %	\end{equation}
Its closure is the Schubert subvariety $X_w\Fdot$,
which has complex codimension $\ell(w)$.
Also, $u\leq w \Leftrightarrow X_u\Fdot \supset X_w\Fdot$.
The {\em Schubert class} $[X_w\Fdot]$ is the cohomology class
Poincar{\'e} dual to the fundamental cycle of $X_w\Fdot$.
These Schubert classes form a basis for the cohomology. 
Schubert polynomials were defined so that 
${\frak S}_w(x_1,\ldots,x_n)=[X_w\Fdot]$.
We write ${\frak S}_w$ for $[X_w\Fdot]$.

If $\Fdot$ and $\Fpdot$ are opposite flags, then 
$X_u\Fdot\bigcap X_v\Fpdot$ is an irreducible, generically transverse
intersection, a consequence of~\cite{Deodhar} 
({\em cf.}~\cite[\S5]{sottile_pieri_schubert}). 
Thus its codimension is $\ell(u)+\ell(v)$, and 
the fundamental cycle of $X_u\Fdot\bigcap X_v\Fpdot$ is Poincar\'e dual to 
${\frak S}_u\cdot{\frak S}_v$. 
Since 
$$
{\Bbb Z}[x_1,\ldots,x_n]\  \longrightarrow \ 
%{\Bbb Z}[x_1,\ldots,x_{n+m}] \twoheadrightarrow 
{\Bbb Z}[x_1,\ldots,x_{n+m}]/\Span{e_i(x_1,\ldots,x_{n+m})
\,|\,i=1,\ldots,n+m},
$$
is an isomorphism on 
${\Bbb Z}\Span{x_1^{a_1}\cdots x_n^{a_n}\,|\, a_i<m}$, 
identities of Schubert polynomials follow from product 
formulas for Schubert classes.
The Schubert basis is self-dual for 
the intersection pairing: 
If $\ell(w) + \ell(v) = {n\choose 2}$, then
$$
{\frak S}_w\cdot {\frak S}_v \quad=\quad\left\{
\begin{array}{lll} {\frak S}_{w_0}&\ & v=w_0 w\\0&&\mbox{otherwise}
\end{array}\right..
$$

Let $\Gr_kV$ be the Grassmannian of $k$-dimensional subspaces 
of $V$, a $k(n{-}k)$-dimensional manifold.
A flag $\Fdot$ induces a cellular decomposition indexed by
partitions $\lambda\subset (n{-}k)^k$.
The closure of the cell indexed by $\lambda$ is the Schubert variety
$\Omega_\lambda\Fdot$:
$$
\Omega_\lambda\Fdot\quad:=\quad
\{H\in \Gr_kV\,|\, \dim H\bigcap F_{n+j-k-\lambda_j} \geq j,\ 
j=1,\ldots,k\}.
$$
The cohomology class Poincar{\'e} dual to the fundamental cycle of
$\Omega_\lambda\Fdot$ is $S_\lambda(x_1,\ldots,x_k)$,
where $x_1,\ldots,x_k$ are negative Chern roots
of the tautological $k$-plane bundle on $\Gr_kV$.
Write $S_\lambda$ for $S_\lambda(x_1,\ldots,x_k)$, if $k$ is
understood.
As with the flag manifold, these {\em Schubert classes} form a basis for
cohomology,  $\mu\subset\lambda\Leftrightarrow
\Omega_\mu\Fdot\supset\Omega_\lambda\Fdot$,
and if $\Fdot,\Fpdot$ are opposite flags, then 
$$
[\Omega_\mu\Fdot \bigcap \Omega_\nu\Fpdot]\ =\ 
[\Omega_\mu\Fdot]\ \cdot\ [\Omega_\nu\Fpdot]\ =\ 
\sum_{\lambda\subset (n-k)^k} c^\lambda_{\mu\,\nu}\, S_\lambda,
$$
where the $c^\lambda_{\mu\,\nu}$ are the Littlewood-Richardson
coefficients~\cite{Fulton_tableaux}.

This Schubert basis is self-dual:
If $\lambda\subset (n-k)^k$, then let $\lambda^c$, the {\em complement} of
$\lambda$,  be the partition $(n-k-\lambda_k,\ldots,n-k-\lambda_1)$. 
Suppose $|\lambda|+|\mu|=k(n-k)$, then 
$$
S_\lambda(x_1,\ldots,x_k)\cdot S_\mu(x_1,\ldots,x_k)
\quad=\quad
\left\{
\begin{array}{ll}  S_{(n-k)^k} & \mbox{ if } \mu=\lambda^c\\
0 & \mbox{ otherwise}
\end{array} \right.   .
$$
We suppress the dependence of $\lambda^c$ on $n$ and $k$, which may be
determined by context.

A map $f: X\rightarrow Y$ between manifolds
induces a homomorphism $f_* : H^*X \rightarrow H^*Y$
of abelian groups via the functorial map on homology and the
Poincar{\'e} duality isomorphism between homology and cohomology.
While $f_*$ is not a map of graded rings, it does satisfy the projection
formula ({\em cf.}~\cite[8.1.7]{Fulton_intersection}):
Let $\alpha\in H^*X$ and $\beta\in H^*Y$, then
$$ %	\begin{equation}\label{eq:projection}
	f_*(f^*\alpha \cap \beta)\quad=\quad
	\alpha\cap f_*\beta.
\eqno(2.3.1)     %\end{equation}
$$
For a(n oriented) manifold $X$ of dimension $d$, 
$H^dX={\Bbb Z}\cdot[\mbox{pt}]$ is generated by the class of a point.
Let $\deg : H^*X\rightarrow {\Bbb Z}$ be the map which selects the
coefficient of $[\mbox{pt}]$
Then $\deg(f_*\beta)=\deg(\beta)$.

Let $\pi_k:{\Bbb F}\ell V\twoheadrightarrow \Gr_kV$ be defined by
$\pi_k(\Edot)=E_k$.
Then $\pi_k^{-1}\Omega_\lambda\Fdot = X_{v(\lambda,k)}\Fdot$ and 
$\pi_k: X_{w_0v(\lambda^c,k)}\Fdot\twoheadrightarrow\Omega_\lambda\Fdot$
is generically one-to-one.
Thus on cohomology,
\begin{eqnarray*}
\pi^* S_\lambda&=& {\frak S}_{v(\lambda,k)}\\
(\pi_k)_* {\frak S}_w &=&\left\{\begin{array}{lll}
S_\lambda&\ & \mbox{if }w = w_0 v(\lambda^c,k)\\
0&& \mbox{otherwise}\end{array}\right..
\end{eqnarray*}

By the K{\"u}nneth formula, the cohomology of 
${\Bbb F}\ell V\times{\Bbb F}\ell W$ ($\dim W=m$) has an integral basis
of classes ${\frak S}_u\otimes{\frak S}_x$ for $u\in{\cal S}_n$ and 
$x\in {\cal S}_m$.
Likewise the cohomology of $\Gr_kV\times \Gr_l W$ has a basis
$S_\lambda\otimes S_\mu$ for $\lambda\subset(n-k)^k$ and
$\mu\subset(m-l)^l$. 

While we use the cohomology rings of complex varieties,
our results and methods are valid for the Chow
rings~\cite{Fulton_intersection} and $l$-adic ({\'e}tale)
cohomology~\cite{Deligne_SGA4.5} of these same varieties over
any field.

\section{Orders on ${\cal S}_\infty$}\label{sec:orders}
\subsection{The $k$-Bruhat order}\label{sec:k-bruhat}
The $k$-Bruhat order,  $\leq_k$,  is a suborder
of the  Bruhat order on ${\cal S}_\infty$
which is linked to the Littlewood-Richardson coefficients 
$c^w_{u\,v(\lambda,k)}$.
It appeared 
in~\cite{Lascoux_Schutzenberger_symmetry}, where it was called the
$k$-coloured Ehresmano{\"e}dre.
Its covers are given by the index of summation in Monk's
formula~\cite{Monk}:
$$
{\frak S}_u \cdot (x_1+\cdots+x_k)\quad =\quad 
\sum_{\stackrel{\mbox{\scriptsize $u\leq_k w$}}{\ell(w)=\ell(u)+1}}
{\frak S}_w
$$
Thus $w$ covers $u$ in the $k$-Bruhat order if $w$ covers $u$ in the
Bruhat order, so that $w = u(a,b)$ and 
$\ell(w)=\ell(u)+1$, with the additional requirement that $a\leq k<b$.
The $k$-Bruhat order has the following non-recursive characterization.
\medskip

\noindent{\bf Theorem~\ref{thm:k-length}.}\ 
{\em
Let $u,w\in {\cal S}_\infty$. 
Then $u\leq_k w$ if and only if
\begin{enumerate}
\item[I.] $a\leq k < b$ implies $u(a)\leq w(a)$ and $u(b)\geq w(b)$.
\item[II.] If $a<b$,  $u(a)<u(b)$,  and $w(a)>w(b)$, then $a\leq k<b$.
\end{enumerate}
}\medskip

\noindent{\bf Proof.  }
The idea is to show that the transitive relation $u\trianglelefteq_k w$
defined by the conditions I and II  coincides with the $k$-Bruhat order.
If $u\lessdot _k u(a,b)$ is a cover, then 
$u\trianglelefteq_k u(a,b)$.
Thus $u\leq _k w$ implies $u\trianglelefteq_k w$.
Algorithm~\ref{alg:chain},
which, given $u\trianglelefteq_k w$
produces a chain in the $k$-Bruhat order from $u$ to $w$,
completes the proof. 
\QEDnoskip

\begin{alg}[Produces a chain in the $k$-Bruhat order]\label{alg:chain}
\mbox{ }

\noindent{\tt input: }Permutations $u,w\in {\cal S}_\infty$ with 
$u\trianglelefteq_k w$.

\noindent{\tt output: }A chain in the $k$-Bruhat order from $w$ to $u$.

Output $w$.
While $u\neq w$, do
\begin{enumerate}
\item[1] Choose $a\leq k$ with $u(a)$ minimal subject to $u(a)< w(a)$.
\item[2] Choose $k< b$ with $u(b)$ maximal subject to $w(b)<w(a)\leq u(b)$.
\item[3] $w:=w(a,b)$, output $w$.   
\end{enumerate}

At every iteration of\/ {\rm 1}, $u\trianglelefteq_k w$.
Moreover, this algorithm terminates in $\ell(w)-\ell(u)$ iterations and the
sequence of permutations produced is a  chain in the $k$-Bruhat
order from $w$ to $u$.
\end{alg}

\noindent{\bf Proof. }
It suffices to consider a single iteration.
We first show it is possible to choose $a$ and $b$, then
$u\trianglelefteq_k w(a,b)$, 
and lastly $w(a,b)\lessdot_kw$ is a cover in the $k$-Bruhat order.

In 1, $u\neq w$, so one may always choose $a$.
Suppose $u\trianglelefteq_k w\in {\cal S}_n$ and 
it is not possible to choose $b$.
In that case, if $j>k$ and $w(j)<w(a)$, then
also $u(j)<w(a)$.
Similarly, if $j\leq k$ and $w(j)<w(a)$, then $u(j)\leq w(j)<w(a)$.
Thus $\alpha<w(a)\Leftrightarrow uw^{-1}(\alpha)<w(a)$, which 
contradicts $uw^{-1}(w(a))=u(a)< w(a)$.

Let  $w':= w(a,b)$.
Note that $w(b)\geq u(a)$ implies Condition I for $(u,w')$.
Suppose $w(b)< u(a)$.
Set $b_1:= u^{-1}w(b)$.
Then $w(b_1)\neq u(b_1)$ and the minimality of $u(a)$ shows that $b_1>k$ and
$w(b_1)<u(b_1)$. 
Similarly, if $b_2:= u^{-1}w(b_1)$, then $b_2>k$ and
$w(b_2)<u(b_2)$. 
Continuing, we obtain a sequence $b_1,b_2,\ldots$ with 
$u(a)>u(b_1)>u(b_2)>\cdots$, a contradiction.

We show $(u,w')$ satisfies Condition II.
Suppose $i<j$ and $u(i)<u(j)$.
If $j\leq k$, then $w(i)<w(j)$.
To show $w'(i)<w'(j)$, it suffices to consider the case $j=a$.
But then $u(i)<u(a)$, and thus $u(i)=w(i)=w'(i)$,
by the minimality of $u(a)$.
Then $w'(i)<u(a)\leq w(b)=w'(a)$.
Similarly, if $k<i$, then  $w'(i)<w'(j)$.

Finally, suppose $w$ does not cover $w'$ in the $k$-Bruhat order.
Since $w(a)>w(b)$, there exists a $c$ with 
$a<c<b$ and $w(a)>w(c)>w(b)$.
If $k<c$, then Condition II implies $u(c)>u(b)$
and then the maximality of $u(b)$ implies $w(a)<w(c)$, a contradiction. 
The case $c\leq k$ similarly leads to a contradiction. 
\QEDnoskip

\begin{rem}
{\em 
 Algorithm~\ref{alg:chain} depends only upon $\zeta=wu^{-1}$.
}

\noindent{\tt input: }A permutation $\zeta\in {\cal S}_\infty$.

\noindent{\tt output: }Permutations
$\zeta,\zeta_1,\ldots,\zeta_m=e$ such that  if $u\leq_k \zeta u$, then
$$
u\ \lessdot_k\ \zeta_{m-1}u\ \lessdot_k\ \cdots\ 
\lessdot_k\ \zeta_1u\ \lessdot_k\ \zeta u\ (=w)
$$
is a saturated chain in the $k$-Bruhat order.

Output $\zeta$.
While $\zeta\neq e$, do
\begin{enumerate}
\item[1] Choose $\alpha$ minimal subject to $\alpha< \zeta(\alpha)$.
\item[2] Choose $\beta$ maximal subject to 
$\zeta(\beta)<\zeta(\alpha)\leq\beta$. 
\item[3] $\zeta:=\zeta(\alpha,\beta)$,  output $\zeta$.
\end{enumerate}

{\em 
To see this is equivalent to Algorithm~\ref{alg:chain}, 
set $\alpha=u(a)$ and $\beta=u(b)$ so that $w(a)=\zeta(\alpha)$ and
$w(b)=\zeta(\beta)$.
Thus 
$w(a, b) = \zeta u(a, b) = \zeta(\alpha,\beta) u$.
}
\end{rem}

More is true, the full interval $[u,w]_k$
depends only upon $wu^{-1}$:

\begin{thm}\label{thm:k-order}
If $u\leq_kw$ and $x\leq_ky$ with $wu^{-1}=zx^{-1}$, then 
the map $v\mapsto v u^{-1}x$ induces an isomorphism of graded posets 
$[u,w]_k \stackrel{\sim}{\longrightarrow}[x,z]_k$.
\end{thm}

This is a consequence of the following lemma.

\begin{lem}\label{lem:k-order}
Let  $u\leq_kw$ and $x\leq_k z$ with $wu^{-1}=zx^{-1}$.
If $u\lessdot_k (\alpha,\beta)u$ is a cover with
$(\alpha,\beta)u \leq_k w$, then 
$x\lessdot_k (\alpha,\beta)x$ is a cover with
$(\alpha,\beta)x \leq_k z$
\end{lem}

\noindent{\bf Proof. }
Let $\zeta=wu^{-1}=zx^{-1}$.
By the position of $\gamma$ in $u$, we mean $u^{-1}(\gamma)$.

Suppose $(\alpha,\beta)x$ does not
cover $x$ in the $k$-Bruhat order,
so there is a $\gamma$ with $\alpha<\gamma<\beta$ and 
$x^{-1}(\alpha)<x^{-1}(\gamma)<x^{-1}(\beta)$.
Then, in one line notation, $x$ and $z$ are as illustrated:
$$
\begin{array}{lc}
x:\ &\ldots\ \ \alpha\ \ \ldots\ \ \gamma\ \ \ldots\ \ \beta\ \ \ldots\\
z:\ &\ldots\zeta(\alpha)\ldots\zeta(\gamma)\ldots\zeta(\beta)\ldots
\end{array}
$$
Since $u\lessdot_k (\alpha,\beta)u$ is a cover in the $k$-Bruhat order,
either $k<u^{-1}(\beta)<u^{-1}(\gamma)$ or else
$u^{-1}(\gamma)<u^{-1}(\alpha) \leq k$.
We illustrate $u$, $(\alpha,\beta)u$, and $w$ for each possibility:
$$
\begin{array}{rccc}
&k<u^{-1}(\beta)<u^{-1}(\gamma)&&u^{-1}(\gamma)<u^{-1}(\alpha)\leq k\\
u:\ &
\ldots\ \ \alpha\ \ \ldots\ \ \beta\ \ \ldots\ \ \gamma\ \ \ldots
&\qquad & 
\ldots\ \ \gamma\ \ \ldots\ \ \alpha\ \ \ldots\ \ \beta\ \ \ldots\\
(\alpha,\beta)u:\ &
\ldots\ \ \beta\ \ \ldots\ \ \alpha\ \ \ldots\ \ \gamma\ \ \ldots&&
\ldots\ \ \gamma\ \ \ldots\ \ \beta\ \ \ldots\ \ \alpha\ \ \ldots\\
w:\ &
\ldots\zeta(\alpha)\ldots\zeta(\beta)\ldots\zeta(\gamma)\ldots&&
\ldots\zeta(\gamma)\ldots\zeta(\alpha)\ldots\zeta(\beta)\ldots
\end{array}
$$
Assume  $k<u^{-1}(\beta)<u^{-1}(\gamma)$ .
Then Theorem~\ref{thm:k-length} and  
$(\alpha,\beta)u\leq_k w$ imply  $\gamma \geq \zeta(\gamma)$
and $\zeta(\beta)<\zeta(\gamma)$,
since $\alpha<\gamma$ and both have positions
greater than $k$ in $(\alpha\,\beta)u$.
Let $c:=x^{-1}(\gamma)$.
If $c\leq k$, then $x\leq_k z$ implies 
$\gamma \leq \zeta(\gamma)$ so $\gamma=\zeta(\gamma)$.
Also, $\alpha<\gamma$ implies $\zeta(\alpha)<\zeta(\gamma)$
and thus $\zeta(\gamma)=\gamma <\beta\leq \zeta(\alpha)$, a contradiction.
Similarly, if $c>k$, then $\gamma<\beta$ implies 
$\zeta(\gamma)<\zeta(\beta)$, another contradiction.
The other possibility, $u^{-1}(\gamma)<u^{-1}(\alpha)$,  leads to a
similar contradiction. 
Thus $x\lessdot_k (\alpha,\beta)x$ is a cover in the $k$-Bruhat order.

To show $y:= (\alpha,\beta)x\leq_k z$, first note that the pair 
$(y,z)$ satisfy condition I of
Theorem~\ref{thm:k-length}, because 
$(\alpha,\beta)u\leq_k w$.
For condition II,  we need only show:
\begin{enumerate}
\item[a)] If $\alpha<\gamma<\beta$ and $x^{-1}(\gamma)<x^{-1}(\alpha)$, 
so that $\gamma=yx^{-1}(\gamma)<yx^{-1}(\alpha)=\beta$, 
then $zx^{-1}(\gamma) = \zeta(\gamma) <\zeta(\beta) = zx^{-1}(\alpha)$,
and
\item[b)] If $\alpha<\gamma<\beta$ and $x^{-1}(\beta)<x^{-1}(\gamma)$, 
so that $\alpha=yx^{-1}(\beta)<yx^{-1}(\gamma)=\gamma$, 
then $\zeta(\alpha) <\zeta(\gamma)$.
\end{enumerate}

If $\alpha<\gamma<\beta$, then one of these two possibilities does occur,
as $x\lessdot_k (\alpha,\beta) x = y$ is a cover in the $k$-Bruhat order.
Suppose $x^{-1}(\gamma)<x^{-1}(\alpha)$, as the other case is similar.

Since $x^{-1}(\gamma)<k$ and $x\leq_k z$, we have $\gamma\leq
\zeta(\gamma)$, by condition I.
If $u^{-1}(\gamma)<u^{-1}(\alpha)$, then 
$(\alpha\,\beta)u\leq_k w\Rightarrow \zeta(\gamma)<\zeta(\alpha)$.
If $u^{-1}(\beta)<u^{-1}(\gamma)$, then $\gamma=\zeta(\gamma)$, and so
$\zeta(\gamma)=\gamma<\beta\leq\zeta(\alpha)$. 
Since $u\leq_k (\alpha\,\beta)u$, we cannot have
$u^{-1}(\alpha)<u^{-1}(\gamma)<u^{-1}(\beta)$.
\QED

Define $\mathrm{up}_\zeta := \{ \alpha\,|\, \alpha<\zeta(\alpha)\}$ and 
$\mathrm{down}_\zeta := \{ \beta\,|\, \beta>\zeta(\beta)\}$.

\begin{thm}\label{thm:lengthenings_exist}
Let $\zeta\in {\cal S}_\infty$.
\begin{enumerate}
\item[({\em i})]  For $u\in {\cal S}_\infty$, $u\leq_k \zeta u$ if
and only if the 
following conditions are satisfied.
\begin{enumerate}
\item $u^{-1} \mathrm{up}_\zeta \subset \{1,\ldots,k\}$,
\item $u^{-1} \mathrm{down}_\zeta \subset \{k+1,k+2,\ldots\}$, and 
\item For all $\alpha,\beta\in \mathrm{up}_\zeta$ 
(respectively $\alpha,\beta\in  \mathrm{down}_\zeta$), 
$\alpha<\beta$ and $u^{-1}(\alpha) < u^{-1}(\beta)$ together imply
$\zeta(\alpha)<\zeta(\beta)$. 
\end{enumerate}
\item[({\em ii})] If $\# \mathrm{up}_\zeta \leq k$, then there is a
permutation $u$ such that $u\leq_k \zeta u$.
\end{enumerate}
\end{thm}

\noindent{\bf Proof. }
Statement ({\em i}) is a consequence of 
Theorem~\ref{thm:k-length}.
For ({\em ii}), let $\{a_1,\ldots,a_k\}\subset{\Bbb N}$ contain
$\mathrm{up}_\zeta$ and possibly some fixed points of $\zeta$,
and let $\{a_{k+1},a_{k+2},\ldots\}$ be the complementary set in 
${\Bbb N}$.
Suppose these sets are indexed so that
$\zeta(a_i)<\zeta(a_{i+1})$ for $i\neq k$.
Define $u\in {\cal S}_\infty$ by $u(i)   = a_i$
Then $\zeta u$ is Grassmannian with descent $k$, and
Theorem~\ref{thm:k-length} implies $u\leq_k\zeta u$.
\QEDnoskip

\subsection{A new partial order on ${\cal S}_\infty$}\label{sec:new:order}

For $\zeta\in {\cal S}_\infty$, define $|\zeta|$ to be the quantity:
$$
\begin{array}{c}
\#\{(\alpha,\beta)\in 
\zeta(\mathrm{up}_\zeta)\times\zeta(\mathrm{down}_\zeta)\,|\,
\alpha>\beta\}\ -\ 
 \#\{a,b\in \mathrm{up}_\zeta\,|\,
a>b \mbox{ and  } \zeta(a)<\zeta(b)\}\\ 
 -\  
\#\{a,b\in \mathrm{down}\zeta\,|\,
a>b \mbox{ and } \zeta(a)<\zeta(b)\}\ -\ 
\#\{(a,b)\in \mathrm{up}_\zeta\times\mathrm{down}_\zeta\,|\,
a>b\}.
\end{array}
$$

\begin{lem} % \label{lemma:new_length}
If $u\leq_k \zeta u$, then 
$\ell(u) + |\zeta|= \ell(\zeta u)$.
\end{lem}

\noindent{\bf Proof. } 
By Theorem~\ref{thm:k-order}, 
$\ell(\zeta u)-\ell(u)$ depends only upon $\zeta$.
Computing this for the permutation $u$ defined in the proof of
Theorem~\ref{thm:lengthenings_exist},  shows it equals $|\zeta|$:
If $c=\zeta(c)$, then the number of inversions involving $c$ in
$u$ equals the number involving $c$ in $\zeta u$.
The first term of the expression for $|\zeta|$ counts the remaining
inversions in $\zeta u$ and the last three terms the remaining
inversions in $u$. 
\QED

By Theorem~\ref{thm:k-order}, the interval
$[u,\zeta u]_k$ depends only upon $\zeta$ if $u\leq_k \zeta u$.
A closer examination of our arguments shows it is independent of
$k$, as well. 
That is, if $x\leq_l \zeta x$,
then the map $v \mapsto xu^{-1}v$ defines an isomorphism
$[u,\,\zeta u]_k\stackrel{\sim}{\longrightarrow}[x,\,\zeta x]_l$.
This motivates the following definition.

\begin{defn}
{\em 
For $\zeta,\eta\in{\cal S}_\infty$, let
$\eta\preceq \zeta$ if 
there exists $u\in {\cal S}_\infty$ and a positive integer $k$ such that 
$u\leq_k\eta u\leq_k\zeta u$.
If $u$ is chosen as in the proof of
Theorem~\ref{thm:lengthenings_exist}, then we see that
$\eta\preceq\zeta$ if 
\begin{enumerate}
\item if $\alpha<\eta(\alpha)$, then $\eta(\alpha)\leq\zeta(\alpha)$,
\item if $\alpha>\eta(\alpha)$, then $\eta(\alpha)\geq\zeta(\alpha)$, and 
\item if
$\alpha,\beta\in\mathrm{up}_\zeta$ (respectively, 
$\alpha,\beta\in\mathrm{down}_\zeta$) with $\alpha<\beta$ and
$\zeta(\alpha)<\zeta(\beta)$, then $\eta(\alpha)<\eta(\beta)$.
\end{enumerate}
}
\end{defn}

\begin{figure}[htb]
$$\epsfxsize=4.in \epsfbox{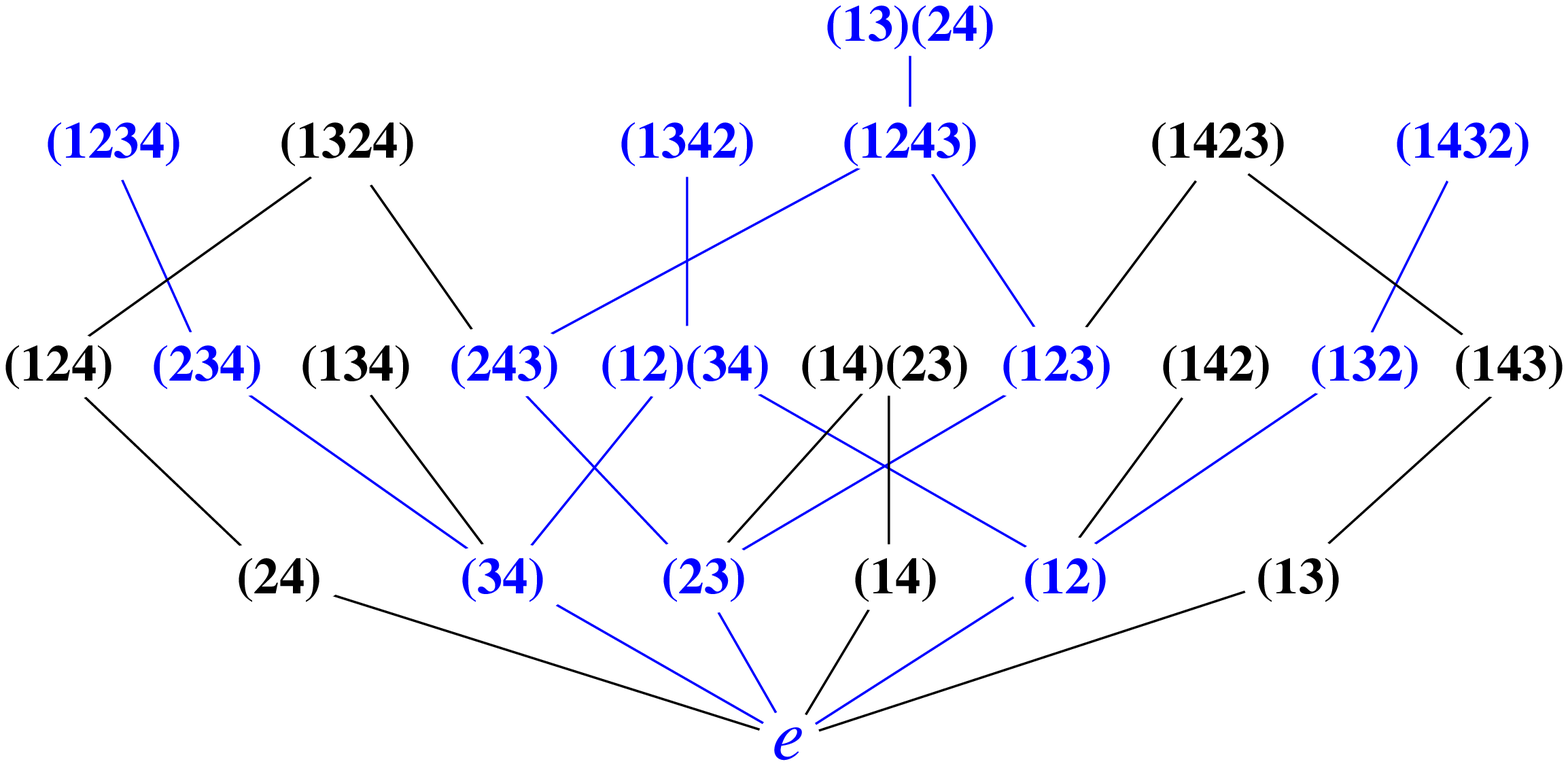}$$
\caption{ $\preceq$ on  ${\cal S}_4$\label{fig:new_order.S_4}}
\end{figure}

Figure~\ref{fig:new_order.S_4} illustrates $\preceq$ on ${\cal S}_4$.
%We  illustrate $\preceq$ on ${\cal S}_4$:

\begin{thm}\label{thm:new_order}
Suppose $u,\zeta,\eta,\xi\in{\cal S}_\infty$.
\begin{enumerate}
\item[({\em i})] $({\cal S}_\infty,\preceq)$ is a graded poset 
with rank function $|\zeta|$.
\item[({\em ii})] The map $\lambda \mapsto v(\lambda,k)$ exhibits Young's
lattice of partitions with at most $k$ parts as an induced suborder of  
$({\cal S}_\infty,\preceq)$.
\item[({\em iii})] If $u\leq_k \zeta u$,
then the map $\eta\mapsto \eta u$ induces an isomorphism
$[e,\,\zeta]_\preceq \stackrel{\sim}{\longrightarrow}
[u,\,\zeta u]_k$.
\item[({\em iv})] 
If $\eta\preceq\zeta$, then the map 
$\xi\mapsto \xi\eta^{-1}$ induces
an isomorphism $[\eta,\zeta]_\preceq \stackrel{\sim}{\longrightarrow}
[e,\zeta\eta^{-1}]_\preceq$.
\item[({\em v})] 
For every infinite set  $P\subset {\Bbb N}$,
$\phi_P:{\cal S}_\infty \rightarrow {\cal S}_\infty$ is an injection of
graded posets.
Thus, if $\zeta, \eta\in{\cal S}_\infty$ are shape equivalent, then 
$[e,\,\zeta]_\preceq \simeq [e,\,\eta]_\preceq$.
\item[({\em vi})] 
The map $\eta\mapsto \eta\zeta^{-1}$ induces an order reversing isomorphism
between $[e,\zeta]_\preceq$ and $[e,\zeta^{-1}]_\preceq$.
\item[({\em vii})] 
The homomorphism $\zeta\mapsto \overline{\zeta}$ on ${\cal S}_n$ induces an
automorphism of $({\cal S}_n,\preceq)$.
\end{enumerate}
\end{thm}

Theorem~\ref{thm:B} ({\em i}) is an immediate consequence of the
definition of $\preceq$ and ({\em v}).

\noindent{\bf Proof. }
Statements ({\em i})--({\em v}) follow from the definitions.
Suppose $u\leq_k \eta u\leq_k \zeta u$ with  
$u,\eta u,\zeta u\in{\cal S}_n$.
If $w:= \zeta u$, then $ w w_0 \leq_{n-k} \eta\zeta^{-1}w w_0
\leq_{n-k} \zeta^{-1}w w_0$, which proves ({\em vi}).  
Similarly, $u\leq_k w \Leftrightarrow  
\overline{u}\leq_{n-k}\overline{w}$  implies ({\em vii}).
\QEDnoskip

\begin{ex}\label{ex:neworder}{\em
Let $\zeta = (24)(153)$ and $\eta=(35)(174)=\phi_{\{1,3,4,5,7\}}(\zeta)$.
Then
$21345 \leq_2 45123 = \zeta\cdot 21345$
and
$3215764 \leq_3 5273461 = \eta\cdot 3215764$.
Figure~\ref{fig:new_order} illustrates the intervals	
$[21342,\,45123]_2$, $[3215764,\,5273461]_3$, and 
$[e,\zeta]_\preceq$.

\begin{figure}[htb]
$$\epsfxsize=5.3in \epsfbox{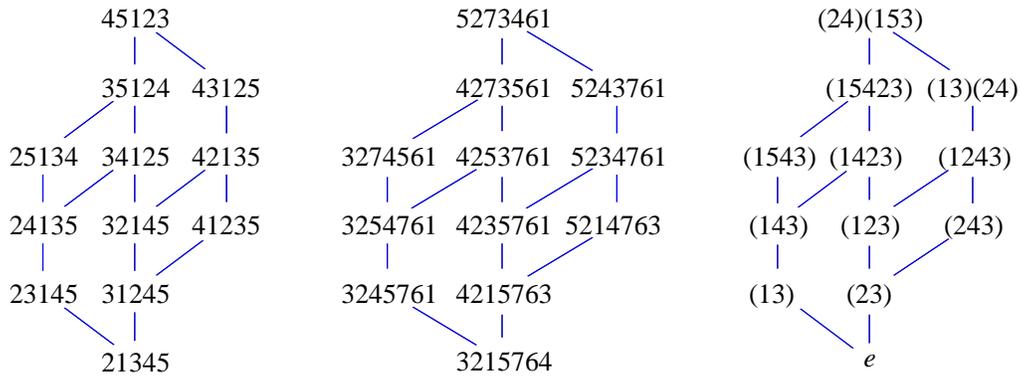}$$
\caption{Isomorphic intervals in $\leq_2, \leq_3$, and
$\preceq$\label{fig:new_order}} 
\end{figure}
}\end{ex}

\subsection{Disjoint permutations}\label{sec:disjoint_permutations}

Let $\zeta\in {\cal S}_n$ and $1,\ldots,n$ be the vertices of a convex
planar $n$-gon numbered consecutively.
Define $\Gamma_\zeta$, a directed geometric graph 
to be the union of directed chords
$\Span{\alpha,\zeta(\alpha)}$ for $\alpha$ in the support,
$\mbox{supp}_\zeta$, of $\zeta$.

Permutations $\zeta$ and $\eta$ are {\em disjoint} if the edge sets of
their geometric graphs $\Gamma_\zeta$ and $\Gamma_\eta$ (drawn on the
same $n$-gon) are disjoint as {\em subsets of the plane}.
This implies (but is not equivalent to) 
supp$_\zeta\bigcap\:$supp$_\eta=\emptyset$.
Disjointness may be rephrased in terms of
partitions~\cite{Stanley_enumerative}  of $[n]$.
Suppose for simplicity, that supp$_\zeta$ and supp$_\eta$ partition
$[n]$.
Then $\zeta$ and $\eta$ are disjoint if and only if
there is a non-crossing partition~\cite{Kreweras} $\pi$ of $[n]$ 
refining the 
partition (supp$_\zeta$, supp$_\eta$), and which is itself refined by
the partition given by the cycles of $\zeta$ and $\eta$.

We compare the  graphs of the pairs of permutations 
$(1782), (345)$ and $(13), (24)$.
$$
\epsfxsize=2.5in \epsfbox{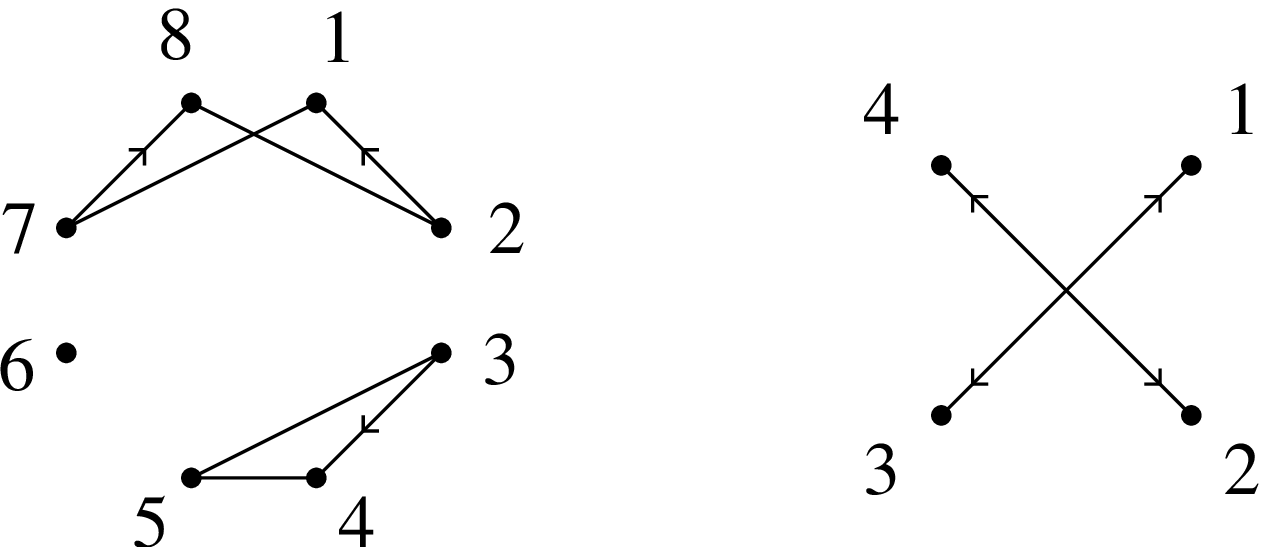}
$$
The first pair is disjoint and the second is not.
We relate this definition to that given in \S\ref{sec:Schur_identities}. 

\begin{lem}
Let $\zeta,\eta\in {\cal S_\infty}$.
Then the edges of $\/\Gamma_\zeta$ are disjoint from the edges of
$\/\Gamma_\eta$ if and only if $\/\zeta$ and $\eta$ have disjoint support
and $|\zeta|+|\eta|=|\zeta\eta|$. 
\end{lem}

\noindent{\bf Proof. }
Suppose $\zeta$ and $\eta$ have disjoint support and 
let  $\Span{a,\zeta(a)}$  be an edge of $\Gamma_\zeta$ and 
$\Span{b,\eta(b)}$ be an edge of $\Gamma_\eta$.
Consider the contribution of the endpoints of these edges to 
$|\zeta\eta|-|\zeta|-|\eta|$. 
This contribution is zero if the edges do not cross, which proves the
forward implication.

For the reverse, suppose these edges cross.
Then the contribution is 1 if $a<\zeta(a)$ and $b>\eta(b)$ (or vice-versa),
and 0 otherwise.
Because each edge is part of a directed cycle in the graph,
if one edge of $\Gamma_\zeta$ crosses an edge of $\Gamma_\eta$, then
there are at least four crossings, one of each type illustrated.
$$
\begin{picture}(100,80)\setlength{\unitlength}{.7pt}\thicklines
\put(25,25){\line(1,1){50}}
\put(37.5,37.5){\line(1,0){5}}\put(37.5,37.5){\line(0,1){5}}
\put(25,75){\line(1,-1){50}}
\put(35,65){\line(1,0){5}}\put(35,65){\line(0,-1){5}}
\put(-5,10){$\zeta(a)$}\put(-5,80){$\eta(b)$}
\put(80,10){$b$}\put(80,80){$a$}
\put(25,25){\circle*{5}}\put(25,75){\circle*{5}}
\put(75,25){\circle*{5}}\put(75,75){\circle*{5}}
\end{picture}%\hspace{-5pt}
\begin{picture}(100,85)\setlength{\unitlength}{.75pt}\thicklines
\put(25,25){\line(1,1){50}}
\put(35,35){\line(1,0){5}}\put(35,35){\line(0,1){5}}
\put(25,75){\line(1,-1){50}}
\put(65,35){\line(-1,0){5}}\put(65,35){\line(0,1){5}}
\put(-5,10){$\zeta(a)$}\put(12,80){$b$}
\put(78,10){$\eta(b)$}\put(80,80){$a$}
\put(25,25){\circle*{5}}\put(25,75){\circle*{5}}
\put(75,25){\circle*{5}}\put(75,75){\circle*{5}}
\end{picture}
\begin{picture}(100,85)\setlength{\unitlength}{.75pt}\thicklines
\put(25,25){\line(1,1){50}}
\put(65,65){\line(-1,0){5}}\put(65,65){\line(0,-1){5}}
\put(25,75){\line(1,-1){50}}
\put(35,65){\line(1,0){5}}\put(35,65){\line(0,-1){5}}
\put(12,10){$a$}\put(-5,80){$\eta(b)$}
\put(80,10){$b$}\put(80,80){$\zeta(a)$}
\put(25,25){\circle*{5}}\put(25,75){\circle*{5}}
\put(75,25){\circle*{5}}\put(75,75){\circle*{5}}
\end{picture}
\begin{picture}(100,85)\setlength{\unitlength}{.75pt}\thicklines
\put(25,25){\line(1,1){50}}
\put(65,65){\line(-1,0){5}}\put(65,65){\line(0,-1){5}}
\put(25,75){\line(1,-1){50}}
\put(65,35){\line(-1,0){5}}\put(65,35){\line(0,1){5}}
\put(12,10){$a$}\put(12,80){$b$}
\put(78,10){$\eta(b)$}\put(80,80){$\zeta(a)$}
\put(25,25){\circle*{5}}\put(25,75){\circle*{5}}
\put(75,25){\circle*{5}}\put(75,75){\circle*{5}}
\end{picture}
$$
Here, the numbers increase in a clockwise direction, with the least
number in the northeast ($\nearrow$).
Thus $|\zeta\eta|>|\zeta|+|\eta|$.
\QEDnoskip

\begin{lem}\label{lem:disjoint_one}
Let $\alpha<\beta$ and $\zeta\in{\cal S}_\infty$ and suppose
$\zeta\precdot(\alpha,\beta)\zeta$ is a cover.
Then
\begin{enumerate}
\item[({\em i})] $\alpha$ and $\beta$ are connected in the geometric
graph $\Gamma_{(\alpha,\beta)\zeta}$.
\item[({\em ii})] If $\Span{c,d}$ is any chord of the $n$-gon  meeting 
$\Gamma_\zeta$, then $\Span{c,d}$ meets $\Gamma_{(\alpha,\beta)\zeta}$.
\item[({\em iii})] If $p$ and $q$ are connnected in $\Gamma_\zeta$, then
they are 
connected in $\Gamma_{(\alpha,\beta)\zeta}$.
\item[({\em iv})] If $\zeta$ and $\eta$ are disjoint permutations and
$\zeta'\preceq\zeta$, then $\zeta'$ and $\eta$ are disjoint.
\end{enumerate}
\end{lem}

\noindent{\bf Proof. }
Suppose $u\in {\cal S}_\infty$ with 
$u\leq_k \zeta u\lessdot_k (\alpha,\beta)\zeta u$.
Define $i$ and $j$ by $\zeta u(i)=\alpha$ and $\zeta u(j)=\beta$, and
set $a=u(i)$ and $b=u(j)$.
Since $\zeta u\lessdot_k (\alpha,\beta)\zeta u$ is a cover, $i\leq k <j$, 
and thus $a\leq \alpha<\beta\leq b$, since 
$u\leq_k \zeta u$.
Thus the edges $\Span{a,\beta}$ and $\Span{b,\alpha}$
of $\Gamma_{(\alpha,\beta)\zeta}$ meet, proving ({\em i}).

For ({\em ii}), note that $\Gamma_{(\alpha,\beta)\zeta}$ differs from
$\Gamma_\zeta$ only by the (possible) deletion of edges
$\Span{a,\alpha}$ and $\Span{b,\beta}$ and the addition of the edges
$\Span{a,\beta}$ and $\Span{b,\alpha}$.
Checking all possibilities for the chords $\Span{c,d}$,
$\Span{a,\alpha}$, and $\Span{b,\beta}$ shows ({\em ii}).

Statement ({\em iii}) follows from ({\em ii}) by considering edges of 
$\Gamma_\zeta-\Span{a,\alpha}-\Span{b,\beta}$.

The contrapositive of ({\em iv}) is also a consequence of ({\em ii}); 
If $\zeta'$ and $\eta$ are not disjoint and $\zeta'\preceq\zeta$, then 
$\zeta$ and $\eta$ are not disjoint.
\QEDnoskip

\begin{lem}\label{lem:disjoint_two}
Suppose $\zeta$ and $\eta$ are disjoint permutations.
For every $u\in {\cal S}_\infty$,
$$
u\ \leq_k\ \zeta \eta u \quad \Longleftrightarrow\quad
u\ \leq_k\ \zeta u \quad\mbox{and}\quad u\ \leq_k\ \eta u.
$$
\end{lem}

\noindent{\bf Proof. }
Suppose $u \leq_k \zeta \eta u$.
Let $i\leq k$ so that $u(i)\leq \zeta\eta u(i)$.
Since $\zeta$ and $\eta$ have disjoint supports, 
$u(i)\leq \zeta u(i)$.
Similarly, if $k<j$, then $u(j)\geq \zeta u(j)$,
showing Condition I of Theorem~\ref{thm:k-length} holds for the pair 
$(u,\zeta u)$.

For Condition II, suppose 
$i<j$, with $u(i)<u(j)$ and $\zeta u(i)>\zeta u(j)$.
If $j\leq k$, this implies $u(i)\in \mbox{supp}_\zeta$.
Since $u \leq_k \zeta \eta u$, and $\zeta,\eta$ have disjoint supports,
we have $\zeta u(i) = \eta\zeta u(i)<\eta\zeta u(j)$,
which implies $u(j)\in\mbox{supp}_\eta$ and so
$$
u(i) \ <\ u(j)\ <\ \zeta u(i) \ <\ \eta u(j).
$$
But then the edge $\Span{u(i),\,\zeta u(i)}$ of $\Gamma_\zeta$ meets the
edge  $\Span{u(j),\,\eta u(j)}$ of $\Gamma_\eta$, a contradiction.
The assumption that $k<i$ leads similarly to a contradiction.
Thus $u\leq_k \zeta u$ and similarly, $u\leq_k \eta u$.

Suppose now that  $u\leq_k \zeta u$ and $u\leq_k \eta u$.
Condition I of Theorem~\ref{thm:k-length} for $(u,\zeta\eta u)$
holds as $\zeta$ and $\eta$
have disjoint support. 
For condition II, let $i<j$ with $u(i)<u(j)$ and suppose that 
$j\leq k$.
If the set $\{u(i),u(j)\}$ meets at most one of supp$_\zeta$ or
supp$_\eta$ , say  supp$_\zeta$, then $u\leq_k \zeta u$
implies $\zeta\eta u(i) < \zeta \eta u(j)$.
Suppose now that $u(i)$ is in the support of $\zeta$ and $u(j)$ is in
the support of $\eta$.
Since $u\leq_k \zeta u$, we have $\zeta u(i)<\zeta u(j) =u(j)$.
But $u\leq_k \eta u$ implies $u(j)\leq \eta u(j)$.
Thus 
$\eta \zeta u(i) = \zeta u(i) < u(j)\leq \eta u(j) = \eta \zeta u(j)$.
Similar arguments suffice when $k<i$.
\QEDnoskip

\begin{thm}\label{thm:disjoint_iso}
Suppose $\zeta$ and $\eta$ are disjoint.
Then the map $[e,\zeta]_\preceq\times [e,\eta]_\preceq \rightarrow 
[e,\zeta\eta]_\preceq$
defined by $(\zeta',\eta') \mapsto \zeta'\eta'$
is an isomorphism of graded posets.
\end{thm}

\noindent{\bf Proof. }
By Lemmas~\ref{lem:disjoint_one} and~\ref{lem:disjoint_two}, 
this map is an injection of graded posets.
For surjectivity, let $\xi\preceq\zeta\eta$.
By Lemma~\ref{lem:disjoint_one} ({\em iii}) and downward induction from
$\zeta\eta$ to $\xi$, $\Gamma_\xi$ has no edges connecting 
supp$_\zeta$ to supp$_\eta$.
Set  $\xi':= \xi|_{\mbox{\scriptsize supp}_\zeta}$, and $\xi'':=
\xi|_{\mbox{\scriptsize supp}_\eta}$.
Then $\xi = \xi' \xi''$ and $\xi'$ and $\xi''$ are
disjoint. 
Surjectivity will follow by showing $\xi'\preceq\zeta$ and
$\xi''\preceq\eta$. 

It suffices to consider the case 
$\xi\precdot (\alpha,\beta)\xi = \zeta\eta$ is a cover. 
By Lemma~\ref{lem:disjoint_one} ({\em i}), $\alpha$ and $\beta$ are
connected 
in $\Gamma_{\zeta\eta}$, so we may assume that
$\alpha,\beta$ are both in the support of $\zeta$.
Then $\xi''=\eta$ and $(\alpha,\beta)\xi'=\zeta$.
We show that $\xi' \preceq (\alpha,\beta)\xi'=\zeta$ is a cover,
which will complete the proof.

Choose $u\in {\cal S}_\infty$ with 
$u\leq_k\xi u\leq_k\zeta\eta u$.
Let $a:=(\xi'u)^{-1}(\alpha)$ and $b:=(\xi' u)^{-1}(\beta)$.
Since $\xi'$ and $\eta$ are disjoint, $\alpha,\beta\notin\mbox{supp}_\eta$
and so $a,b\notin\mbox{supp}_\eta$.
Thus $(\alpha,\beta)\xi'\eta u = \xi'\eta u(a,b)$, showing
$a\leq k <b$, as $\xi'\eta u\lessdot_k (\alpha,\beta)\xi'\eta u$.

Since $\xi'$ and $\eta$ are disjoint and $\xi=\xi'\eta$,
Lemma~\ref{lem:disjoint_two} implies $u\leq_k \xi' u$.
Thus $|\xi'|+\ell(u)=\ell(\xi' u)$.
But since $\xi'$ and $\eta$ are disjoint and 
$\xi'\eta \precdot \zeta\eta$ is a cover, we have
$$
|\zeta|+|\eta|\  =\ |\zeta\eta|\ =\  
1+ |\xi'\eta|\ =\ 1 + |\xi'| + |\eta|,
$$
so $\ell(\xi'u)+1 = \ell(\xi' u(a,b))$.
Since $a\leq k<b$ and $\zeta u=\xi'u(a,b)$, this implies
$\xi' u\lessdot_k \zeta u$.
\QEDnoskip

\begin{ex}
{\em 
Let $\zeta=(2354)$ and $\eta=(176)$, which are disjoint.
Let $u= 2316745$.
Then 
$$
u\leq_3 \zeta\eta u = 3571624, \quad
u\leq_3 \zeta u = 3516724, \quad \mbox{ and }
u\leq_3 \eta u = 2371645.
$$
The intervals $[u,\zeta u]_3$,  $[u,\eta u]_3$, and 
$[u,\;\zeta\eta u]_3$ are illustrated in 
Figure~\ref{fig:disjoint_interval}. 
\begin{figure}[htb]
$$\epsfxsize=5.2in \epsfbox{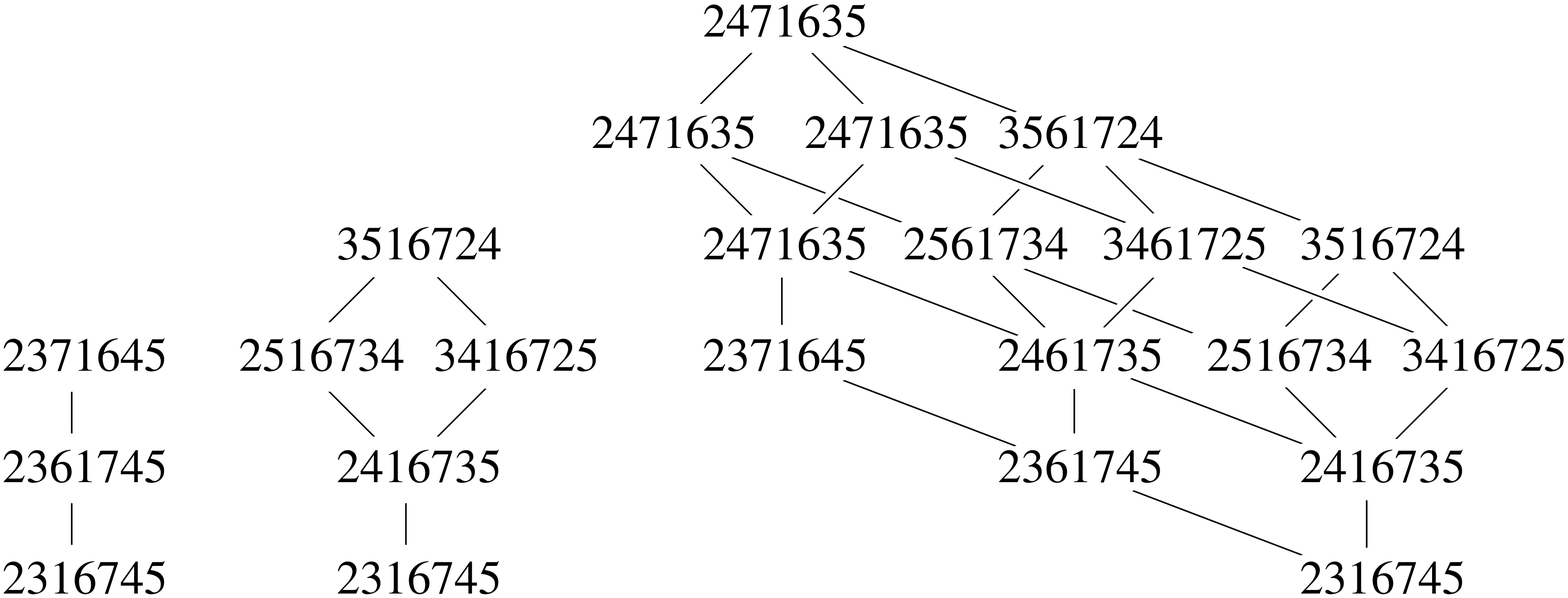}$$
\caption{Intervals of disjoint permutations\label{fig:disjoint_interval}}
\end{figure}

}\end{ex}

\section{Cohomological formulas and identities for the $c^w_{u\,v}$}
\subsection{Two maps on ${\cal S}_\infty$}\label{sec:fixed_points} 
For positive integers $p,q$  and $w\in {\cal S}_\infty$, define 
$\varepsilon_{p,q}(w)\in {\cal S}_\infty$
by
$$
\varepsilon_{p,q} (w)(j) \quad =\quad \left\{\begin{array}{lcl}
w(j)     && j < p \mbox{ and } w(j)  <  q\\
w(j)+1   && j < p \mbox{ and } w(j)\geq q\\
q        && j = p\\
w(j-1)   && j > p \mbox{ and } w(j)  <  q\\
w(j-1)+1 && j > p \mbox{ and } w(j)\geq q
\end{array}\right..
$$
Note that $\varepsilon_{p,p}=\phi_{{\Bbb N}-\{p\}}$.
However,  if $p\neq q$, then this injection, 
$\varepsilon_{p,q}: {\cal S}_\infty \hookrightarrow {\cal S}_\infty$,
is not a group homomorphism.
The map $\varepsilon_{p,q}$ has a left inverse 
$/_p : {\cal S}_\infty\rightarrow {\cal S}_\infty$:
For $x\in {\cal S}_\infty$, define $x/_p$ by
$$
x/_p(j)\quad =\quad  \left\{\begin{array}{lcl}
x(j)    &&j  <  p\mbox{ and }x(j) < x(p)\\
x(j)-1  &&j  <  p\mbox{ and }x(j) > x(p)\\
x(j+1)  &&j\geq p\mbox{ and }x(j) < x(p)\\
x(j+1)-1&&j\geq p\mbox{ and }x(j) > x(p)
\end{array}\right..
$$
Representing permutations as matrices, the 
effect of $/_p$ on $x$ is to erase the $p$th row and $x(p)$th column.
The effect of $\varepsilon_{p,q}$ is to expand the matrix by adding a
new $p$th row and $q$th column consisting mostly of zeroes, but with a 1
in the $(p,q)$th position.
For example,
$$
\varepsilon_{3,3}(23154)\ =\ 243165 \quad\mbox{and}\quad
264351/_3\ =\ 25341 
$$
These maps have
some order-theoretic properties.

\begin{lem}\label{lem:expanding_bruhat}
Suppose $x\leq z$ and $p,q$ are positive integers.
Then
\begin{enumerate}
\item[({\em i})] $\varepsilon_{p,q}(x)\leq \varepsilon_{p,q}(z)$.
\item[({\em ii})] If $\ell(z)-\ell(x) = 
\ell(\varepsilon_{p,q}(z))-\ell(\varepsilon_{p,q}(x))$,
then 
$$
\varepsilon_{p,q}\ :\ [x,z] \ \stackrel{\sim}{\longrightarrow}
[\varepsilon_{p,q}(x),\varepsilon_{p,q}(z)].
$$
\item[({\em iii})] If $x,z\in{\cal S}_n$ and either of $p$ or $q$ is
equal to either 1 or $n+1$, then  
$\ell(z)-\ell(x) = 
\ell(\varepsilon_{p,q}(z))-\ell(\varepsilon_{p,q}(x))$.

\item[({\em iv})]If $x\leq_k z$ and $x(p)=z(p)$, then 
$x/_p \leq_{k'} z/_p$ and $[x,z]_k\simeq[x/_p,z/_p]_{k'}$, 
where $k'$ is equal to $k$ if $k<p$ and $k-1$ otherwise. 
Furthermore,
$zx^{-1} = \varepsilon_{x(p),x(p)}(z/_p(x/_p)^{-1})$. 
\end{enumerate}
\end{lem}

\noindent{\bf Proof. }
Suppose $x\lessdot x(a,b)$ is a cover.
Then $\varepsilon_{p,q}(x)< \varepsilon_{p,q}(x(a,b))$ is a cover
if either $p\leq a$ or $b<p$, or else $a<p\leq b$ and either 
$q\leq x(a)$ or $x(b)<q$.
If however, $a<p\leq b$ and $x(a)<q\leq x(b)$, then there is a chain of
length 3 from $\varepsilon_{p,q}(x)$ to 
$\varepsilon_{p,q}(x(a,b))=\varepsilon_{p,q}(x) (a,b{+}1)$:
$$
\varepsilon_{p,q}(x)\  \lessdot \ 
\varepsilon_{p,q}(x)(a,p)\  \lessdot \ 
\varepsilon_{p,q}(x) (a,b{+}1,p) \  \lessdot \  
\varepsilon_{p,q}(x) (a,b{+}1).
$$

The lemma follows from this observation.
For example, under the hypothesis of ({\em ii}), the number of
inversions in $\varepsilon_{p,q}(z)$ involving $q$ equals the number of
inversion in  $\varepsilon_{p,q}(x)$ involving $q$.
Thus, if $\varepsilon_{p,q}(x)\leq u\leq \varepsilon_{p,q}(z)$,
then $u(p)=q$.
\QEDnoskip 

\subsection{An embedding of flag manifolds}\label{sec:emdedd}
Let $W\subset V$ with $W\simeq {\Bbb C\,}^n$ and 
$V\simeq {\Bbb C\,}^{n+1}$.
Suppose $f\in V- W$ so that $V=\Span{W,f}$.
For $p\in [n{+}1]$ define the injection 
$\psi_p: {\Bbb F}\ell W\hookrightarrow {\Bbb F}\ell V$ by
$$
\left( \psi_p \Edot\right)_j\quad =\quad \left\{ \begin{array}{lcl}
E_j &&\mbox{if } j<p\\
\Span{E_{j-1},f} &&\mbox{if }j\geq p \end{array}\right.
$$

\begin{prop}[\cite{sottile_pieri_schubert},
Lemma~12]\label{prop:embedding}  
Let $\Edot\in {\Bbb F}\ell W$ and $w\in {\cal S}_n$.
Then, for every $p,q\in [n{+}1]$,
$$
\psi_p X_w \Edot \quad \subset \quad
X_{\varepsilon_{p,q}(w)} \psi_{n+2-q}\Edot.
$$
\end{prop}

Recall that $e$ is the identity permutation.

\begin{cor}\label{cor:geometric_pushforward}
Let $w\in {\cal S}_n$ and $\Edot,\Epdot \in {\Bbb F}\ell W$ be opposite flags.
Then $\psi_1 \Edot$ and $\psi_{n+1}\Epdot$ are opposite flags in 
${\Bbb F}\ell V$ and 
$$
\psi_p X_w\Edot \ =\ 
X_{\varepsilon_{p,1}(w)}\psi_{n+1}\Edot \bigcap 
X_{\varepsilon_{p,n+1}(e )}\psi_1\Epdot
\ =\ 
X_{\varepsilon_{p,1}(e )}\psi_{n+1}\Epdot \bigcap 
X_{\varepsilon_{p,n+1}(w)}\psi_{1}\Edot.
$$
\end{cor}

\noindent{\bf Proof. }
Since $X_e \Epdot = {\Bbb F}\ell W$,
Proposition~\ref{prop:embedding} with $q=1$ or $n+1$
implies $\psi_p X_w\Edot$
is a subset of either intersection:
$$
X_{\varepsilon_{p,1}(w)}\psi_{n+1}\Edot \bigcap 
X_{\varepsilon_{p,n+1}(e )}\psi_1\Epdot
\qquad\mbox{or}\qquad
X_{\varepsilon_{p,1}(e )}\psi_{n+1}\Epdot \bigcap
X_{\varepsilon_{p,n+1}(w)}\psi_{1}\Edot.
$$
Since $\Edot$ and $\Epdot$ are opposite flags, 
$\psi_{n+1}\Edot$ and $\psi_1\Epdot$ are
opposite flags, so both  
intersections are generically transverse and irreducible.
Since 
$$
\ell(\varepsilon_{p,1}(w))\ =\ \ell(w) + p-1\ \qquad \mbox{and}\ \qquad 
\ell(\varepsilon_{p,n+1}(w))\ =\ \ell(w) + n+1-p,
$$
both intersections have the same dimension 
as $\psi_p X_w\Edot$, proving equality.
\QED

Since $\varepsilon_{p,n+1}(e)=v(n{+}1{-}p,\,p)$, where $n+1-p$ is the
partition of $n+1-p$ into a single part, we see that 
${\frak S}_{\varepsilon_{p,n+1}(e )}=h_{n+1-p}(x_1,\ldots,x_p)$, 
the complete symmetric polynomial of degree $n+1-p$ in 
$x_1,\ldots,x_p$.
Similarly, 
${\frak S}_{\varepsilon_{p,1}(e )}=e_{p-1}(x_1,\ldots,x_{p-1}) =
x_1\cdots x_{p-1}$, 
as $\varepsilon_{p,1}=v(1^{p-1},p-1)$, where $1^{p-1}$ is the partition
of $p{-}1$ into $p{-}1$ equal parts, each of size 1.

\begin{cor}\label{cor:equal_products}
Let $w\in {\cal S}_n$. 
In $H^*{\Bbb F}\ell V$, 
$$
{\frak S}_{\varepsilon_{p,1}(w)} 
\cdot h_{n+1-p}(x_1,\ldots,x_p)\quad =\quad 
{\frak S}_{\varepsilon_{p,n+1}(w)} \cdot x_1\cdots x_{p-1}
$$
and these products are equal to $(\psi_p)_* {\frak S}_w$. 
\end{cor}

We use this to compute $\psi_p^*$.
The Pieri-type formulas of~\cite{sottile_pieri_schubert}
show that 
if $u\in {\cal S}_n$ and $k,m\leq n$ positive integers, then
\begin{eqnarray*}
\hspace{.8in}
{\frak S}_u \cdot{\frak S}_{w_0w} \cdot e_m(x_1\cdots x_k) &=& 
\left\{\begin{array}{ll}1&u\stackrel{c_{k,m}}{\llra} w\\ 
0&\mbox{otherwise}
\end{array}\right. %\label{eq:pieri_elementary}
\makebox[.1in]{\qquad\qquad} \hspace{1.35in}(4.2.1)\\
{\frak S}_u\cdot{\frak S}_{w_0w}  \cdot h_{n+1-m}(x_1,\ldots,x_k) &=&
\left\{\begin{array}{ll}1&u\stackrel{r_{k,m}}{\llra} w\\ 0&\mbox{otherwise}
\end{array}\right.,    \hspace{1.38in}(4.2.2)  %\label{eq:pieri_homogeneous}
\end{eqnarray*}
where $u\stackrel{c_{k,m}}{\llra} w$ if there is a (saturated)
chain in the $k$-Bruhat order
from $u$ to $w$:
$$
u\ \lessdot_k\ (\alpha_1,\beta_1)u\ \lessdot_k\ \cdots\ \lessdot_k \ 
(\alpha_m,\beta_m)\cdots(\alpha_1,\beta_1)u\ =\ w
$$
such that $\beta_1>\cdots>\beta_m$.
When $k=m$, it follows that
$\{\alpha_1,\ldots,\alpha_k\}=\{u(1),\ldots,u(k)\}$.
When $k=m=p-1$, we write $\cpp$ for this relation. 
Similarly, $u\stackrel{r_{k,m}}{\llra} w$ if there is a chain in the
$k$-Bruhat order: 
$$
u\ \lessdot_k\ (\alpha_1,\beta_1)u\ \lessdot_k\ \cdots\ \lessdot_k\ 
(\alpha_{n+1-m},\beta_{n+1-m})\cdots(\alpha_1,\beta_1)u\ =\ w
$$
such that $\beta_1<\beta_2<\cdots<\beta_{n+1-m}$.

\begin{thm}\label{thm:projection} 
Let $x\in {\cal S}_{n+1}$. 
In $H^*{\Bbb F}\ell_n$,
\begin{enumerate}
\item[({\em i})]
${\displaystyle
\psi_p^*{\frak S}_x \quad=\quad
\sum_{\stackrel{\mbox{\scriptsize$w\in {\cal S}_n$}}%
{x\cpp\varepsilon_{p,1}(w)}} {\frak S}_w
\quad=\quad
\sum_{\stackrel{\mbox{\scriptsize$w\in {\cal S}_n$}}%
{x\rpp\varepsilon_{p,n+1}(w)}} {\frak S}_w}$.
\item[({\em ii})] ${\displaystyle \psi_p^*(x_i)\ =\ \left\{
\begin{array}{lll}x_i&& i<p\\
0&& i=p\\x_{i-1}&& i>p\end{array}\right.}$.
\end{enumerate}
\end{thm}

\noindent{\bf Proof. }
In $H^*{\Bbb F}\ell_n$,
$$
\psi_p^* {\frak S}_x \quad =\quad 
\sum_{w\in {\cal S}_n}\,
\deg({\frak S}_{w_0w}\cdot\psi_p^* {\frak S}_x)\,{\frak S}_w.
$$
By the projection formula~(2.3.1) %\ref{eq:projection}) 
and Corollary~\ref{cor:equal_products}, we have
$$
\deg({\frak S}_{w_0w}\cdot\psi_p^* {\frak S}_x)\quad=\quad
\deg( {\frak S}_x\cdot(\psi_{p})_*{\frak S}_{w_0w})\quad=\quad
\deg({\frak S}_x\cdot
{\frak S}_{\varepsilon_{p,n+1}(w_0w)}\cdot
x_1\cdots x_{p-1}).
$$
Note that $\varepsilon_{p,n+1}(w_0w)=w_0^{(n+1)}\varepsilon_{p,1}(w)$.
By~(4.2.1), 	%\ref{eq:pieri_elementary}),
the triple product 
$$
{\frak S}_x\cdot
{\frak S}_{\varepsilon_{p,n+1}(w_0w)}\cdot
x_1\cdots x_{p-1}
$$
is zero unless 
$x\cpp \varepsilon_{p,1}(w)$, and in this case
it equals ${\frak S}_{w_0^{(n+1)}}$.
This establishes the first equality of ({\em i}).
For the second, use the other 
formula for $(\psi_p)_*{\frak S}_w$ from
Corollary~\ref{cor:equal_products} and~(4.2.2). %\ref{eq:pieri_homogeneous}).

For ({\em ii}), let ${\cal F}_{\DOT}$ be the tautological flag on 
${\Bbb F}\ell_{n+1}$, ${\cal E}_{\DOT}$ the  tautological flag on 
${\Bbb F}\ell_n$, and 1 the trivial line bundle.
Then 
$$
\psi_p^*({\cal F}_i/{\cal F}_{i-1})\ =\ 
\left\{\begin{array}{lll}{\cal E}_i/{\cal E}_{i-1}&&\mbox{ if } i<p\\
1&&\mbox{ if } i=p\\
{\cal E}_{i-1}/{\cal E}_{i-2}&&\mbox{ if } i>p\end{array}\right.,
$$
Since $-x_i$ is the Chern class of both 
${\cal F}_i/{\cal F}_{i-1}$
and ${\cal E}_i/{\cal E}_{i-1}$, we are done.
\QEDnoskip

\subsection{The endomorphism $x_p\mapsto 0$}\label{sec:endomorphism}
For $p\in{\Bbb N}$ and $x\in {\cal S}_\infty$,
define
$$
A_p(x)\ :=\ 
\{u\in{\cal S}_\infty\,|\, 
x\cpp\varepsilon_{p,1}(u)\}.
$$
\begin{lem}\label{lem:index_sets}
If $x\in {\cal S}_n$ and $p\leq n$, then 
${\displaystyle 
A_p(x)\ =\ 
\{u\in{\cal S}_n\,|\, 
x\stackrel{r_{p,n+1-p}}{\llllra}\varepsilon_{p,n+1}(u)\}}$.
\end{lem}

\noindent{\bf Proof. }
If $x\in{\cal S}_n$, $p\leq n$, and 
$x\cpp w$,
then $w\in{\cal S}_{n+1}$, so $A_p(x)\subset{\cal S}_n$.
But then $A_p(x)$ and 
$\{u\in{\cal S}_n\,|\, 
x\stackrel{r_{p,n+1-p}}{\llllra}\varepsilon_{p,n+1}(u)\}$
index the two equal sums
in Theorem~\ref{thm:projection}({\em i}).
\QED

Let $\Psi_p:{\Bbb Z}[x_1,x_2,\ldots]\rightarrow{\Bbb Z}[x_1,x_2,\ldots]$
be defined by
$$
\Psi_p(x_i)\ =\ \left\{
\begin{array}{lll}x_i&&\mbox{ if } i<p\\
0&&\mbox{ if } i=p\\x_{i-1}&&\mbox{ if } i>p\end{array}\right..
$$
\begin{thm}\label{thm:theorem_A_iii}
For $x\in {\cal S}_\infty$,  and $p\in{\Bbb N}$, 
${\displaystyle \Psi_p{\frak S}_x = 
\sum_{u\in A_p(x)} {\frak S}_u}$.
\end{thm}

\noindent{\bf Proof. }
For $p\leq n+1$, the homomorphism
$\Psi_p$ induces the map 
$\psi_p^*: H^*{\Bbb F}\ell_{n+1} \rightarrow H^*{\Bbb F}\ell_n$, by
Theorem~\ref{thm:projection} ({\em ii}).
Choosing $n$ large enough completes the proof.
\QEDnoskip

\begin{cor}\label{cor:more_identities}
For $w, x, y\in {\cal S}_\infty$ and $p\in{\Bbb N}$, 
$$
\sum_{u\in A_p(x)} \quad \sum_{v\in A_p(y)}
c^w_{u\,v} 
\qquad = \qquad
\sum_{\stackrel{\mbox{\scriptsize $z$}}{w \in A_p(z)}}c^z_{x\, y}.
$$
\end{cor}

\noindent{\bf Proof. }
Apply $\Psi_p$ to the identity 
${\frak S}_x\cdot{\frak S}_y = \sum_zc^z_{x\,y}{\frak S}_z$
to obtain:
$$
\sum_{u\in A_p(x)} \quad \sum_{v\in A_p(y)}
{\frak S}_u \cdot{\frak S}_v
\quad=\quad
\sum_z c^z_{x\, y}
\sum_{w \in A_p(z)}c^z_{x\, y}\,
{\frak S}_w.
$$
Expanding the product ${\frak S}_u \cdot{\frak S}_v$ and equating the
coefficients of ${\frak S}_w$ proves the identity.
\QEDnoskip

\begin{ex}
{\em
We illustrate the effect of  $\Psi_3$ with an example.
Since
\begin{eqnarray*}
{\frak S}_{413652} &=&
x_1^4x_2x_4x_5 + x_1^3x_2^2x_4x_5 + x_1^3x_2x_4^2x_5 +\\
&\ &
x_1^4x_2x_3x_4 + x_1^4x_2x_3x_5 + x_1^4x_3x_4x_5 + 
x_1^3x_2^2x_3x_4 + x_1^3x_2^2x_3x_5 + x_1^3x_2x_3^2x_4 +\\
&\ & 
x_1^3x_2x_3^2x_5 + x_1^3x_2x_3x_4^2 + 
x_1^3x_3^2x_4x_5 + x_1^3x_3x_4^2x_5 + 
2\cdot x_1^3x_2x_3x_4x_5,
\end{eqnarray*}
we have
$$
\Psi_3({\frak S}_{413652}) \ =\ 
x_1^4x_2x_3x_4 + x_1^3x_2^2x_3x_4 + x_1^3x_2x_3^2x_4.
$$
However, 
\begin{eqnarray*}
{\frak S}_{52341} &=& x_1^4x_2x_3x_4\\
{\frak S}_{42531} &=& x_1^3x_2^2x_3x_4 + x_1^3x_2x_3^2x_4,
\end{eqnarray*}
which shows
$$
\Psi_3({\frak S}_{413652})\ =\ {\frak S}_{52341}+{\frak S}_{42531}.
$$
To see this agrees with 
Theorem~\ref{thm:theorem_A_iii}, compute the permutations
$w$ such that 
$x\stackrel{c_3}{\longrightarrow}w$:
$$
\epsfxsize=2.in \epsfbox{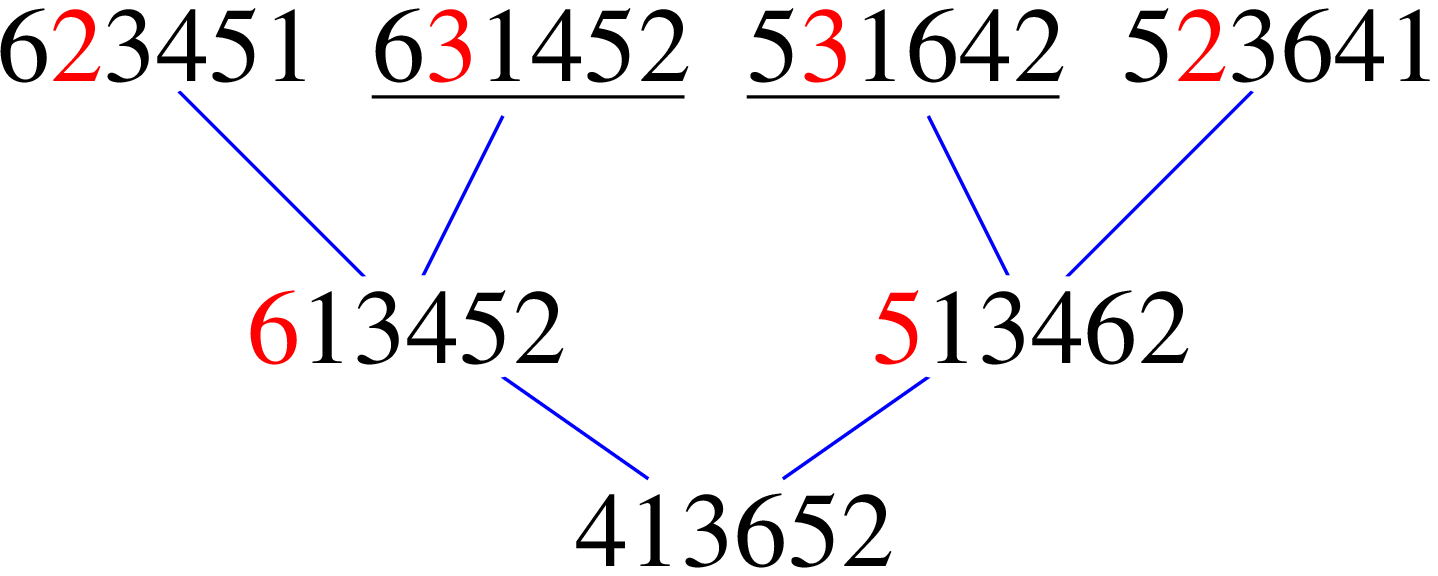}
$$
Of these, only the two underlined permutations are of the form
$\varepsilon_{3,1}(u)$: 
$$
631452\ =\ \varepsilon_{3,1}(52341)\quad\mbox{and}\quad
531642\ =\ \varepsilon_{3,1}(42531).
$$
}\end{ex}

\begin{lem}\label{lem:restriction}
Let $\lambda$ be a partition and $p,k$ positive integers.
Then $A_p(v(\lambda,k)) = \{(v(\lambda,k')\}$, where
$k'=k-1$ if $p\leq k$ and $k$ otherwise.
\end{lem}

\noindent{\bf Proof. }
By the combinatorial definition of Schur
functions~\cite[\S 4.4]{Sagan}, 
$\Psi_p({\frak S}_{v(\lambda,k)}) = {\frak S}_{v(\lambda,k')}$.
\QED

Lemma~\ref{lem:restriction} implies that $v(\lambda,k')$ is the only
solution $x$ to the 
equation $v(\lambda,k)\cpp\varepsilon_{p,1}(x)$, a statement about 
chains in the Bruhat order.

\subsection{Identities for $c^z_{x\,y}$ when
$x(p)=z(p)$}\label{sec:fixed_point_identities} 

\begin{lem}\label{lemma:fixed_pts}
Let $x,z\in {\cal S}_{n+1}$ with $x(p)=z(p)$ for some $p\in [m+1]$
and suppose $\ell(z)-\ell(x)= \ell(z/_p)-\ell(x/_p)$.
In $H^*{\Bbb F}\ell_{n+1}$, 
$$
(\psi_p)_* \left({\frak S}_{x/_p}\cdot
{\frak S}_{w_0^{(n)}(z/_p)}\right) 
\quad = \quad 
{\frak S}_x \cdot {\frak S}_{w_0^{(n+1)}z}.
$$
\end{lem}

\noindent{\bf Proof. }
Let $\Edot,\Epdot$ be opposite flags in $W$.
By Proposition~\ref{prop:embedding},
$$	%\begin{equation}\label{eq:push_forward}
	\psi_p\left(
	X_{w_0^{(n)}(z/_p)}\Edot \bigcap X_{x/_p} \Epdot \right)
	\quad = \quad
	X_{w_0^{(n+1)}z}\psi_{z(p)}\Edot \bigcap 
	X_x \psi_{n+2-x(p)}\Epdot
\eqno(4.4.1)
$$	%\end{equation}
Note that $w_0^{(n)}(z/_p)= (w_0^{(n+1)}z)/_p$.
Since $x(p)=z(p)$,
the flags $\psi_{z(p)}\Edot$ and $\psi_{n+2-x(p)}\Epdot$ are opposite  
in $V$.
Moreover, as $\ell(z)-\ell(x)= \ell(z/_p)-\ell(x/_p)$,
both sides of (4.4.1) %\ref{eq:push_forward}) 
have the same dimension,
so they are equal, proving the lemma.
\QEDnoskip

\begin{thm}\label{thm:coeff_sum}
Let $x,z\in {\cal S}_\infty$ with $x(p)=z(p)$ and suppose that  
$\ell(z)-\ell(x)= \ell(z/_p)-\ell(x/_p)$.
Then, for every $y\in {\cal S}_\infty$ and positive integer $p$,
$$
c^z_{x\, y}\quad =\quad 
\sum_{v\in A_p(y)} c^{z/_p}_{x/_p\: v}
$$
\end{thm}

\noindent{\bf Proof. }
It suffices to  compute this in 
$H^*{\Bbb F}\ell_{n+1}$, for $n$ such that 
$p\leq n$, $y\in {\cal S}_{n+1}$ and $A_p(y)\subset {\cal S}_n$.
By Lemma~\ref{lemma:fixed_pts}, 
$$
{\frak S}_x\cdot {\frak S}_{w_0^{(n+1)}z}
\ =\  (\psi_p)_* \left({\frak S}_{x/_p}\cdot 
{\frak S}_{w_0^{(n)}(z/_p)}\right)
\ =\  (\psi_p)_* \left(\sum_{v\in {\cal S}_n}
c^{w_0^{(n)}v}_{x/_p\ \,w_0^{(n)}(z/_p)}\:{\frak S}_{w_0^{(n)}v}\right).
$$
Since $c^{w_0^{(n)}v}_{u\:w_0^{(n)}w} = c^w_{u\,v}$ for 
$u,v,w\in {\cal S}_n$ and 
$\varepsilon_{p,1}(w_0^{(n)}v)=w_0^{(n+1)}\varepsilon_{p,1}(v)$,
\begin{eqnarray*}
{\frak S}_x\cdot {\frak S}_{w_0^{(n+1)}z}
&=&
 \sum_{v\in {\cal S}_n}
c^{z/_p}_{x/_p\:v} (\psi_p)_*\left({\frak S}_{w_0^{(n)}v}\right)\\
&=&
 \sum_{v\in {\cal S}_n}
c^{z/_p}_{x/_p\:v}
{\frak S}_{w_0^{(n+1)}\varepsilon_{p,1}(v)} \cdot x_1\cdots x_{p-1},
\end{eqnarray*}
by Corollary~\ref{cor:equal_products}.
Thus
\begin{eqnarray*}
c^z_{x\,y}&=& \deg\left({\frak S}_x\cdot {\frak S}_{w_0^{(n+1)}z}
\cdot{\frak S}_y\right)\\ 
&=& \sum_{v\in {\cal S}_n}c^{z/_p}_{x/_p\: v} \cdot
\deg\left({\frak S}_{w_0^{(n+1)}\varepsilon_{p,1}(v)} \cdot
(x_1\cdots x_{p-1})\cdot {\frak S}_y\right)\\
&=& 
\sum_{v\in A_p(y)}c^{z/_p}_{x/_p\: v}.\ \ \ \ \ \ \QEDnoskip
\end{eqnarray*}

When $p=1$, this has the following consequence:

\begin{cor}\label{cor:1_restrict}
If $x(1)=z(1)$, then $c^z_{x\,y}=0$ unless $y=1\times v$.
In that case,
$c^z_{x\;1\times v}=c^{z/\!_1}_{x/\!_1\: v}$.
\end{cor}

\subsection{Products of flag manifolds}
\label{sec:flag_products}
Let $P,Q\in {[n+m]\choose n}$, that is, $P,Q\subset [n+m]$ and each
has order $n$.
Index the sets $P,Q$ and their complements $P^c,Q^c$ as
follows: 
$$
\begin{array}{lcl}
P\ =\ p_1<\cdots<p_n &\qquad& P^c\ :=\ [n+m]-P\ =\ p^c_1<\cdots<p^c_m\\
Q\ =\ q_1<\cdots<q_n && Q^c\ :=\ [n+m]-Q\ =\ q^c_1<\cdots<q^c_m
\end{array}
$$
Define a function 
$\varepsilon_{P,Q}:
{\cal S}_n\times{\cal S}_m \hookrightarrow{\cal S}_{m+n}$
by:
\begin{eqnarray*}
\varepsilon_{P,Q}(v,w)(p_i) \ =\ q_{v(i)} && i=1,\ldots,n\\
\varepsilon_{P,Q}(v,w)(p^c_j) \ =\ q^c_{w(j)} && i=1,\ldots,m.
\end{eqnarray*}
As permutation matrices, 
$\varepsilon_{P,Q}(v,w)$ is obtained from $v$ and $w$ by placing the 
entries of $v$ in the blocks $P\times Q$ and those of $w$ in the
blocks $P^c\times Q^c$.
If $P=[n+1]- \{p\}$ and $Q=[n+1]-\{q\}$, then 
$\varepsilon_{P,Q}(v,e ) = \varepsilon_{p,q}(v)$.

Suppose $V\simeq {\Bbb C}\,^n$, $W\simeq {\Bbb C}\,^m$, and
$P\in{[n+m]\choose n}$.
Define a map
$$
\psi_P\ :\ {\Bbb F}\ell V\times{\Bbb F}\ell W
\quad \hookrightarrow \quad{\Bbb F}\ell(V\oplus W)
$$
by $\psi_P(\Edot,\Fdot)_j  =
\Span{E_i,F_{i'}\,|\, p_i, p^c_{i'} \leq j}$.
Equivalently, if $e_1,\ldots,e_n$ is a basis for $V$ and
$f_1,\ldots,f_m$ a basis for $W$, then
$\psi_P(\SPan{e_1,\ldots,e_n},\SPan{f_1,\ldots,f_m})=
\SPan{g_1,\ldots,g_{n+m}}$, where $g_{p_i}=e_i$ and $g_{p^c_i}=f_i$. 
From this, it follows that if $\Edot,\Epdot\in{\Bbb F}\ell V$ and
$\Fdot,\Fpdot\in{\Bbb F}\ell W$ are pairs of opposite flags, then 
$\psi_P(\Edot,\Fdot)$ and $\psi_{w_0^{(n+m)}P}(\Epdot,\Fpdot)$ are
opposite flags in $V\oplus W$.

\begin{lem}\label{lem:product_subset}
Let $P,Q\in{[n+m]\choose n}$,  $v\in {\cal S}_n$, and 
$w\in {\cal S}_m$.
Then, for $\Edot\in{\Bbb F}\ell V$ and $\Fdot\in{\Bbb F}\ell W$,
\begin{eqnarray*}
\psi_P\left(X_{w_0^{(n)}v}\Edot \times X_{w_0^{(m)}w}\Fdot\right)
&\subset& X_{w_0^{(n+m)}\varepsilon_{P,Q}(v,w)}\psi_Q(\Edot,\Fdot)\\
\psi_P\left(X_v\Edot \times X_w\Fdot\rule{0pt}{13pt}\right)
&\subset& X_{\varepsilon_{P,Q}(v,w)}\psi_{w_0^{(n+m)}Q}(\Edot,\Fdot).
\end{eqnarray*}
\end{lem}

\noindent{\bf Proof. }
For a flag $\Gdot$, define $G^\circ_j:= G_j-G_{j-1}$.
By the definition of $\psi_Q$, we have 
$E^\circ_i\subset \psi_Q(\Edot,\Fdot)^\circ_{q_i}$ and 
$F^\circ_i\subset \psi_Q(\Edot,\Fdot)^\circ_{q^c_i}$.
Since 
$$
w_0^{(n+m)}Q\quad =\quad n+m+1-q_n\ <\ \cdots\ <\ 
n+m+1-q_1,
$$ 
$E_{n+1-j} \subset \psi_{w_0^{(n+m)}Q}(\Edot,\Fdot)_{n+m+1-q_j}$,
and $F_{n+1-j} \subset \psi_{w_0^{(n+m)}Q}(\Edot,\Fdot)_{n+m+1-q'_j}$, 
the lemma is a consequence of the definitions of Schubert varieties and
$\psi_P$.
 \QEDnoskip

\begin{cor}\label{cor:geometry_product}
Let $\Edot,\Epdot\in{\Bbb F}\ell V$ and $\Fdot,\Fpdot\in{\Bbb F}\ell W$
be pairs of opposite flags and let $P\in{[n+m]\choose n}$.
Set $Q=\{m+1,\ldots,m+n\}$.
Then, for every $v\in {\cal S}_n$ and $w\in {\cal S}_m$,
\begin{eqnarray*}
\psi_P\left(X_v\Edot\times X_w\Fdot\right)  &=& 
X_{\varepsilon_{P,[n]}(v,w)}\psi_{Q}(\Edot,\Fdot)
\bigcap X_{\varepsilon_{P,Q}(e ,e )}
\psi_{[n]}(\Epdot,\Fpdot)\\
 &=& 
X_{\varepsilon_{P,[n]}(v,e )}\psi_{Q}(\Edot,\Fpdot)
\bigcap X_{\varepsilon_{P,Q}(e ,w)}
\psi_{[n]}(\Epdot,\Fdot)\\
 &=& 
X_{\varepsilon_{P,[n]}(e ,w)}\psi_{Q}(\Epdot,\Fdot)
\bigcap X_{\varepsilon_{P,Q}(v,e )}
\psi_{[n]}(\Edot,\Fpdot)\\
 &=& 
X_{\varepsilon_{P,[n]}(e ,e )}\psi_{Q}(\Epdot,\Fpdot)
\bigcap X_{\varepsilon_{P,Q}(v,w)}
\psi_{[n]}(\Edot,\Fdot).
\end{eqnarray*}
\end{cor}

\noindent{\bf Proof. }
Since $w_0^{(n+m)}[n] = Q$,
$X_e \Edot={\Bbb F}\ell V$,  and 
$X_e \Fdot={\Bbb F}\ell W$, 
Lemma~\ref{lem:product_subset} shows that 
$\psi_P\left(X_v\Edot\times X_w\Fdot\right)$
is a subset of any of the four intersections.
Equality follows as they have the same dimension.
Indeed, for $x,z\in{\cal S}_n$ and $y,u\in{\cal S}_m$,
\begin{eqnarray*}
\ell(\varepsilon_{P,[n]}(x,y))&=&
\ell(x)+\ell(y)+\#\{i\in[n],j\in[m]\,|\,p_i>p^c_j\}\\
\ell(\varepsilon_{P,Q}(z,u))&=&
\ell(z)+\ell(u)+\#\{i\in[n],j\in[m]\,|\,p^c_j>p_i\}.
\end{eqnarray*}
Thus $\ell(\varepsilon_{P,[n]}(x,y)) + \ell(\varepsilon_{P,Q}(z,u))
= \ell(x)+\ell(y)+\ell(z)+\ell(u) + n\cdot m$
and so 
$$
{n+m\choose 2} -\ell(\varepsilon_{P,[n]}(x,y)) - 
\ell(\varepsilon_{P,Q}(z,u))\ =\ 
{n\choose 2}+{m\choose 2}-\ell(x)-\ell(y)-\ell(z)-\ell(u).
$$
If $(x,y,z,u)$ is one of $(v,w,e ,e ), (v,e ,e ,w),
(e ,w,v,e ), (e ,e ,v,w)$, then these are, respectively, the 
the dimension of one of the intersections and the dimension of 
$X_v\Edot\times X_w\Fdot$.
\QEDnoskip

\begin{cor}\label{cor:product_push_forward}
Let $Q=\{m+1,\ldots,m+n\}=w_0^{(n+m)}[n]$.
For every $v\in{\cal S}_n$, $w\in{\cal S}_m$, and 
$P\in{[n+m]\choose n}$, 
the following identities hold in $H^*{\Bbb F}\ell_{n+m}$:
\medskip

\noindent$\hspace{.8in}
{\frak S}_{\varepsilon_{P,[n]}(v,w)}\cdot
{\frak S}_{\varepsilon_{P,Q}(e ,e )} \quad=\quad
{\frak S}_{\varepsilon_{P,[n]}(v,e )}\cdot
{\frak S}_{\varepsilon_{P,Q}(e ,w)} \quad=\quad$\smallskip

\hfill$
{\frak S}_{\varepsilon_{P,[n]}(e ,w)}\cdot
{\frak S}_{\varepsilon_{P,Q}(v,e )} \quad=\quad
{\frak S}_{\varepsilon_{P,[n]}(e ,e )}\cdot
{\frak S}_{\varepsilon_{P,Q}(v,w)},\hspace{.8in}
$\medskip

\noindent
and this common cohomology class is 
$(\psi_P)_*({\frak S}_v\otimes{\frak S}_w)$.
\end{cor}

\begin{thm}\label{thm:many_identities}
Let  $x\in {\cal S}_{n+m}$ and $P\in{[n+m]\choose n}$.
Then
\begin{enumerate}
\item[({\em i})]\begin{minipage}[t]{5in} \mbox{ }\vspace{-18pt}
\begin{eqnarray*}
\psi_P^*{\frak S}_x &=&
\sum_{v\in{\cal S}_n,\ w\in{\cal S}_m}
c^{\varepsilon_{P,[n]}(v,w)}_{\varepsilon_{P,[n]}(e ,e )\ x}
\ {\frak S}_v\otimes{\frak S}_w\\
 &=&
\sum_{v\in{\cal S}_n,\ w\in{\cal S}_m}
c^{\varepsilon_{P,[n]}(v,w_0^{(m)})}_%
{\varepsilon_{P,[n]}(e ,w_0^{(m)}w)\ x}
\ {\frak S}_v\otimes{\frak S}_w\\
 &=&
\sum_{v\in{\cal S}_n,\ w\in{\cal S}_m}
c^{\varepsilon_{P,[n]}(w_0^{(n)},w)}_%
{\varepsilon_{P,[n]}(w_0^{(n)}v,e )\ x}
\ {\frak S}_v\otimes{\frak S}_w\\
 &=&
\sum_{v\in{\cal S}_n,\ w\in{\cal S}_m}
c^{\varepsilon_{P,[n]}(w_0^{(n)},w_0^{(m)})}_%
{\varepsilon_{P,[n]}(w_0^{(n)}v,w_0^{(m)}w)\ x}
\ {\frak S}_v\otimes{\frak S}_w
\end{eqnarray*}\end{minipage}\medskip
\item[({\em ii})] Let $Q=\{m+1,\ldots,m+n\}$.
For every $v\in{\cal S}_n$ and $w\in {\cal S}_m$, we have
$$
\begin{array}{ccccccc}
c^{\varepsilon_{P,[n]}(v,w)}_{\varepsilon_{P,[n]}(e ,e )\ x} & = &
c^{\varepsilon_{P,[n]}(v,w_0^{(m)})}_%
{\varepsilon_{P,[n]}(e ,w_0^{(m)}w)\ x} & = &
c^{\varepsilon_{P,[n]}(w_0^{(n)},w)}_%
{\varepsilon_{P,[n]}(w_0^{(n)}v,e )\ x} & = &
c^{\varepsilon_{P,[n]}(w_0^{(n)},w_0^{(m)})}_%
{\varepsilon_{P,[n]}(w_0^{(n)}v,w_0^{(m)}w)\ x} \\
||&&||&&||&&||\\
c^{\varepsilon_{P,Q}(v,w)}_{\varepsilon_{P,Q}(e ,e )\ x} & = &
c^{\varepsilon_{P,Q}(v,w_0^{(m)})}_%
{\varepsilon_{P,Q}(e ,w_0^{(m)}w)\ x} & = &
c^{\varepsilon_{P,Q}(w_0^{(n)},w)}_%
{\varepsilon_{P,Q}(w_0^{(n)}v,e )\ x} & = &
c^{\varepsilon_{P,Q}(w_0^{(n)},w_0^{(m)})}_%
{\varepsilon_{P,Q}(w_0^{(n)}v,w_0^{(m)}w)\ x}
\end{array}
$$
\end{enumerate}
\end{thm}

\begin{rem}
{\em 
Each structure constant in ({\em ii}) is of the form $c^{\zeta y}_{y\, x}$, 
where $\zeta$ is, respectively, 
$v\times w, v\times \overline{w}^{-1}, \overline{v}^{-1}\times w$, and
$\overline{v}^{-1}\times\overline{w}^{-1}$.
Each interval $[y,\zeta y]$
is isomorphic to $[e ,v]\times[e ,w]$.
This is consistent with the expectation that the $c^z_{y\,x}$
should only depend upon $[y,z]$ and $x$.
}\end{rem}

\noindent{\bf Proof. }
In ({\em ii}), the second row is a consequence
of the first as  
$c^z_{y\,x} = c^{w_0^{(n+m)}y}_{w_0^{(n+m)}z\: x}$, for 
$x,y,z\in {\cal S}_{n+m}$.
The first row of equalities is a consequence of the
identities in ({\em i}).
For ({\em i}), there exist integral constants $d^{v\,w}_x$ 
defined by the identity
$$
\psi^*_P {\frak S}_x \quad=\quad
\sum d^{v\,w}_x\ {\frak S}_v\otimes{\frak S}_w.
$$
Since the Schubert basis is self-dual with respect to the  intersection
pairing, we have 
\begin{eqnarray*}
d^{v\,w}_x &=& \deg \left(\psi^*_P {\frak S}_x\cdot
({\frak S}_{w_0^{(n)}v}\otimes{\frak S}_{w_0^{(m)}w})\right)\\
&=& \deg\left({\frak S}_x\cdot
(\psi_P)_*({\frak S}_{w_0^{(n)}v}\otimes{\frak S}_{w_0^{(m)}w})\right).
\end{eqnarray*}
Each expression for 
$(\psi_P)_*({\frak S}_{w_0^{(n)}v}\otimes{\frak S}_{w_0^{(m)}w})$
of Corollary~\ref{cor:product_push_forward} yields one of the sums in
({\em i}).
For example, the last expression in
Corollary~\ref{cor:product_push_forward} yields
\begin{eqnarray*}
d^{v\,w}_x &=&
\deg\left({\frak S}_x\cdot
{\frak S}_{\varepsilon_{P,[n]}(e ,e )}\cdot
{\frak S}_{w_0^{(n+m)}\varepsilon_{P,[n]}(v,w)}\right)\\
&=& c^{\varepsilon_{P,[n]}(v,w)}_{\varepsilon_{P,[n]}(e ,e )\ x},
\end{eqnarray*}
since 
$w_0^{(n+m)}\varepsilon_{P,[n]}(v,w) = 
\varepsilon_{P,w_0^{(n+m)}[n]}\left(w_0^{(n)}v,w_0^{(m)}w\right)$.
\QEDnoskip

\begin{cor}\label{cor:v_times_w}
Let $u,v,w\in {\cal S}_n$ and $x,y,z\in {\cal S}_\infty$.
Then $c^{ w\times z}_{u\times x\ \,v\times y} = 
c^w_{u\,v}\cdot c^z_{x\,y}$.
\end{cor}

\noindent{\bf Proof. }
Choose $m$ so that $x,y,z\in{\cal S}_m$.
Since $\varepsilon_{[n],[n]}(u,x)=u\times x$, the first identity of
Theorem~\ref{thm:many_identities} ({\em i}) implies 
$\psi^*_{[n]} {\frak S}_{u\times x} =
{\frak S}_u \otimes {\frak S}_x$.
Then
\begin{eqnarray*}
c^w_{u\,v}\cdot c^z_{x\,y} &=&
\deg\left( ({\frak S}_u \otimes {\frak S}_x)\cdot
({\frak S}_v \otimes {\frak S}_y)\cdot 
({\frak S}_{w_0^{(n)}w} \otimes {\frak S}_{w_0^{(m)}z})\right)\\
 &=&
\deg\left( \psi^*_{[n]}({\frak S}_{u\times x}\cdot
{\frak S}_{v\times y})\cdot 
({\frak S}_{w_0^{(n)}w} \otimes {\frak S}_{w_0^{(m)}z})\right)\\
 &=&
\deg\left( {\frak S}_{u\times x}\cdot
{\frak S}_{v\times y}\cdot 
{\frak S}_{w_0^{(n+m)}(w\times z)}\right)\\
 &=&c^{ w\times z}_{u\times x\ \,v\times y},
\end{eqnarray*}
as
$(\psi_{[n]})_*{\frak S}_{w_0^{(n+m)}(w\times z)}
= {\frak S}_{w_0^{(n)}w} \otimes {\frak S}_{w_0^{(m)}z}$,
by Corollary~\ref{cor:product_push_forward}.
\QEDnoskip

\subsection{Maps ${\Bbb Z}[x_1,x_2,\ldots]\rightarrow
{\Bbb Z}[y_1,y_2,\ldots,z_1,z_2,\ldots]$}
\label{sec:substitution}

Let $P\subset {\Bbb N}$, define $P^c:= {\Bbb N}-P$, and suppose 
 $P^c$ is infinite.
Enumerate $P$ and $P^c$ as follows:
$$
\begin{array}{rcl}
P&:& p_1\ <\ p_2\ <\ \left\{\begin{array}{lll}
\cdots <p_s&\quad&\mbox{if } \#P=s\\
\cdots &&\mbox{otherwise}\end{array}\right.\\
P^c&:&p^c_1\ <\ p^c_2\ <\ \cdots\end{array}
$$
Define 
$\Psi_P: {\Bbb Z}[x_1,x_2,\ldots]\rightarrow
{\Bbb Z}[y_1,y_2,\ldots,z_1,z_2,\ldots]$
by
$$
x_{p_i}\ \longmapsto\ y_i \qquad 
x_{p^c_i}\ \longmapsto\ z_i.
$$
Then there exist integer constants $d_w^{u\,v}(P)$ for 
$u,v,w\in {\cal S}_\infty$ defined by the identity:
$$
\Psi_P({\frak S}_w(x))\ =\ 
\sum_{u,v} d_w^{u\,v}(P)\:{\frak S}_u(y)\;{\frak S}_v(z).
$$
For $l,d\in{\Bbb N}$ and $R\subset \{d+1,\ldots,d+2l\}$ with $\#R=l$,
define $\overline{P}(l,d,R):= (P\bigcap [d])\bigcup R$.

\begin{thm}\label{thm:substitution_constants}
Let $P\subset{\Bbb N}$ and $w\in {\cal S}_\infty$.
For any integers $l>\ell(w)$ and $d$ exceeding the last descent of $w$
and any subset $R$ of $\{d+1,\ldots,d+2l\}$ of cardinality $l$, set 
$n:=\#\overline{P}(l,d,R)$, $m:=d+2l-n$, and 
$\pi :=\varepsilon_{\overline{P}(l,d,R),\:[n]}(e,e)$.
Then $d_w^{u\,v}(P) = 0$ unless $u\in {\cal S}_n$ and $\in {\cal S}_m$,
and in that case,
$$
d_w^{u\,v}(P)\ =\ c^{(u\times v)\pi}_{\pi\; w}.
$$
Moreover, $d_w^{u\,v}(P) \neq 0$ implies that $a:= \# P\bigcap [d]$
exceeds the last descent of $u$ and $b:=d-a$ exceeds the last descent of
$v$.
\end{thm}

\begin{rem}\label{rem:I_P}
{\em 
Theorem~\ref{thm:substitution_constants}
generalizes~\cite[1.5]{Lascoux_Schutzenberger_structure_de_Hopf} 
(see also~\cite[4.19]{Macdonald_schubert}) where
it is shown that  $d_w^{u\,v}([a])\geq 0$.
Define $I_P$ to be
$$
\{\varepsilon_{\overline{P}(l,l,R),[n]}(e,e)\;|\;
l\in {\Bbb N},\ n=l+\#(P\bigcap [l]),\mbox{ and }
 R\subset\{l+1,\ldots,3l\},\# R=l\}.
$$
For $w\in {\cal S}_n$, let $N$ be an integer such that $N/3$ exceeds
both the last descent the length of $w$.
If $\pi\in I_P$ with $\pi\not\in{\cal S}_N$, then 
$\pi = \varepsilon_{\overline{P}(l,d,R),[n]}(e,e)$ for $l,d,R$ satisfying
the conditions of Theorem~\ref{thm:substitution_constants}
and so 
$d_w^{u\,v}(P)\ =\ c^{(u\times v)\pi}_{\pi\; w}$
for every $\pi \in I_P-{\cal S}_N$, which establishes
Theorem~\ref{thm:substitution}.}
\end{rem}

Apply the ring homomorphism $\Psi_P$ to both sides of the
product:
$$
{\frak S}_w(x)\;{\frak S}_\gamma(x)\quad=\quad
\sum_\zeta c^\zeta_{w\,\gamma}\,{\frak S}_\zeta(x).
$$
If we expand this in terms of 
${\frak S}_\eta(y)\; {\frak S}_\xi(z)$ and 
equate the coefficients, we get a corollary.

\begin{cor}\label{cor:brutally_complicated}
Let $w,\gamma,\eta,\xi\in{\cal S}_\infty$, and $P\subset{\Bbb N}$.
Then there exists an integer $N\in {\Bbb N}$ such that if $\pi\in
I_P-{\cal S}_N$, then 
$$
\sum_{\zeta} c^{(\eta\times\xi)\pi}_{\pi\  \zeta}\:
c^\zeta_{w\;\gamma} \quad = \quad
\sum_{u,v,\alpha,\beta} c^{(u\times v)\pi}_{\pi\ w}\:
 c^{(\alpha\times\beta)\pi}_{\pi\ \gamma}\:
c^\eta_{u\;\alpha}\:c^\xi_{v\;\beta}.
$$
\end{cor}

\noindent{\bf Proof of Theorem~\ref{thm:substitution_constants}. }
First, a Schubert polynomial 
${\frak S}_\pi(x)\in{\Bbb Z}[x_1,\ldots,x_s]$ if and only if $s$ exceeds
the last descent of
$\pi$~\cite{Lascoux_Schutzenberger_polynomes_schubert} 
(see also~\cite[4.13]{Macdonald_schubert}).
Thus, ${\frak S}_w(x)\in{\Bbb Z}[x_1,\ldots,x_d]$, and 
if $d^{u\,v}_w(P)\neq 0$, then 
${\frak S}_u(y)\in{\Bbb Z}[y_1,\ldots,y_a]$ and 
${\frak S}_v(z)\in{\Bbb Z}[z_1,\ldots,z_b]$,
hence $a$, respectively, $b$, exceeds the last descent of $u$,
respectively $v$.
Since $\deg{\frak S}_w(x)\leq l$, both $\deg{\frak S}_u(y)$ and 
 $\deg{\frak S}_v(z)$ are at most $l$.
Consider the commutative diagram
$$
\begin{picture}(395,70)
\put(0,50){${\Bbb Z}[x_1,\ldots,x_d]$}
\put(110,50){${\Bbb Z}[x_1,\ldots,x_{n+m}]$}
\put(130,5){$H^* {\Bbb F}\ell_{n+m}$}
\put(258,50){${\Bbb Z}[y_1,\ldots,y_n,z_1,\ldots,z_m]$}
\put(280,5){$H^*{\Bbb F}\ell_n\otimes H^*{\Bbb F}\ell_m$}
\put(72,54){\oval(4,4)[l]}\put(72,52){\vector(1,0){33}}
\put(86,57){$\iota$}
\put(195,52){\vector(1,0){58}}\put(215,57){$\overline{\Psi_P}$}
\put(185,7){\vector(1,0){90}}\put(215,14){$\psi^*_{\overline{P}}$}
\put(150,43){\vector(0,-1){20}}\put(150,43){\vector(0,-1){25}}
\put(320,43){\vector(0,-1){20}}\put(320,43){\vector(0,-1){25}}
\end{picture}
$$
Here, $\overline{\Psi_P}$ is the restriction of $\Psi_{\overline{P}}$ to 
${\Bbb Z}[x_1,\ldots,x_{n+m}]$.
The vertical arrows are injective on the module
${\Bbb Z}\Span{x_1^{\alpha_1}\cdots x_d^{\alpha_d}\;|\; 
\alpha_i \leq l}$
and its image
$$
{\Bbb Z}\Span{y_1^{\beta_1}\cdots y_a^{\beta_a}
z_1^{\gamma_1}\cdots z_b^{\gamma_b}
\;|\; \beta_i,\gamma_j \leq l}
\quad\subset\quad{\Bbb Z}[y_1,\ldots,y_n,z_1,\ldots,z_m].
$$
Moreover, since $P\bigcap[d]=\overline{P}\bigcap [d]$, the composition,
$\Psi_{\overline{P}}\circ\iota$, of the top row coincides with 
$\Psi_{P}\circ\iota$.
Since ${\frak S}_w(x)\in  {\Bbb Z}
\Span{x_1^{\alpha_1}\cdots x_d^{\alpha_d}\;|\; \alpha_i \leq l}$,
the cohomological formula for $\psi_{\overline{P}}^*({\frak S}_w)$  in
Theorem~\ref{thm:many_identities}  computes $\Psi_P({\frak S}_w(x))$.
\QED

In the statement of the Theorem~\ref{thm:substitution_constants},
$\ell(w)$ could be replaced by  $\max_i\{\deg_{x_i}({\frak S}_w(x))\}$.

\subsection{Products of Grassmannians}
Let $k\leq n$ and $l\leq m$ be integers, $V\simeq {\Bbb C}^n$, and
$W\simeq{\Bbb C}^m$.
Define 
$\varphi_{k,l}:\Gr_k V\times \Gr_l W \hookrightarrow \Gr_{k+l}(V\oplus W)$
by 
$$
\varphi_{k,l}\ :\ (H,K)\quad \longmapsto\quad H\oplus K.
$$

\begin{thm}\label{thm:grassmann_product}
\mbox{ }
\begin{enumerate}
\item[({\em i})] For every Schubert class $S_\lambda\in
H^*\Gr_{k+l}V\oplus W$, 
$$
\varphi^*_{k,l}(S_\lambda) \quad=\quad
\sum_{\mu,\nu}
c^\lambda_{\mu\,\nu} S_\mu\otimes S_\nu.
$$
\item[({\em ii})] If $S_{\mu^c}\otimes S_{\nu^c}\in
H^*\Gr_k V\otimes  H^*\Gr_l W$, then 
$$
(\varphi_{k,l})_*(S_{\mu^c}\otimes S_{\nu^c}) \quad =\quad
\sum_\lambda c^\lambda_{\mu\,\nu} S_{\lambda^c},
$$
where $\lambda^c,\mu^c$, and $\nu^c$ are defined by
$\mu^c_i=n-k-\mu_{k+1-i}$,
$\nu^c_i=m-l-\nu_{l+1-i}$, and 
$\lambda^c_i=m+n-k-l-\lambda_{k+l+1-i}$.
\end{enumerate}
\end{thm}

\begin{rem}{\em 
If $-x_1,\ldots,-x_k$ are the Chern roots of the tautological $k$-plane
bundle over $\Gr_k V$, and $-y_1,\ldots,-y_l$ those of the tautological
$l$-plane bundle over $\Gr_l W$, 
and $f\in H^*\Gr_{k+l}V\oplus W$
(which is a symmetric polynomial in the negative Chern roots of the
tautological bundle over $\Gr_{k+l}V\oplus W$).
Then
$$
\varphi_{k,l}^* f\ = \ f(x_1,\ldots,x_k,y_1,\ldots,y_l).
$$
Let $\Lambda=\Lambda(z)$ be the ring of symmetric functions, 
which is the inverse limit (in the category of graded
rings) of the rings of symmetric
polynomials in the variables $z_1,\ldots,z_n$. 
Fixing $\lambda$ and choosing $k,l,n$, and $m$ large enough 
gives a new proof of~\cite[I.5.9]{Macdonald_symmetric}:
}\end{rem}

\begin{prop}[{\cite[I.5.9]{Macdonald_symmetric}}]
Let $\lambda$ be a partition and $x$, $y$ be infinite sets of variables.
Then
$$
S_\lambda(x,y)\ =\ \sum_{\mu,\nu} 
c^\lambda_{\mu\,\nu}\; S_\mu(x)\cdot S_\nu(y),
$$
where $S_\mu$ denotes the Schur function basis of the ring
$\Lambda$ of symmetric functions.
\end{prop}

If we define a linear map $\Delta:\Lambda(z) \rightarrow 
\Lambda(x)\otimes_{\Bbb Z}\Lambda(y)$ by 
$\Delta(f(z))= f(x,y)$, then $\Delta$ is induced by the maps
$\varphi^*_{k,l}$. 
Moreover, the obvious commutative diagrams of spaces give a new proof
of~\cite[I.5.24]{Macdonald_symmetric}, that 
$\Lambda$ is a cocommutative Hopf algebra with comultiplication
$\Delta$.\medskip

%  \begin{picture}(375,70)
%  \put(0,55){$\Gr_i{\Bbb C}^j\times\Gr_k{\Bbb C}^n\times\Gr_l{\Bbb C}^m$}
%  \put(20,5){$\Gr_{i+k}{\Bbb C}^{j+n}\times\Gr_l{\Bbb C}^m$}
%  \put(240,55){$\Gr_i{\Bbb C}^j\times\Gr_{k+l}{\Bbb C}^{n+m}$}
%  \put(260,5){$\Gr_{i+k+l}{\Bbb C}^{j+n+m}$}
%  \put(175,58){\vector(1,0){60}}
%  \put(157,8){\vector(1,0){100}}
%  \put(85,48){\vector(0,-1){30}}
%  \put(307,47){\vector(0,-1){32}}
%  \put(190,62){\scriptsize$1\times \varphi_{k,l}$}
%  \put(191,12){\scriptsize$\varphi_{i+k,l}$}
%  \put(53,33){\scriptsize$\varphi_{i,k}\times 1$}
%  \put(311,33){\scriptsize$\varphi_{i,k+l}$}
%  \end{picture}

\noindent{\bf Proof of Theorem~\ref{thm:grassmann_product}.}
The first statement is  a consequence of the second:
Schubert classes form a basis for the cohomology ring, so there exist
integral constants $d^{\mu\,\nu}_\lambda$ such that 
$$
\varphi^*_{k,l}(S_\lambda)\quad=\quad
\sum_{\mu,\nu} d^{\mu\,\nu}_\lambda S_\mu\otimes S_\nu.
$$
Since the Schubert basis diagonalizes the intersection pairing,
$$
d^{\mu\,\nu}_\lambda \quad=\quad 
\deg(\varphi^*_{k,l}(S_\lambda)\cdot (S_{\mu^c}\otimes S_{\nu^c})).
$$
Apply $(\varphi_{k,l})_*$ and use the second assertion to 
obtain
\begin{eqnarray*}
d^{\mu\,\nu}_\lambda&=&
\deg(S_\lambda\cdot(\varphi_{k,l})_*(S_{\mu^c}\otimes S_{\nu^c}))\\
&=& S_\lambda \cdot \sum_\kappa c^\kappa_{\mu\,\nu} S_{\kappa^c}\\
&=& c^\lambda_{\mu\,\nu}.
\end{eqnarray*}

The second assertion is a consequence of the following lemma.

\begin{lem}\label{lem:grassmann_computation}
Suppose $\mu,\nu$ are partitions with $\mu\subset (n-k)^k$ and
$\nu\subset(m-l)^l$. 
Let $\Edot\in {\Bbb F}\ell V$ and  $\Fdot\in{\Bbb F}\ell W$
and let $\Gpdot$ be any flag opposite to $\psi_{[n]}(\Edot,\Fdot)$ with
$G'_m=W$.
Then 
$$	%\begin{equation}\label{eq:grass_product}
\varphi_{k,l}\left(\Omega_{\mu^c}\Edot\times
\Omega_{\nu^c}\Fdot\right) \quad =\quad
\Omega_{\rho^c}\psi_{[n]}(\Edot,\Fdot)
\bigcap \Omega_{(m-l)^k}\Gpdot,
\eqno(4.7.1)
$$	%\end{equation}
where $\rho$ is the partition
$$
\nu_1+(n-k)\  \geq\  \cdots\ \geq\  \nu_l+(n-k)\  \geq \ 
\mu_1\ \geq \cdots\ \geq \mu_k.
$$
\end{lem}

We finish the proof of Theorem~\ref{thm:grassmann_product}.
Lemma~\ref{lem:grassmann_computation} implies 
\begin{eqnarray*}
\left(\varphi_{k,l}\right)_*\left(S_{\mu^c}\otimes S_{\nu^c}\right) 
&=&
\left[ \Omega_{\rho^c}\psi_{[n]}(\Edot,\Fdot) \bigcap 
\Omega_{(n-k)^l}\Gpdot\right]\\
&=&
\sum_\lambda c^{\lambda^c}_{\rho^c\: (n-k)^l}\; S_{\lambda^c}.
\end{eqnarray*}
Since $\deg(S_\alpha\cdot S_\beta\cdot S_\gamma) =
c^{\alpha^c}_{\beta\,\gamma}$, we see that 
$$
c^{\lambda^c}_{\rho^c\: (n-k)^l} \ =\ 
c^\rho_{(n-k)^l\; \lambda} \ =\ 
c^{\rho/(n-k)^l}_\lambda\ =\ 
c^{\mu\coprod \nu}_\lambda \ =\ 
c^\lambda_{\mu\,\nu}.
$$
Here, $\mu \coprod\nu$ is a skew partition with two
components $\mu$ and $\nu$ and the last equality is a special case
of~(1.3.1) %\ref{eq:disjoint_skew}) 
in \S\ref{sec:Schur_identities}. 
\QED

\noindent{\bf Proof of Lemma~\ref{lem:grassmann_computation}. }
Since 
$\Omega_{(n-k)^l}\Gpdot 
= \{M\in \Gr_{k+l}V\oplus W\,|\, \dim M\bigcap G'_m\geq l\}$
and $G'_m=W$, we see that 
$\varphi_{k,l}(\Gr_kV\times\Gr_lW)\subset\Omega_{(n-k)^l}\Gpdot$.
It is an exercise in the definition of the Schubert varieties involved
and of $\psi_{[n]}(\Edot,\Fdot)$ to see that
$$
\varphi_{k,l}(\Omega_{\mu^c}\Edot\times\Omega_{\nu^c}\Fdot)
\ \subset\ \Omega_{\rho^c}\psi_{[n]}(\Edot,\Fdot),
$$
which shows the inclusion $\subset$ in~(4.7.1). %\ref{eq:grass_product}).
Equality follows as they have the same dimension:
The intersection has dimension
$|\rho|-|(n-k)^l|\ =\ |\mu|+|\nu|$,
the dimension of $\Omega_{\mu^c}\Edot\times\Omega_{\nu^c}\Fdot$.
\QEDnoskip

\section{Symmetries of the Littlewood-Richardson coefficients}
\label{sec:geometry}

\subsection{Proof of Theorem~\ref{thm:B} ({\em ii})}\label{sec:proof_thm_B}

Combining Lemma~\ref{lem:restriction} with Theorem~\ref{thm:coeff_sum},
we deduce:

\begin{lem}\label{lem:shape_reduction}
Suppose $x\leq_k z$ and $x(p)=z(p)$.
Let $k'=k-1$ if $p<k$ and $k'=k$ otherwise.
Then for all partitions $\lambda$, we have
$$
c^z_{x\;v(\lambda,k)}\quad =\quad 
c^{z/_p}_{x/_p\:v(\lambda,k')}.
$$
\end{lem}

Note that $zx^{-1}$ and $z/_p(x/_p)^{-1}$ are shape-equivalent, 
by Lemma~\ref{lem:expanding_bruhat}~({\em iv}).

\begin{lem}\label{lem:skew_critical}
Let $x,z,u,w\in{\cal S}_n$.
Suppose $x\leq_kz$, $u\leq_kw$, and $zx^{-1}=wu^{-1}$.
Further suppose that $w$ is Grassmannian with descent $k$, the
permutation $wu^{-1}$ has no fixed points, and, for $k<i\leq n$,
$u(i)=x(i)$.
Then, for all partitions $\lambda$ with at most $k$ parts,
$$
c^w_{u\,v(\lambda,k)}\quad =\quad c^z_{x\,v(\lambda,k)}.
$$
\end{lem}

\noindent{\bf Proof of Theorem~\ref{thm:B}~({\em ii}) using
Lemma~\ref{lem:skew_critical}. }
We reduce Theorem~\ref{thm:B}~({\em ii}) to the special case of
Lemma~\ref{lem:skew_critical}. 
First, by Lemma~\ref{lem:shape_reduction}, it suffices to prove 
Theorem~\ref{thm:B}~({\em ii}) when 
$x,z,u,w\in{\cal S}_n$, $k=l$, with $wu^{-1}=zx^{-1}$ and the
permutation $wu^{-1}$ has no fixed points.

Define $s\in {\cal S}_n$ by
$$
s(i)\quad:=\quad\left\{\begin{array}{lcl}u(i)&&1\leq i\leq k\\
x(i)&& k<i\leq n\end{array}\right.
$$
and set $t:=wu^{-1}s$.
Then $s\leq_k t$ and 
$$
t(i)\quad=\quad\left\{\begin{array}{lcl}w(i)&&1\leq i\leq k\\
z(i)&& k<i\leq n\end{array}\right..
$$
It suffices to show separately that $c^w_{u\:v(\lambda,k)}$ and
$c^t_{s\:v(\lambda,k)}$ each equal $c^z_{t\:v(\lambda,k)}$.
Thus we may further assume $u(i)=x(i)$ for $1\leq i\leq k$
or $u(i)=x(i)$  for $k<i\leq n$.

Suppose that $u(i)=x(i)$ for  $1\leq i\leq k$.
If  for $v\in {\cal S}_n$, $\overline{v}:=w_0vw_0$,
$$
c^z_{x\:v(\lambda,k)}\quad=\quad c^w_{u\:v(\lambda,k)}
\quad\Longleftrightarrow\quad
c^{\overline{z}}_{\overline{x}\:\overline{v(\lambda,k)}}
\quad=\quad c^{\overline{w}}_{\overline{u}\:\overline{v(\lambda,k)}}.
$$
Set $l=n-k$ and $\lambda^t$ the partition conjugate to $\lambda$.
Then 
$\overline{x}\leq_{k'}\overline{z}$, 
$\overline{u}\leq_{k'}\overline{w}$, 
$\overline{z}(\overline{x}^{-1})=\overline{w}{\overline{u}^{-1}}$, 
$\overline{v(\lambda,k)}= v(\lambda^t,l)$, and 
$\overline{x}(i)=\overline{u}(i)$ for $l<i\leq n$.
Thus we may assume $x(i)=u(i)$ for $1\leq i\leq k$.

Finally,  there is a Grassmannian permutation $t\in {\cal S}_n$  with
descent $k$  and  a permutation $s\in{\cal S}_n$ such that 
$t=wu^{-1}s$.
Thus it suffices to further assume that $w$ is Grassmannian with
descent $k$, the situation of Lemma~\ref{lem:skew_critical}.
\QED

We prove Lemma~\ref{lem:skew_critical} by studying two 
intersections of Schubert varieties and their image under the projection 
${\Bbb F}\ell V\twoheadrightarrow \Gr_k V$.
Let $e_1,\ldots,e_n$ be a basis for $V$ and set 
$\Fdot=\SPan{e_1,\ldots,e_n}$.
Let $M(w)\subset M_{n\times n}{\Bbb C}$ be the set of matrices
satisfying the conditions:
\begin{enumerate}
\item[(a)]  $M(w)_{i,w(i)}=1$
\item[(b)] $M(w)_{i,j}=0$ if either $w(i)<j$ or else $w^{-1}(j)<i$.
\end{enumerate}
Then $M(w)\simeq{\Bbb C\,}^{\ell(w)}$ as the only unconstrained
entries of $M(w)$ are $M(w)_{i,j}$ when $j<w(i)$ and $i<w^{-1}(j)$, and
there are $\ell(w)$ such entries.

\begin{ex}
{\em 
Let $w=25134\in {\cal S}_5$, a Grassmannian permutation with descent 2.
Then $M(w)$ is the set of matrices
$$
\left\{
\left[\begin{array}{ccccc}
  a   & 1 &  0   &  0   & 0\\
  b   & 0 &  c   &  d   & 1\\
  1   & 0 &  0   &  0   & 0\\
  0   & 0 &  1   &  0   & 0\\
  0   & 0 &  0   &  1   & 0\end{array}\right]
\ \raisebox{-30pt}{\rule{.5pt}{68pt}}\, \mbox{ for every } 
(a,b,c,d)\in{\Bbb C}\,^4\right\}
$$
}\end{ex}

Fix a basis $e_1,\ldots,e_n$ for $V$.
For $\alpha\in M(w)$, and $1\leq i\leq n$, define the vector
$f_i(\alpha) := \sum_j \alpha_{i,j}e_j$.
Then $ f_1(\alpha),\ldots,f_n(\alpha)$ are the 
`row vectors' of the matrix $\alpha$ and they form a basis for $V$ as
$\alpha$ has determinant  $(-1)^{\ell(w)}$.
Set $\Edot(\alpha)=\SPan{f_1(\alpha),\ldots,f_n(\alpha)}$.
Since $f_i(\alpha)\in F_{w(i)}- F_{w(i)-1}$, we see that 
$\Edot(\alpha)\in X^\circ_{w_0w}\Fdot$.
Moreover, given $\Edot\in X^\circ_{w_0w}\Fdot$, restricted row
reduction on a basis for $\Edot$ shows there is a unique 
$\alpha\in M(w)$ with $\Edot=\Edot(\alpha)$.
In the case that $w$ is Grassmannian with descent $k$, matrices in 
$M(w)$ have a simple form: if $k<i$, then
$f_i(\alpha) = e_{w(i)}$.

For opposite flags $\Fdot,\Fpdot$, ${\frak S}_{w_0w}\cdot{\frak S}_u$
is the class Poincar\'e dual to the
fundamental cycle of $X_{w_0w}\Fdot\bigcap X_u\Fpdot$.
We  use the projection formula~(2.3.1) %\ref{eq:projection}) 
to compute the coefficient
$c^w_{u\,v(\lambda,k)}$:
\begin{eqnarray*}
c^w_{u\,v(\lambda,k)} 
&=& \deg(S_\lambda(x_1,\ldots,x_k)\cdot{\frak S}_{w_0w}\cdot{\frak S}_u)\\
 &=& \deg (\pi_k)_*(S_\lambda(x_1,\ldots,x_k)
\cdot{\frak S}_{w_0w}\cdot{\frak S}_u)\\
 &=& \deg(S_\lambda \cdot
(\pi_k)_*({\frak S}_{w_0w}\cdot{\frak S}_u))
\end{eqnarray*}

Thus Lemma~\ref{lem:skew_critical}
is a consequence of the following calculation:

\begin{lem}\label{lem:skew_calculation}
Let $u,w,x,z$ satisfy the hypotheses of Lemma~\ref{lem:skew_critical}.
Then, if $\Fdot$ and $\Fpdot$ are opposite flags in $V$,
$$
\pi_k \left(X_{w_0w}\Fdot\bigcap X_u\Fpdot\right) \quad =\quad
\pi_k \left(X_{w_0z}\Fdot\bigcap X_x\Fpdot\right).
$$
\end{lem}

\noindent{\bf Proof. }
Let $e_1,\ldots,e_n$ be a basis for $V$ such that
$\Fdot=\SPan{e_1,\ldots,e_n}$ and $\Fpdot=\SPan{e_n,\ldots,e_1}$,
and define $M(w)$ as before.
Let $A\subset M(w)$ consist of those matrices $\alpha$
such that $\Edot(\alpha)\in X^\circ_u\Fpdot$.
If $j>k$, set $g_j(\alpha)=f_j(\alpha)=e_{w(j)}$.
For $j\leq k$ construct $g_j(\alpha)$ inductively, setting
$g_j(\alpha)$ to be the intersection of $F'_{n+1-u(j)}$ and the affine
space
$f_j(\alpha) + \Span{g_i(\alpha)\,|\,i<j\mbox{ and } u(i)<u(j)}$.
Since $\Edot(\alpha)\in X_u\Fdot$ and 
$\dim E_j(\alpha)\bigcap F'_{n+1-u(j)}=
\#\{i\leq j\,|\, u(i)>u(j)\}$, this intersection consists of a single,
non-zero vector, $g_j(\alpha)$.

Then the algebraic map
$A\ni \alpha \mapsto (g_1(\alpha),\ldots,g_n(\alpha))\in V^n$  gives a
parameterized basis of $V$. 
Moreover, $\Edot(\alpha)=\SPan{g_1(\alpha),\ldots,g_n(\alpha)}$ for all 
$\alpha\in A$, and if $1\leq j\leq k$, then 
$g_j(\alpha)\in F'_{n+1-u(j)}\bigcap F_{w(j)}$.
Note that for $\alpha\in A$,
$$	%\begin{equation}\label{eq:new_flag}
\Gdot(\alpha) \ :=\
\SPan{g_{u^{-1}x(1)}(\alpha),\ldots,g_{u^{-1}x(n)}(\alpha)}
\ \in\ X_{w_0z}\Fdot\bigcap X_x\Fpdot.
\eqno(5.1.1)
$$	%\end{equation}
and thus $A$ also parameterizes a subset of $X_{w_0z}\Fdot\bigcap
X_x\Fpdot$.
Indeed, for $1\leq j\leq k$,
$g_{u^{-1}x(j)}\in F'_{n+1-x(j)}\bigcap F_{y(j)}$.
Also for $j>k$, we have $u^{-1}x(j)=j=w^{-1}y(j)$, thus
$g_j(\alpha)=f_j(\alpha)=e_{y(j)}$
and $G_j(\alpha)=E_j(\alpha)$.
The definition of Schubert cells for the flag
manifold  in \S\ref{sec:flag} then implies (5.1.1). %\ref{eq:new_flag}).

Both  cycles $X_{w_0w}\Fdot\bigcap X_u\Fpdot$ and
$X_{w_0 z}\Fdot\bigcap X_x\Fpdot$ are irreducible and have the same
dimension, $\ell(w)-\ell(u)= | w u^{-1}|$.
Since $\Gdot(\alpha)=\Gdot(\beta)$ if and only if $\alpha=\beta$,
the loci of flags $\{\Gdot(\alpha)\,|\, \alpha\in A\}$ is dense in 
$X_x\Fpdot\bigcap X_{w_0y}\Fdot$.
Finally, for $\alpha\in A$, we have $G_k(\alpha)=E_k(\alpha)$,
as $u^{-1}x$ permutes $\{1,\ldots,k\}$,
which completes the proof.
\QEDnoskip

\subsection{Proof of Theorem~\ref{thm:C}~({\em ii}). }\label{sec:proof_C}
We show that if $\zeta$ and $\eta$ are disjoint permutations and
$\lambda$ any partition,  then
$$
c^{\zeta\eta}_\lambda\quad=\quad
\sum_{\mu,\nu} c^\lambda_{\mu\,\nu}\;c^\zeta_\mu\; c^\eta_\nu.
$$

\begin{lem}\label{lem:long}
Let $\zeta,\eta\in{\cal S}_{n+m}$ be disjoint permutations.
Suppose $k\geq \#\mbox{up}_\zeta$,  $l\geq \#\mbox{up}_\eta$,
$n\geq\#\mbox{supp}_\zeta$, and $m\geq \#\mbox{supp}_\eta$.
Let $u\in{\cal S}_{n+m}$ be a permutation such that 
$u\leq _{k+l}\zeta\eta u$.
Let $Q$ be any element of ${[n+m]-\mbox{supp}_\eta\choose n}$ which
contains  $\mbox{supp}_\zeta$ for which 
$k=\# u^{-1}(Q)\bigcap [k+l]$.
Set $Q^c:= [n+m]-Q$.

Define $\zeta'\in{\cal S}_n$ and $\eta'\in{\cal S}_m$ by
$\phi_Q(\zeta')=\zeta$ and $\phi_{Q^c}(\eta')=\eta$.
Set $P=u^{-1}(Q)$, $P^c=u^{-1}(Q^c)$,  and define 
$v\in{\cal S}_n$ and $w\in{\cal S}_m$ by
$u(p_i)=q_{v(i)}$ and $u(p^c_i)=q^c_{w(i)}$,
where 
$$
\begin{array}{lll}
P\ =\ p_1\ <\ p_2\ <\ \cdots\ <\ p_n &\qquad&
P^c\ =\ p^c_1\ <\ p^c_2\ <\ \cdots\ <\ p^c_m\\
Q\ =\ q_1\ <\ q_2\ <\ \cdots\ <\ q_n &\qquad&
Q^c\ =\ q^c_1\ <\ q^c_2\ <\ \cdots\ <\ q^c_m\end{array}
$$
Then
\begin{enumerate}
\item[({\em i})]  $v\leq_k \zeta' v$ and $w\leq_l \eta' w$,
\item[({\em ii})] $u=\varepsilon_{P,Q}(v,w)$ and 
$\zeta\eta u =\varepsilon_{P,Q}(\zeta' v,\eta' w)$, and
\item[({\em iii})] For all pairs of opposite flags 
$\Edot,\Epdot\in{\Bbb F}\ell_n$
and $\Fdot,\Fpdot\in{\Bbb F}\ell_m$,
\smallskip

${\displaystyle \psi_P\left[
\left(X_{w_0^{(n)}\zeta' v}\Edot \bigcap X_v\Epdot\right) \times
\left(X_{w_0^{(m)}\eta' w}\Fdot \bigcap X_w\Fpdot\right)\right]\ =}$

\hfill ${\displaystyle
X_{w_0^{(n+m)}\zeta\eta u}\psi_Q(\Edot,\Fdot) \bigcap
X_u\psi_{w_0^{(m+n)}Q}(\Epdot,\Fpdot)}$.\qquad
\end{enumerate}
\end{lem}

\noindent{\bf Proof. }
Since $u\leq_{k+l} \zeta\eta u$,  ({\em i}) follows from
Theorem~\ref{thm:k-length} and the definitions.
The second statement is also immediate.
For ({\em iii}), Lemma~\ref{lem:product_subset}  shows the 
inclusion $\subset$.
Since $\zeta'$ is shape equivalent to $\zeta$, $\eta'$ to $\eta$, and
$\zeta$ and $\eta$ are disjoint,
$|\zeta\eta|=|\zeta'|+|\eta'|$, showing both cycles have the same
dimension, and hence are equal, as 
$\psi_Q(\Edot,\Fdot)$ and $\psi_{w_0^{(m+n)}Q}(\Epdot,\Fpdot)$ are
opposite flags.
\QED

Note that if $u\leq_k \zeta u$, then 
$$
c^\zeta_\lambda \quad=\quad
\deg(S_\lambda \cdot 
(\pi_k)_*({\frak S}_{w_0\zeta u}\cdot{\frak S}_u)).
$$
Thus the skew Littlewood-Richardson coefficients $c^\zeta_\lambda$ 
are defined by the identity in $H^* \Gr_k V$:
$$	%\begin{equation}\label{eq:homology_coefficients}
(\pi_k)_*({\frak S}_{w_0\zeta u}\cdot{\frak S}_u)
\quad=\quad \sum_{\lambda\subset (n-k)^k} 
c^\zeta_\lambda\, S_{\lambda^c}.
\eqno(5.2.1)
$$	%\end{equation}

\noindent{\bf Proof of Theorem~\ref{thm:C}~({\em ii}). }
We use the notation of Lemma~\ref{lem:long}.
The following diagram
commutes since $[k+l]=\{p_1,\ldots,p_k,p^c_1,\ldots,p^c_l\}$.
$$
\setlength{\unitlength}{2.2pt}
\begin{picture}(107,34)
\put(0,3){$\Gr_k{\Bbb C}^n\times\Gr_l{\Bbb C}^m$}
\put(75,3){$\Gr_{k+l}{\Bbb C}^{n+m}$}
\put(13.8,26){${\Bbb F}\ell_n\times{\Bbb F}\ell_m$}
\put(82,26){${\Bbb F}\ell_{n+m}$}
\put(58,7){$\varphi_{k,l}$}     \put(60,30){$\psi_P$}
\put(6.5,16){$\pi_k\times\pi_l$}  \put(91,16){$\pi_{k+l}$}
\put(53,4.5){\vector(1,0){20}}
\put(40,27.5){\vector(1,0){40}}
\put(25.5,23){\vector(0,-1){14}}
\put(89,23){\vector(0,-1){14}}
\end{picture}
$$
{}From this and Lemma~\ref{lem:long}, we see that 
$$
\pi_{k+l}\left(X_{w_0^{(n+m)}\zeta\eta u}\psi_Q(\Edot,\Fdot)
\bigcap X_u\psi_{w_0^{(m+n)}Q}(\Epdot,\Fpdot)\right)
$$
is equal to 
$$
\varphi_{k,l}\left(
\pi_k\left(X_{w_0^{(n)}\zeta' v}\Edot\bigcap X_v\Epdot\right) \times
\pi_l\left(X_{w_0^{(m)}\eta' w}\Fdot\bigcap X_w\Fpdot\right) \right).
$$
Thus
$(\pi_{k+l})_*\left({\frak S}_{w_0^{(n+m)}\zeta\eta u}
\cdot {\frak S}_u\right)$ is equal to 
$$
(\varphi_{k,l})_*\left(
(\pi_k)_*\left({\frak S}_{w_0^{(n)}\zeta' v}
\cdot{\frak S}_v\right)\otimes
(\pi_l)_*\left({\frak S}_{w_0^{(m)}\eta' w}
\cdot{\frak S}_w\right)\right).
$$
This, together with (5.2.1), %\ref{eq:homology_coefficients}), 
gives
\begin{eqnarray*}
\sum_\lambda c^{\zeta\eta}_\lambda S_{\lambda^c}
&=&
(\pi_{k+l})_*\left({\frak S}_{w_0^{(n+m)}\zeta\eta u}
\cdot {\frak S}_u\right) 
\\
&=& (\varphi_{k,l})_*\left(
\sum_\mu c^{\zeta'}_\mu S_{\mu^c} \otimes
\sum_\nu c^{\eta'}_\nu S_{\nu^c}\right)
\\
&=& 
\sum_{\mu,\nu} c^{\zeta'}_\mu c^{\eta'}_\nu
(\varphi_{k,l})_*\left(S_{\mu^c} \otimes \ S_{\nu^c}\right)
\\
&=& 
\sum_{\mu,\nu} c^{\zeta'}_\mu c^{\eta'}_\nu
\sum_\lambda c^\lambda_{\mu\,\nu} S_{\lambda^c}.
\end{eqnarray*}
This completes the proof, as $\zeta',\zeta$ and $\eta',\eta$ are shape
equivalent pairs. 
\QEDnoskip

\subsection{Theorem~\ref{thm:D}$'$ (Cyclic Shift) }\label{sec:thmd}\ 
{\em
Let $u,w,x,z\in{\cal S}_\infty$ with $u\leq_k w$ and $x\leq_l z$.
Suppose $wu^{-1}\in {\cal S}_n$ and 
$zx^{-1}$ is shape equivalent to $(wu^{-1})^{(1\,2\,\ldots\,n)^t}$,
for some  $t$.
For every partition $\lambda$, 
$$
c^w_{u\,v(\lambda,k)}\quad =\quad c^z_{x\,v(\lambda,l)}.
$$
}\medskip

\noindent{\bf Proof. }
By Theorem~\ref{thm:B}~({\em ii}), it suffices to
prove a restricted case. 
Suppose $u,w\in{\cal S}_n$, $u\leq_k w$, and $w$ is Grassmannian with
descent $k$.
The idea is to construct permutations 
$x,z\in {\cal S}_n$ with  $x\leq_k z$ and 
$zx^{-1}=(wu^{-1})^{(1\,2\,\ldots\,n)}$
for which 
$$	%\begin{equation}\label{eq:cyclic_shift}
	\pi_k\left(X_{w_0w}\Fdot\bigcap X_u\Fpdot\right)
	\quad=\quad
	\pi_k\left(X_{w_0z}\Gdot\bigcap X_x\Gpdot\right),
\eqno(5.3.1)
$$	%\end{equation}
where $e_1,\ldots,e_n$ be a basis for $V$
and the flags $\Fdot,\Fpdot,\Gdot$, and $\Gpdot$ are 
$$
\begin{array}{rclcrcl}
\Fdot & = & \SPan{e_1,\ldots,e_n} &&
\Fpdot & = & \SPan{e_n,\ldots,e_1}\\
\Gdot & = & \SPan{e_n,e_1,\ldots,e_{n-1}} &&
\Gpdot & = & \SPan{e_{n-1},\ldots,e_1,e_n}.
\end{array}
$$
Then (5.3.1) %\ref{eq:cyclic_shift}) 
implies the identity
$c^w_{u\: v(\lambda,k)} = c^z_{x\: v(\lambda,k)}$,
which completes the proof. 

If $wu^{-1}(n)= n$, then 
$zx^{-1} = 1\times wu^{-1}$, which is shape equivalent to 
$wu^{-1}$, and the result follows by Theorem~\ref{thm:B}~({\em ii}).
Assume $wu^{-1}(n)\neq n$. 
Then $w(k)=n$ and $u(k)<n$, as 
$w$ is Grassmannian with descent $k$.
Set $m:=u(k)$, $p:= u^{-1}(n) (>k)$ ,  and $l:=w(p)$.
Define $x\in {\cal S}_n$ by:
$$
x(j)\quad=\quad\left\{\begin{array}{lcl}
u(j)+1  &&  1\leq j< k\mbox{ or } p<j\\
1       &&  j=k\\
m+1     &&  j=k+1\\
u(j-1)+1&&  k+1<j\leq p   \end{array}\right..
$$
Then $x\leq _k z:= (wu^{-1})^{(1\,2\,\ldots\,n)}x$ where 
$$
z(j)\quad=\quad\left\{\begin{array}{lcl}
w(j)+1  &&  1\leq j<k\mbox{ or } p<j\\
l+1     &&  j=k\\
1       &&  j=k+1\\
w(j-1)+1&&  k+1<j\leq p   \end{array}\right..
$$

To show (5.3.1), %\ref{eq:cyclic_shift}), 
let $g_1(\alpha),\ldots,g_n(\alpha)$
for $\alpha\in A$ be the parameterized basis for flags 
$\Edot(\alpha)\in X^\circ_u\Fpdot\bigcap X^\circ_{w_0w}\Fdot$
constructed in the proof of Lemma~\ref{lem:skew_calculation}.
Since $g_k(\alpha)\in F'_{n+1-u(k)}\bigcap F_{w(k)}$,
$u(k)=m$, and $w(k)=n$,  
there exist regular functions $\beta_j(\alpha)$ on $A$ such that 
$$
g_k(\alpha)= e_n + \sum_{j=m}^{n-1} \beta_j(\alpha)e_j.
$$
Since 
$F'_1=\Span{e_n}\subset E_p(\alpha)- E_{p-1}(\alpha)$ and
$g_p(\alpha)=e_l$, there exist regular functions $\delta_j(\alpha)$ on
$A$ with $\delta_p(\alpha)$ nowhere vanishing such that 
\begin{eqnarray*}
e_n &=&
\sum_{j=1}^p \delta_j(\alpha)g_j(\alpha)\\
&=& g_k(\alpha)+ 
\sum_{j=1}^{k-1} \delta_j(\alpha)g_j(\alpha)
+\sum_{j=k+1}^p \delta_j(\alpha)e_{w(j)},
\end{eqnarray*}
as $g_k(\alpha)$ is the only vector among the $g_j(\alpha)$ in which 
$e_n$ has a non-zero coefficient.
Thus
$$
e_n- \sum_{j=k+1}^p \delta_j(\alpha) e_{w_j}\ =\ 
g_k(\alpha) + \sum_{j=1}^{k-1}\delta(\alpha)g_j(\alpha)
$$
is a vector in $E_k(\alpha) - E_{k-1}(\alpha)$.

Define a basis $h_1(\alpha),\ldots,h_n(\alpha)$ for $V$ by
$$
h_j(\alpha)\quad=\quad
\left\{\begin{array}{lcl}
g_j(\alpha) && 1\leq j< k \mbox{ or } p<j\\
e_n-\left( 
\sum_{j=k+1}^p \delta_j(\alpha)e_{w(j)}\right)
&& j=k\\
e_n && j=k+1\\
g_{j-1}(\alpha) && k+1<j\leq p 
\end{array}\right..
$$ 

We claim 
$\Epdot(\alpha):=\SPan{h_1(\alpha),\ldots,h_n(\alpha)}$
is a flag in $X_{w_0z}\Gdot\bigcap X_x\Gpdot$, which implies 
(5.3.1): %\ref{eq:cyclic_shift}):
Since 
$h_k(\alpha)\in E_k(\alpha) - E_{k-1}(\alpha)$
and $h_j(\alpha)=g_j(\alpha)$ for $j<k$, we have
$$
E'_k(\alpha)\quad=\quad \Span{h_1(\alpha),\ldots,h_k(\alpha)}
\quad=\quad E_k(\alpha).
$$
Thus if $\alpha\neq\alpha'$, then $\Epdot(\alpha)\neq\Epdot(\alpha')$ 
and so $\{\Epdot(\alpha)\,|\,\alpha\in A\}$ is a subset of 
the intersection 
$X_{w_0z}\Gdot\bigcap X_x\Gpdot$ of dimension equal to 
$\dim A=\ell(w)-\ell(u)=\ell(z)-\ell(x)$, the dimension of 
$X_{w_0z}\Gdot\bigcap X_x\Gpdot$.
Thus $\{\Epdot(\alpha)\,|\,\alpha\in A\}$ is dense, and so
$E'_k(\alpha)=E_k(\alpha)$ implies (5.3.1). %\ref{eq:cyclic_shift}).

For notational convenience, set $G^\circ_j:=G_j-G_{j-1}$, and similarly
for $F^\circ_j$.
To establish this claim, we first show that
$h_j(\alpha)\in G^\circ_{z(j)}$ for $j=1,\ldots,n$,
which shows $h_1(\alpha),\ldots,h_n(\alpha)$ is a parameterized
basis for $V$ 
and  $\Epdot(\alpha)\in X_{w_0z}\Gdot$.
Then, for a fixed $\alpha\in A$, we construct $h'_1,\ldots,h'_n$
which satisfy $\Epdot(\alpha)=\SPan{h'_1,\ldots,h'_n}$ and 
$h'_j\in G'_{n+1-x(j)}$ for $j=1,\ldots,n$, showing 
$\Epdot(\alpha)\in X_x\Gpdot$.

Note that if $i<n$, then $G_{i+1}=\Span{e_n,F_i}$.
Thus $h_j(\alpha)\in F^\circ_{w(j)}\subset G^\circ_{z(j)}$ for 
$1\leq j<k$  and $p<j$, and if $k+1<j\leq p$, then 
$h_j(\alpha)\in F^\circ_{w(j-1)}\subset G^\circ_{z(j)}$.
Then, since $G_1=\Span{e_n}$, we see that 
$h_{k+1}(\alpha)=e_n\in G'_1=G'_{n+1-x(k+1)}$.
Finally, since $w$ is Grassmannian of descent $k$, if $k+1\leq i\leq p$,
then $w(i)\leq w(p)=l$, which shows 
$h_k(\alpha)\in G^\circ_{l+1}=G^\circ_{z(k)}$.
Thus $\Epdot(\alpha)\in X^\circ_{w_0z}\Gdot$.

We now show that $\Epdot(\alpha)\in X_x\Gpdot$.
Note that if $a\leq b<n$, then 
$F'_{n+1-a}\bigcap F_b\subset G'_{n-a}\cap G_{b+1}$.
Thus if $1\leq j< k$, 
$h_j(\alpha)=g_j(\alpha)\in F'_{n+1-u(j)}\bigcap F_{w(j)}
\subset G'_{n+1-x(j)}$.
Since $x(k)=1$, we see that $h_k(\alpha)\in G'_{n+1-x(k)}=V$.
Fix $\alpha\in A$ and set $h'_j=h_j(\alpha)$ for $1\leq j\leq k$.
Define 
$$
h'_{k+1}\ :=\ g_k(\alpha) - e_n\ =\  
\sum_{j=m}^{n-1} \beta_j(\alpha)e_j\ \in\ 
G'_{n+1-(m+1)}\ =\  G'_{n+1-x(k+1)}.
$$
Since $h'_{k+1} + h_{k+1}(\alpha)= g_k(\alpha)$, we see that 
$E'_{k+1}(\alpha) = \Span{E'_k(\alpha),H'_{k+1}}$.

Finally, since $\Edot(\alpha)\in X_u\Fpdot$, 
if $k<j$ there exists a vector 
$$
g'_j\ :=\ \sum_{i\leq j}\gamma_{i,j}g_i(\alpha)
\ \in\  F'_{n+1-u(j)}
$$ 
such that $\Span{E_{j-1}(\alpha),g'_j}= E_j(\alpha)$.
For $k+1<j\leq p$, set 
$$
h'_j\ =\ g'_{j-1}-\gamma_{k,j-1}e_n\ \in\ 
\Span{e_{n-1}\ldots,e_{n+1-u(j-1)}}\ =\ G'_{n+1-x(j)},
$$
as as $g_k(\alpha)$ is the only vector among 
$\{g_1(\alpha),\ldots,g_n(\alpha)\}$ which is not in the span
of $e_1,\ldots,e_{n-1}$.
If $p<j$, set $h'_j=g'_j-\gamma_{k,j}e_n\in G'_{n+1-x(j)}$.
Then $\SPan{h'_1,\ldots,h'_n}=\Epdot(\alpha)$,
completing the proof.
\QEDnoskip

\section{Formulas for some Littlewood-Richardson
coefficients} 

\subsection{A chain-theoretic interpretation}\label{sec:chain_description}
We give a chain-theoretic interpretation for some Littlewood-Richardson
coefficients $c^\zeta_\lambda$ in terms chains in either the
$k$-Bruhat order or the $\preceq$-order, similar to the main results
of~\cite{sottile_pieri_schubert}. 
If either $u\lessdot_k (\alpha, \beta) u$  or 
$\zeta\precdot(\alpha,\beta)\zeta$
is a cover, label that edge in the Hasse diagram
with the integer $\beta=\max\{\alpha,\beta\}$.
Given a saturated chain in the $k$-Bruhat order from $u$ to
$\zeta u$, equivalently, a saturated $\preceq$-chain from $e$ to
$\zeta$, the {\em word} of that chain is its  sequence of edge labels. 
Given a word $\omega=a_1.a_2\ldots a_m$, Schensted
insertion~\cite{Schensted} or~\cite[\S 3.3]{Sagan} of $\omega$ into the
empty tableau 
gives a pair $(S,T)$ of Young tableaux, where $S$ is the {\em insertion
tableau} and $T$ the {\em recording tableau} of $\omega$.

Let $\mu\subset \lambda$ be partitions.
A permutation $\zeta$ is {\em shape-equivalent} to a skew
Young diagram $\lambda/\mu$ if there is a $k$ such that 
$\zeta$ is shape-equivalent to 
$v(\lambda,k)\cdot v(\mu,k)^{-1}$.
It follows that $\zeta$ is shape equivalent to some skew partition
$\lambda/\mu$  if and only if whenever 
$\alpha,\beta\in \mbox{\rm up}_\zeta$ or  
$\alpha,\beta\in \mbox{\rm down}_\zeta$,
$$
\alpha\ <\ \beta \quad
\Longleftrightarrow \quad \zeta(\alpha)\ <\ \zeta(\beta).
$$
We prove a stronger version of Theorem~\ref{thm:skew_shape}:

\begin{thm}\label{thm:skew_shape_prime}
Let $\mu\subset\lambda$ be partitions and 
suppose $\zeta\in {\cal S}_\infty$  is shape equivalent to
$\lambda/\mu$.
Then, for every partition $\nu$
\begin{enumerate}
\item[({\em i})] $c^\zeta_\nu = c^{\lambda/\mu}_\nu$, and 
\item[({\em ii})] For every standard Young tableau $T$ of shape $\nu$,
$$
c^\zeta_\nu\ =\ 
\#\left\{\begin{array}{cc}\mbox{$\preceq$-chains from $e $ to
$\zeta$ whose} \\  \mbox{word has 
recording tableau $T$} \end{array} \right\}
$$
Equivalently, if $u\leq_k w$ and $wu^{-1}=\zeta$, then 
$$
c^w_{u\,v(\nu,k)}\ =\ 
\#\left\{\begin{array}{cc}\mbox{Chains in $k$-Bruhat order from $u$ to} 
\\\mbox{$w$  whose word has 
recording tableau $T$} \end{array}\right\}
$$
\end{enumerate}
\end{thm}

\begin{rem}{\em
Theorem~\ref{thm:skew_shape_prime} ({\em ii}) gives a combinatorial
proof of Proposition~\ref{prop:chains}, when $wu^{-1}$ is shape
equivalent to a skew partition.
Theorem~\ref{thm:skew_shape_prime} ({\em ii}) is similar in form to
Theorem 8 of~\cite{sottile_pieri_schubert}: 
\medskip

\noindent{\bf Theorem 8~\cite{sottile_pieri_schubert}. }}
Suppose $\nu=(p,1^{q-1})$, a partition of `hook' shape.
Then for every $u,w\in {\cal S}_\infty$ and $k\in{\Bbb N}$, the constant
$c^w_{u\, v(\nu,k)}$ counts either set
\begin{enumerate}
\item[({\em i})] ${\displaystyle 
\left\{\begin{array}{c}\makebox[3.1in][c]{Chains in $k$-Bruhat order 
from $u$ to $w$ with}\\\makebox[3.1in][c]{\ word \  
$a_1<\cdots<a_p>a_{p+1}>\cdots>a_{p+q-1}$. }\end{array}\right\}
}$.
\item[({\em ii})]${\displaystyle 
\left\{\begin{array}{c}\makebox[3.1in][c]{Chains in $k$-Bruhat order 
from $u$ to $w$ with}\\\makebox[3.1in][c]{word \ 
$a_1>\cdots>a_q<a_{q+1}<\cdots<a_{p+q-1}$.} \end{array}\right\}
}$.
\end{enumerate}
\medskip

{\em The recording tableaux of words in ({\em i}) each have 
the integers $1,2,\ldots, p$ in the first row and
$1,p+1,\ldots,p{+}q{-}1$ in the first column.
Furthermore, these are the only words with this recording tableau.
Similarly, the  recording tableaux of  words in ({\em ii})
all have the integers $1,2,\ldots,q$ in the first column and 
$1,q{+}1,\ldots,p{+}q{-}1$ in the first row.
However,  Theorem~\ref{thm:skew_shape} is {\em not} a generalization
of this result:
The permutation $\zeta:=(143652)$ is not shape equivalent to any
skew partition as  $4,5\in \mbox{down}_\zeta$ but $\zeta(4)>\zeta(5)$.
Nevertheless, $c^\zeta_{(4,1)}=1$.
Interestingly, $\zeta$ satisfies the conclusions of
Theorem~\ref{thm:skew_shape}. 

While the hypothesis of Theorem~\ref{thm:skew_shape} is not
necessary for the conclusion to hold, some hypotheses are necessary:
Let $\zeta=(162)(354)$, a product of two disjoint 3-cycles.
Then $\zeta^{(1\,\ldots\,6)}=(132)(465)=v(\,%
\begin{picture}(12,6)   \put(0,0){\line(0,1){6}} 
\put(0,0){\line(1,0){12}}\put(3,0){\line(0,1){6}}
\put(0,3){\line(1,0){12}}\put(6,0){\line(0,1){6}}
\put(0,6){\line(1,0){6}}\put(9,0){\line(0,1){3}}
\put(12,0){\line(0,1){3}}\end{picture}\,,\, 2)\cdot
v(%
\begin{picture}(6,3)   \put(0,0){\line(0,1){3}} 
\put(0,0){\line(1,0){6}}\put(3,0){\line(0,1){3}}
\put(0,3){\line(1,0){6}}\put(6,0){\line(0,1){3}}
\end{picture}\,,\, 2)^{-1}$.
Hence, by Theorem~\ref{thm:D}, we have:
$$
c^\zeta_{\begin{picture}(12,3)   \put(0,0){\line(0,1){3}} 
\put(0,0){\line(1,0){12}}\put(3,0){\line(0,1){3}}
\put(0,3){\line(1,0){12}}\put(6,0){\line(0,1){3}}
\put(9,0){\line(0,1){3}}\put(12,0){\line(0,1){3}}
\end{picture}}
\quad=\quad
c^\zeta_{\,\begin{picture}(9,6)   \put(0,0){\line(0,1){6}} 
\put(0,0){\line(1,0){9}}\put(3,0){\line(0,1){6}}
\put(0,3){\line(1,0){9}}\put(6,0){\line(0,1){3}}
\put(0,6){\line(1,0){3}}\put(9,0){\line(0,1){3}}
\end{picture}}
\quad=\quad
c^\zeta_{\,\begin{picture}(6,6)   \put(0,0){\line(0,1){6}} 
\put(0,0){\line(1,0){6}}\put(3,0){\line(0,1){6}}
\put(0,3){\line(1,0){6}}\put(6,0){\line(0,1){6}}
\put(0,6){\line(1,0){6}}\end{picture}}
\quad=\quad 1.
$$
(This may also be seen as a consequence of Theorem~\ref{thm:C} and the
form of the Pieri-type formula
in~\cite{Lascoux_Schutzenberger_polynomes_schubert}, or 
of the main result, Theorem~5, of~\cite{sottile_pieri_schubert}.)
If $u=312645$, then $\zeta u=561234$ and the labeled Hasse diagram of
$[u,\zeta u]_2$ is:
$$
\epsfxsize=1.8in \epsfbox{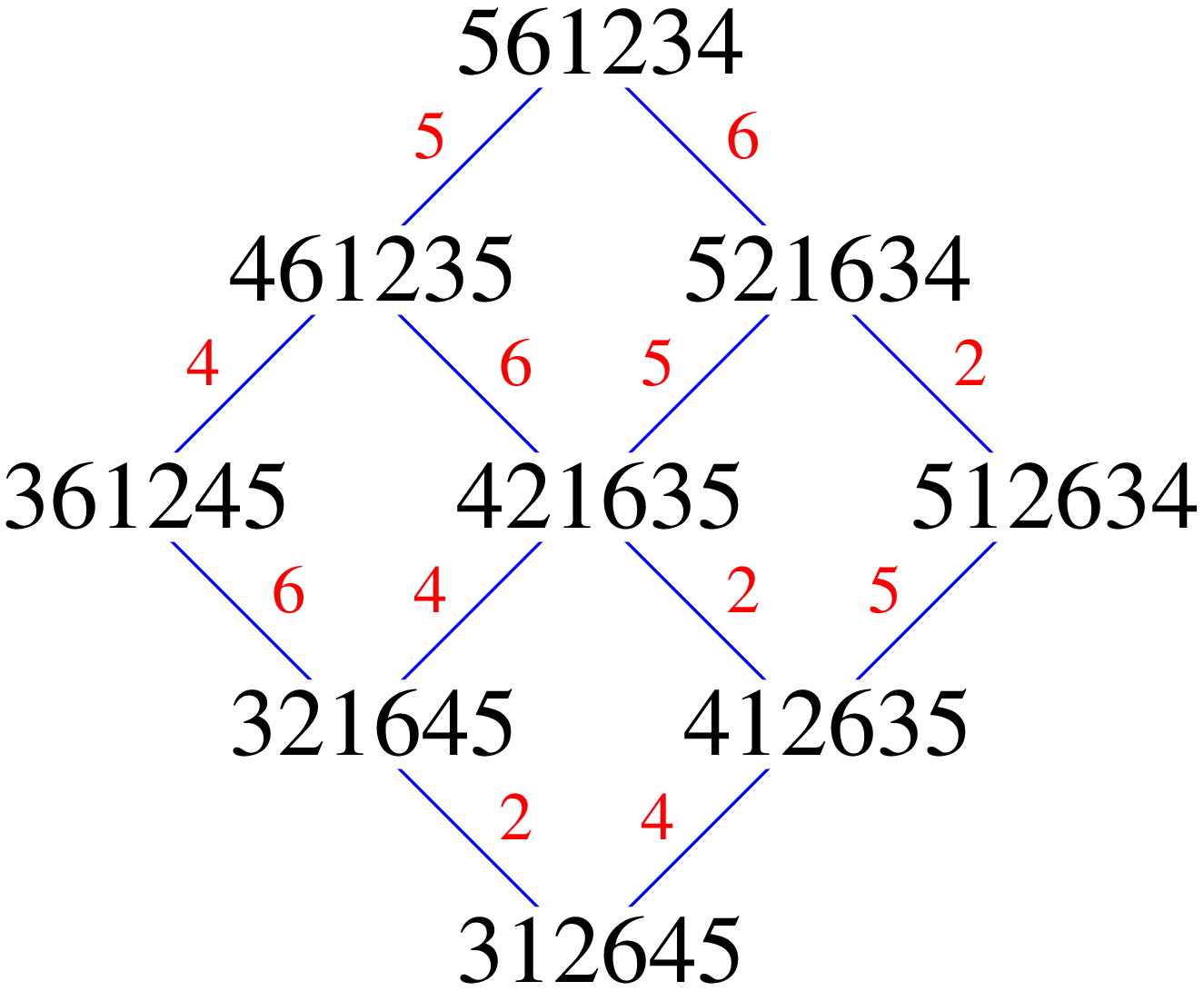}
$$
The labels of the six chains are:
$$
2456,\ 2465,\ 2645,\ 4526,\ 4256,\ 4265
$$
and these have (respective) recording tableaux:
$$
\epsfxsize=4.4in \epsfbox{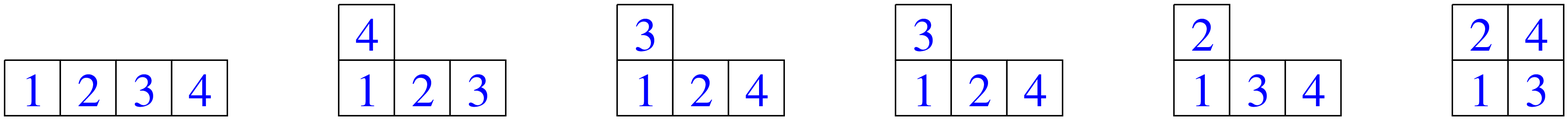}\ \raisebox{6pt}{.}
$$
This list omits the tableau
\begin{picture}(15,15)(0,2)   \put(0,0){\line(0,1){15}} 
\put(0,0){\line(1,0){15}}\put(7.5,0){\line(0,1){15}}
\put(0,7.5){\line(1,0){15}}\put(15,0){\line(0,1){15}}
\put(0,15){\line(1,0){15}}
\put(1.5,1.5){\tiny 1} \put(9,1.5){\tiny 2}
\put(1.5,9){\tiny 3} \put(9,9){\tiny 4}
\end{picture},
and the third and fourth tableaux are identical.
}
\end{rem}

\noindent{\bf Proof of Theorem~\ref{thm:skew_shape}. }
Suppose first that $\zeta=v(\lambda,k)\cdot v(\mu,k)^{-1}$.
Then 
$$
[e,\zeta]_\preceq\quad \simeq \quad  [v(\mu,k),\:v(\lambda,k)]_k
\quad  \simeq \quad  [\mu,\lambda]_\subset.
$$
The first isomorphism preserves the edge labeling of the Hasse diagrams, 
and in the second the labels of the $k$-Bruhat order correspond to 
diagonals in a Young diagram:
If $\nu\;\subsetdot\nu'$ is a cover in Young's lattice, there is a unique $i$
such that $\nu_i \neq\nu'_i$.
In that case $\nu_i+1=\nu'_i$ and
the label of the corresponding edge
in the $k$-Bruhat order is $k-i+\nu'_i$,
the diagonal on which the new box in $\nu'$ lies.

A chain in Young's lattice from $\mu$ to $\lambda$ is a standard skew
tableau $R$ of shape $\lambda/\mu$.
Consider the word, $a_1\ldots a_m$, of that chain as a two-rowed array:
$$
w\quad = \quad \left(\begin{array}{cccc}1&2&\cdots&m\\
a_1&a_2&\cdots&a_m\end{array}\right).
$$
Then the  entry $i$ of $R$ is in the
$a_i$th diagonal.

Let $S$ and $T$ be, respectively, the insertion and recording tableaux
for that two-rowed array.
Consider the two-rowed array consisting of the columns 
${a_i\choose i}$ arranged in lexicographic
order: that is, ${a_i\choose i}$ is to the left of 
${a_j\choose j}$ if either $a_i<a_j$ or $a_i=a_j$
and $i<j$.
Then the insertion and
recording tableaux of this new array are $T$ and $S$,
respectively~\cite{Schutzenberger_insertion,Knuth}.

The second row of this new array, the word
inserted to obtain $T$, is the `diagonal' word of the skew tableau $R$.
That is,  the entries of $R$ read lexicographically by diagonal.
By Lemma~\ref{lem:diagonal_word} (proven below), the diagonal word is
Knuth-equivalent to the original word.
Thus $T$ is the unique tableau of partition shape Knuth-equivalent to
$R$. 
This gives a combinatorial bijection
$$
\left\{\begin{array}{cc}\mbox{$\preceq$-chains from $e $ to
$\zeta$ whose} \\  \mbox{word has 
recording tableau $T$} \end{array} \right\}
\quad\Longleftrightarrow\quad
\left\{\begin{array}{cc}\mbox{Skew tableaux $R$ of shape}\\ 
\mbox{$\lambda/\mu$ Knuth-equivalent to $T$} \end{array} \right\},
$$
proving the theorem in this case, as it is well-known that (see,
for example~\cite[\S 4.9]{Sagan}), 
$$
c^{\lambda/\mu}_\nu \quad =\quad 
\#\left\{\begin{array}{cc}\mbox{Skew tableaux $R$ of shape}\\
\mbox{$\lambda/\mu$ Knuth-equivalent to $T$} \end{array} \right\}.
$$

Now suppose $\zeta$ is shape-equivalent to 
$v(\lambda,k)\cdot v(\mu,k)^{-1}$. 
By Theorem~\ref{thm:B}~({\em ii}), $c^\zeta_\nu=c^{\lambda/\mu}_\nu$,
proving ({\em i}). 
Assume $\lambda, \mu$, and $k$ have been chosen so that
$\zeta=\phi_P\left(v(\lambda,k)\cdot v(\mu,k)^{-1}\right)$, 
for some $P$.
By Theorem~\ref{thm:new_order}~({\em iii}), 
$\phi_P$ induces an isomorphism
$$
\phi_P\ :\ [e ,\,v(\lambda,k)\cdot v(\mu,k)^{-1}]_\preceq
\ \stackrel{\sim}{\longrightarrow}\ 
[e ,\,\zeta]_\preceq.
$$
Moreover, if $\eta\precdot(\alpha,\beta)\eta$ is a cover in 
$[e ,\,v(\lambda,k)\cdot v(\mu,k)^{-1}]_\preceq$, then 
$\phi_P\eta\precdot(\phi_P(\alpha\,\beta))\phi_P\eta$
is a cover in $[e,\zeta]_\preceq$ which has label $p_\beta$, where
$P=p_1<p_2<\cdots$.
Thus, if $\gamma$ is a chain in 
$[e,\,v(\lambda,k)\cdot v(\mu,k)^{-1}]_\preceq$
whose word $a_1,\ldots,a_m$ has recording tableau $T$, then 
$\phi_P(\gamma)$ is a chain in $[e ,\zeta]_\preceq$ 
with word $p_{a_1},\ldots,p_{a_m}$, which 
also has recording tableau $T$.
\QED

Order the diagonals of a skew Young tableau $R$ beginning with the
diagonal incident to the the end of the first column of $R$.
The {\em diagonal word} of $R$ is the entries of $R$ listed
in lexicographic order by diagonal, with magnitude breaking ties.
The tableau on the left below has diagonal word 
$7\,58\,379\,148\,26\,26\,5\,8$.
If we apply Schensted insertion to the initial segment $7\,58\,379\,148$,
(those diagonals incident upon the first column), we obtain the tableau
on the right, whose row word equals this initial segment.
$$
\epsfxsize=1.5in \epsfbox{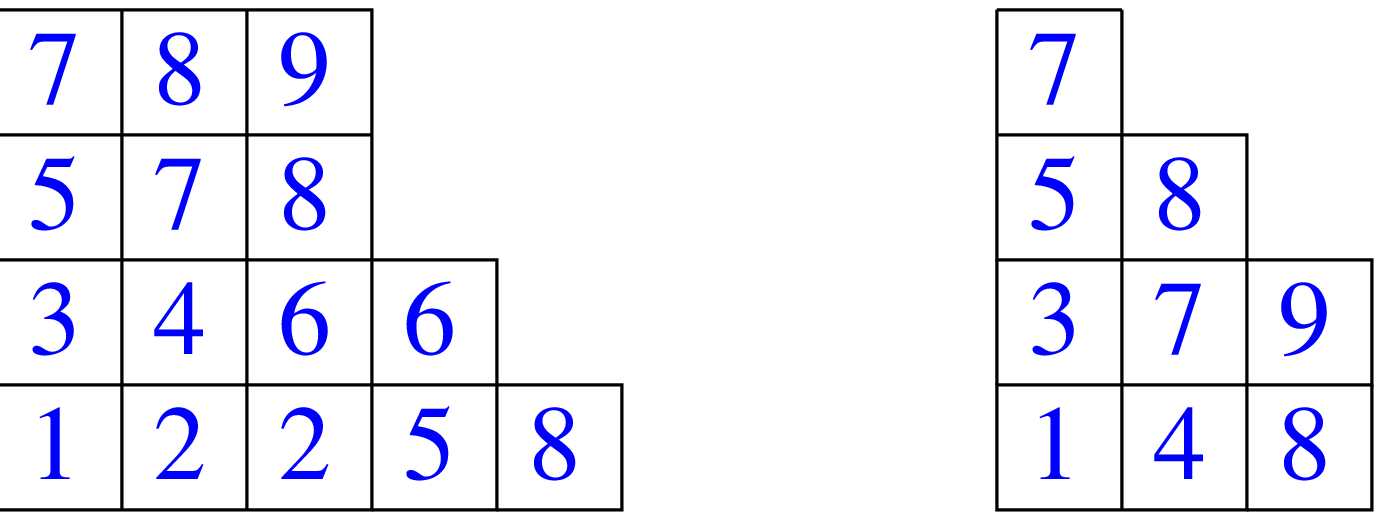}
$$
This observation is the key to the proof of the following lemma.

\begin{lem}\label{lem:diagonal_word}
The diagonal word of a skew tableau is Knuth-equivalent to its column
word. 
\end{lem}

\noindent{\bf Proof. }
For a skew tableau $R$, let $d(R)$ be its diagonal word.
We show that $d(R)$ is Knuth equivalent to the word $c.d(R')$, 
the concatenation of the first column $c$ of $R$ and $d(R')$, where $R'$
is obtained from $R$ by the removal of its first column.
An induction completes the proof.

Suppose the first column of $R$ has length $b$ and 
$R$ has $r$ diagonals.
For $1\leq j\leq b$ let $w_j:= a^j_1\ldots a^j_{s_j}$ be the subword of
$d(R)$ consisting of the $j$th diagonal.
Then $a^j_1<\cdots<a^j_{s_j}$, $s_1\leq s_2\leq \cdots\leq s_b$, and 
if $k\leq s_j$, then 
$a^j_k>a^{j+1}_k>\cdots>a^b_k$, as these are consecutive entries in the
$k$th column of $R$.

Consider the insertion tableau $T_l$ of the word
$w_1.w_2.\ldots w_l$ for $1\leq l\leq b$.
Then the $k$th column of $T_l$ is 
$a^j_k>\cdots>a^l_k$, where $s_{j-1}<k\leq s_j$.
Hence  $c.d(R') = c.\mbox{\em row}(T'). w_{b+1}\ldots  w_r$,
where $\mbox{\em row}(T')$ is the row word of the tableau obtained from
$T_b$ by removing its first column, which is $c$.
Since the column word of a tableau is Knuth-equivalent to its row word,
we are done.
\QEDnoskip

\subsection{Skew permutations}
Define the set of {\em skew permutations} to be the smallest set of
permutations containing all permutations 
$v(\lambda,k)\cdot v(\mu,k)^{-1}$ which is closed under:
\begin{enumerate}
\item[1.] Shape equivalence.  
If $\eta$ is shape equivalent to a skew
permutation $\zeta$, then $\eta$ is skew.
\item[2.] Cyclic shift.
If $\zeta\in{\cal S}_n$ is skew, then so is $\zeta^{(1\,2\,\ldots\,n)}$.
\item[3.] Products of disjoint permutations.
If $\zeta,\eta$ are disjoint and skew, then $\zeta\eta$ is skew. 
\end{enumerate}

A {\em shape} of a skew permutation $\zeta$ is a (non-unique!) skew
partition $\theta$ which is defined inductively.
If $\zeta$ is shape equivalent to $\lambda/\mu$, then 
$\zeta$ has shape $\lambda/\mu$.
If $\zeta\in {\cal S}_n$ is a skew permutation with shape $\theta$, then 
$\zeta^{(1\,2\,\ldots\,n)}$ has shape $\theta$.
If $\zeta$ and $\eta$ are disjoint skew permutations with respective
shapes $\rho$ and $\sigma$, then $\zeta\eta$ has skew shape
$\rho\coprod \sigma$.

\begin{thm}\label{thm:skew_permutation}
Let $\zeta$ be a skew permutation with shape $\theta$, then 
\begin{enumerate}
\item[({\em i})] For all partitions $\nu$,
$$
c^\zeta_\nu\quad=\quad c^\theta_\nu.
$$
\item[({\em ii})] The number of chains in the interval
$[e,\,\zeta]_\preceq$ is equal to  
the number of standard Young tableaux of shape $\theta$.
\end{enumerate}
\end{thm}

\noindent{\bf Proof. }
The number of standard skew tableaux of shape $\theta$ is
$\sum_\lambda f^\lambda c^\theta_\lambda$, 
hence ({\em ii}) is consequence of ({\em i}) and 
Proposition~\ref{prop:chains}.
To show ({\em i}), we need only consider the last part (3.) of the
recursive definition of skew permutations, by
Theorems~\ref{thm:B}~({\em ii}) and~\ref{thm:D}.
Suppose $\zeta$ and $\eta$ are disjoint skew permutations with
respective shapes 
$\rho$ and $\sigma$, and for all partitions $\nu$,
$c^\zeta_\nu=c^\rho_\nu$ and $c^\eta_\nu=c^\sigma_\nu$.
Then by Theorem~\ref{thm:C}~({\em ii}),
\begin{eqnarray*}
c^{\zeta\eta}_\nu &=&
\sum_{\lambda,\mu} c^\nu_{\lambda\,\mu}\; c^\zeta_\lambda\;c^\eta_\mu\\
&=&\sum_{\lambda,\mu} 
c^\nu_{\lambda\,\mu}\; c^\rho_\lambda\;c^\sigma_\mu\\
&=& c^{\rho \coprod \sigma}_\nu.
\makebox[.1in]{\hspace{1in}\QEDnoskip}
\end{eqnarray*}

\begin{ex}\label{example:disjoint}
{\em 
Consider the geometric graph of the permutation $(1978)(26354)$:
$$
\epsfxsize=1.1in \epsfbox{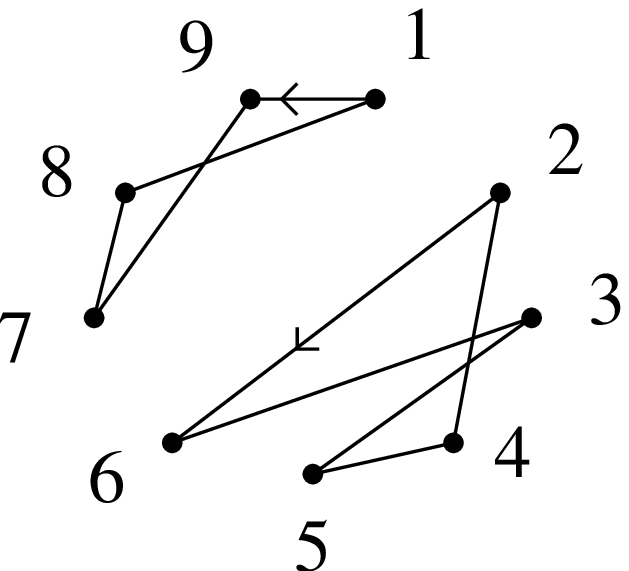}
$$
Thus the two cycles $\zeta=(1978)$ 
and $\eta=(26354)$ are disjoint.

Note that $\zeta$ is shape equivalent to $(1423)$ and 
$(1423)^{(1234)} = (1342)$.
Similarly, $\eta$ is shape equivalent to $(15243)$ and 
$(15243)^{(12345)} = (13542)$.
Both of these cycles, $(1423)$ and  $(15243)$, are skew partitions:
Let $\lambda=\;$%
\setlength{\unitlength}{1.3pt}\begin{picture}(3,6)
\put(0,0){\line(1,0){3}}\put(3,0){\line(0,1){3}}
\put(3,3){\line(-1,0){3}}\put(0,3){\line(0,-1){3}}
\end{picture}\,, 
$\mu=\;$%
\begin{picture}(6,6)
\put(0,0){\line(1,0){6}}\put(6,0){\line(0,1){6}}
\put(6,6){\line(-1,0){6}}\put(0,6){\line(0,-1){6}}
\put(3,0){\line(0,1){6}}\put(0,3){\line(1,0){6}}
\end{picture}\,, 
$\nu=\;$%
\begin{picture}(9,6)
\put(0,0){\line(0,1){6}}\put(0,0){\line(1,0){9}}
\put(3,0){\line(0,1){6}}\put(0,6){\line(1,0){6}}
\put(6,0){\line(0,1){6}}\put(0,3){\line(1,0){9}}
\put(9,0){\line(0,1){3}}
\end{picture}\,.
Then
$$
v(\lambda,2)\quad=\quad 13245,\qquad
v(\mu,2)\quad=\quad 34125,\qquad
v(\nu,2)\quad=\quad 35124
$$
and
$$
\begin{array}{rcccl}
v(\lambda,2)&\leq_2&(1342)\cdot v(\lambda,2)&=& v(\mu,2),\\
v(\lambda,2)&\leq_2&(13542)\cdot v(\lambda,2)&=& v(\nu,2).\\
\end{array}
$$

Hence, for every partition $\kappa$,
$c^\zeta_\kappa=c^{\mu/\lambda}_\kappa$ and 
$c^\eta_\kappa=c^{\nu/\lambda}_\kappa$.
Thus it follows that 
$c^{\zeta\eta}_\kappa = c^\rho_\kappa$, where $\rho$ is any of the
four skew partitions:
$$
\setlength{\unitlength}{.75pt}\begin{picture}(50,40)\thicklines
\put(30, 0){\line(1,0){20}}
\put(20,10){\line(1,0){30}}
\put( 0,20){\line(1,0){40}}
\put( 0,30){\line(1,0){20}}
\put( 0,40){\line(1,0){10}}
\put( 0,20){\line(0,1){20}}
\put(10,20){\line(0,1){20}}
\put(20,10){\line(0,1){20}}
\put(30, 0){\line(0,1){20}}
\put(40, 0){\line(0,1){20}}
\put(50, 0){\line(0,1){10}}
\end{picture}\qquad
\begin{picture}(50,40)\thicklines
\put(30, 0){\line(1,0){20}}
\put(20,10){\line(1,0){30}}
\put(10,20){\line(1,0){30}}
\put( 0,30){\line(1,0){20}}
\put( 0,40){\line(1,0){20}}
\put( 0,30){\line(0,1){10}}
\put(10,20){\line(0,1){20}}
\put(20,10){\line(0,1){30}}
\put(30, 0){\line(0,1){20}}
\put(40, 0){\line(0,1){20}}
\put(50, 0){\line(0,1){10}}
\end{picture}\qquad
\begin{picture}(50,40)\thicklines
\put(40, 0){\line(1,0){10}}
\put(30,10){\line(1,0){20}}
\put(10,20){\line(1,0){40}}
\put( 0,30){\line(1,0){30}}
\put( 0,40){\line(1,0){20}}
\put( 0,30){\line(0,1){10}}
\put(10,20){\line(0,1){20}}
\put(20,20){\line(0,1){20}}
\put(30,10){\line(0,1){20}}
\put(40, 0){\line(0,1){20}}
\put(50, 0){\line(0,1){20}}
\end{picture}\qquad
\begin{picture}(50,40)\thicklines
\put(30, 0){\line(1,0){20}}
\put(30,10){\line(1,0){20}}
\put(10,20){\line(1,0){30}}
\put( 0,30){\line(1,0){30}}
\put( 0,40){\line(1,0){20}}
\put( 0,30){\line(0,1){10}}
\put(10,20){\line(0,1){20}}
\put(20,20){\line(0,1){20}}
\put(30, 0){\line(0,1){30}}
\put(40, 0){\line(0,1){20}}
\put(50, 0){\line(0,1){10}}
\end{picture}
$$

}\end{ex}

\subsection{Further remarks}\label{sec:further}
For small symmetric groups, it is instructive to examine all permutations
and determine to which class they belong.
Here, we enumerate each class in ${\mathcal S}_4$, 
${\mathcal S}_5$, and ${\mathcal S}_6$:

\begin{center}
\begin{tabular}{|l||c|c|c|}
\hline
&\begin{tabular}{c}skew\\ partitions\end{tabular}& 
\begin{tabular}{c}shape equivalent to\\ a skew partition\end{tabular} & 
\begin{tabular}{c}skew\\ permutation\end{tabular}\\
\hline
\hline
${\cal S}_4$ & 14 & 21 & 24\\
\hline
${\cal S}_5$ & 42 & 79 & 120 \\
\hline
${\cal S}_6$ &  132 & 311 & 678\\
\hline
\end{tabular}
\end{center}

If $\zeta$ is one of the 42 permutations in ${\cal S}_6$ which are not
skew permutations, and $\zeta$ is not among
$$
(125634),\  (145236),\  (143652),\  (163254),\ (153)(246),\ \mbox{or}\ 
(135)(264), 
$$
then there is a skew partition $\theta$ such that 
$c^\zeta_\nu = c^\theta_\nu$ for all partitions $\nu$.
It would be interesting to understand why this occurs for all but these 6 
permutations.
Is there a wider class of permutations $\zeta$ such that there exists a
skew partition $\theta$ with $c^\zeta_\nu = c^\theta_\nu$ for all
partitions $\nu$? 

For the six `exceptional' permutations $\zeta$, there is a skew partition
$\theta$ for which $c^\zeta_\nu = c^\theta_\nu$ for all 
$\nu\subset a^b$, where $a=\#$up$_\zeta$ and $b=\#$down$_\zeta$.
For these, $\theta\not\subset a^b$.
For example, let $\zeta=(153)(246)$.
If $u=214365$, then $u\leq_3 \zeta u$ and there are 42 chains in 
$[u,\,\zeta u]_3$.
Also
$$
c^\zeta_{\,\begin{picture}(6,4) \put(0,0){\line(0,1){4}}
\put(0,0){\line(1,0){6}}      \put(2,0){\line(0,1){4}}
\put(0,2){\line(1,0){6}}      \put(4,0){\line(0,1){4}}
\put(0,4){\line(1,0){6}}      \put(6,0){\line(0,1){4}}
\end{picture}}\ =\ 1, \qquad
c^\zeta_{\,\begin{picture}(6,6) \put(0,0){\line(0,1){6}}
\put(0,0){\line(1,0){6}}      \put(2,0){\line(0,1){6}}
\put(0,2){\line(1,0){6}}      \put(4,0){\line(0,1){4}}
\put(0,4){\line(1,0){4}}      \put(6,0){\line(0,1){2}}
\put(0,6){\line(1,0){2}} \end{picture}}\ =\ 2, \qquad\mbox{and}\qquad
c^\zeta_{\,\begin{picture}(4,6)(0,0) \put(0,0){\line(1,0){4}}
\put(0,0){\line(0,1){6}}      \put(0,2){\line(1,0){4}}
\put(2,0){\line(0,1){6}}      \put(0,4){\line(1,0){4}}
\put(4,0){\line(0,1){6}}      \put(0,6){\line(1,0){4}}
\end{picture}}\ =\ 1,
$$
which verifies Proposition~\ref{prop:chains} as
$f^{\,\begin{picture}(6,4) \put(0,0){\line(0,1){4}}
\put(0,0){\line(1,0){6}}      \put(2,0){\line(0,1){4}}
\put(0,2){\line(1,0){6}}      \put(4,0){\line(0,1){4}}
\put(0,4){\line(1,0){6}}      \put(6,0){\line(0,1){4}}
\end{picture}}=5$, $f^{\,\begin{picture}(6,6) \put(0,0){\line(0,1){6}}
\put(0,0){\line(1,0){6}}      \put(2,0){\line(0,1){6}}
\put(0,2){\line(1,0){6}}      \put(4,0){\line(0,1){4}}
\put(0,4){\line(1,0){4}}      \put(6,0){\line(0,1){2}}
\put(0,6){\line(1,0){2}} \end{picture}}=16$, 
and $f^{\,\begin{picture}(6,4) \put(0,0){\line(0,1){4}}
\put(0,0){\line(1,0){6}}      \put(2,0){\line(0,1){4}}
\put(0,2){\line(1,0){6}}      \put(4,0){\line(0,1){4}}
\put(0,4){\line(1,0){6}}      \put(6,0){\line(0,1){4}}
\end{picture}}=5$.
In this case, $\theta = \begin{picture}(12,12) 
\put(6,0){\line(1,0){6}}      \put(0,6){\line(0,1){6}}
\put(6,3){\line(1,0){6}}      \put(3,6){\line(0,1){6}}
\put(0,6){\line(1,0){9}}      \put(6,0){\line(0,1){9}}
\put(0,9){\line(1,0){6}}      \put(9,0){\line(0,1){6}}
\put(0,12){\line(1,0){3}}     \put(12,0){\line(0,1){3}}
\end{picture}\,$.
Since ${\rm up}_\zeta=\{1,2,4\}$ and ${\rm down}_\zeta=\{6,5,3\}$, we
see that $a=b=3$, however
$\theta\not\subset  \begin{picture}(9,9) 
\put(0,0){\line(1,0){9}}      \put(0,0){\line(0,1){9}}
\put(0,3){\line(1,0){9}}      \put(3,0){\line(0,1){9}}
\put(0,6){\line(1,0){9}}      \put(6,0){\line(0,1){9}}
\put(0,9){\line(1,0){9}}      \put(9,0){\line(0,1){9}}
\end{picture}\, = a^b$.

A combinatorial interpretation of the Littlewood-Richardson coefficients
$c^w_{u\,v(\lambda,k)}$  should also give a bijective proof of
Proposition~\ref{prop:chains}.
We show a partial converse to this, that a function $\tau$ from chains
to standard Young tableaux satisfying some extra conditions will provide
a combinatorial interpretation of the Littlewood-Richardson coefficients
$c^w_{u\,v(\lambda,k)}$.

Let $\Ch[u,w]_k$ denote the set of (saturated) chains in the interval
$[u,w]_k$. 
For a partition $\mu$ and integer $m$, let $\mu*m$ be the set of
partitions $\lambda$ with $\lambda-\mu$ a horizontal strip of length
$m$.
These partitions arise in the classical Pieri's formula:
$$
S_\mu(x_1,\ldots,x_k)\cdot h_m(x_1,\ldots,x_k)
\ =\ \sum_{\lambda\in \mu*m} S_\lambda(x_1,\ldots,x_k).
$$
If $T$ is a standard tableau of shape $\mu$ and $m$ and integer,
let $T*m$ be the set of tableaux $U$ which contain $T$ as an initial
segment such that  $U-T$ is a horizontal strip whose entries increase
from left to right.

\begin{thm}\label{thm:combinatorial}
Suppose that for every $u\leq_k w$, there is a map
\begin{eqnarray*}
\Ch[u,w]_k& \longrightarrow&
\left\{
\mbox{\begin{minipage}[c]{2.5in}
\begin{center} Standard  Young tableau $T$ whose shape is a partition of
$\ell(w)-\ell(u)$\end{center}\end{minipage}}\right\}\\
\gamma&\longmapsto& \tau(\gamma)
\end{eqnarray*}
such that 
\begin{enumerate}

\item 
$d^w_{u\:v(\lambda,k)}:=
\#\{\gamma\in\mbox{ch}[u,w]_k\:|\: \tau(\gamma)=T\}$ depends only
upon the shape $\lambda$ of the standard tableau $T$.

\item 
If $\gamma = \delta.\varepsilon$ is the concatenation of two chains
$\delta$ and $\varepsilon$, then 
$\tau(\delta)$ is a subtableau of $\tau(\gamma)$.
(This means that $\tau(\gamma)$ is a recording tableau.)

\item
Suppose $\gamma = \delta.\varepsilon$ 
with $\delta\in\Ch[u,x]_k$, and hence 
$\varepsilon\in\Ch[x,w]_k$.
Then 
$\tau(\delta.\varepsilon)\in \tau(\delta)*m$ only if $x\rkm w$,
and $\varepsilon(\delta):=\varepsilon\in\Ch[x,w]_k$ is unique for this 
to occur.
\end{enumerate}

Then, for every standard tableau $T$ of shape $\lambda$ and $u\leq_k w$,
$$
c^w_{u\;v(\lambda,k)}\ =\ d^w_{u\:v(\lambda,k)}.
$$
\end{thm}

Such a map $\tau$ is a generalization of Schensted insertion.
In that respect, the existence of such a map would generalize
Theorem~\ref{thm:skew_shape_prime}. 
\medskip

\noindent{\bf Proof. }
We induct on $\lambda$.
Assume the theorem holds for all $u,w$, and partitions $\pi$
either with fewer rows than $\lambda$, or if $\lambda$ and $\pi$
have the same number of rows, then the last row of $\pi$ is shorter 
than the last row of $\lambda$.

The form of the Pieri-type formulas expressed
in~\cite{sottile_pieri_schubert,Winkel_multiplication}  
(also \S\ref{sec:emdedd})
and condition (3) prove the theorem when 
$\lambda$ consists of a single row.
Assume that $\lambda$ has more than one row and set $\mu$ to be
$\lambda$ with its last row removed.
Let $m$ be the length of the last row of $\lambda$
and $T$ be any tableau of shape $\mu$.
Recall that  $U\mapsto \mbox{shape}(U)$ gives 
a one-to-one correspondence between $T*m$
and $\mu*m$.

By the definition of $c^y_{u\; v(\mu,k)}$, we have
$$
{\frak S}_u\cdot S_\mu(x_1,\ldots,x_k)\ =\ 
\sum_{u\leq_k y}c^y_{u\; v(\mu,k)}\: {\frak S}_y.
$$
By the Pieri formula for Schubert polynomials, 
$$
{\frak S}_u\cdot S_\mu(x_1,\ldots,x_k)\cdot h_m(x_1,\ldots,x_k)\ =\  
\sum_w\left(
\sum_{\stackrel{\mbox{\scriptsize $u\leq_k y$}}{y\rrkm w}}
c^y_{u\; v(\mu,k)}\right) {\frak S}_w.
$$
By the classical Pieri formula, this also equals
$$
{\frak S}_u\cdot \sum_{\pi\in \mu*m}S_{\pi}(x_1,\ldots,x_k)
\ =\  \sum_w\left(\sum_{\pi\in \mu*m}
c^w_{u\; v(\pi,k)}\right) {\frak S}_w.
$$
Hence 
$$
\sum_{\pi\in \mu*m}c^w_{u\; v(\pi,k)}\ =\ 
\sum_{\stackrel{\mbox{\scriptsize $u\leq_k y$}}{y\rrkm w}}
c^y_{u\; v(\mu,k)}.
$$

We exhibit a bijection between the two sets
$$
M_{T,{k,m}}\ :=\ 
\coprod_{\stackrel{\mbox{\scriptsize $u\leq_k y$}}{y\rrkm w}}
\{ \delta\in\Ch[u,y]_k\: |\: \tau(\delta)=T\}
$$
and $\coprod_{\pi\in\mu*m}L_{\pi}$, where
$$
L_{\pi}\ :=\ 
\{ \gamma\in\Ch[u,w]_k\: |\; \tau(\gamma)\in T*m \mbox{ and 
$\tau(\gamma)$ has shape $\pi$}\}.
$$

This will complete the proof.
Indeed, by the induction hypothesis
$$
\#M_{T,{k,m}}\ =\  
\sum_{\stackrel{\mbox{\scriptsize $u\leq_k y$}}{y\rrkm w}}
c^y_{u\:v(\mu,k)}
$$
and for $\pi\in \mu*m$ with $\pi\neq \lambda$, 
$$
\#L_{\pi}\ = \ c^w_{u\:v(\pi,k)}.
$$
Thus the bijection shows
$$
c^w_{u\:v(\lambda,k)}\ =\ 
\sum_{\stackrel{\mbox{\scriptsize $u\leq_k y$}}{y\rrkm w}}
c^y_{u\:v(\mu,k)} - \sum_{\pi\in\mu*m,\  \pi\neq \lambda}
c^w_{u\:v(\pi,k)}\ =\ \#L_{\lambda},
$$
which is $\#\tau^{-1}(U)$, for any $U$ of shape $\lambda$.

To construct the desired bijection, consider first the map
$$
M_{T,{k,m}}\ \longrightarrow\ \coprod_{\pi\in\mu*m} L_\pi
$$
defined by $\delta\in\Ch[u,y]_k \mapsto \delta.\varepsilon(\delta)$.
By property 3, $\tau(\delta.\varepsilon(\delta))\in T*m$, so this
injective map has the stated range. 
To see it is surjective, let $\pi\in \mu*m$ and $\gamma\in L_\pi$.
Let $\delta$ be the first $|\mu|$ steps in the chain $\gamma$, so that 
$\gamma = \delta.\varepsilon$ and suppose $\delta\in\Ch[u,y]_k$.
Then $\tau(\delta)=T$ so $\tau(\delta.\varepsilon)\in \tau(\delta)*m$.
By 3, this implies $y\rkm w$, and hence 
$\delta\in M_{T,{k,m}}$.
\QEDnoskip

\appendix

\section{Illustrating the geometric theorems}

These appendices are intended for informal distribution
with this manu\-script and will not appear in the published version.
They contain no results, only examples which we hope may illustrate some
of the main results of this manuscript.
This appendix is intended to illustrate the geometric
results in the previous sections, particularly of
Section~\ref{sec:geometry}. 
We hope this may help others think about Schubert varieties and
intersections of Schubert varieties.

Throughout, let $e_1,\ldots,e_n$ be a fixed, ordered basis
for the vector space ${\Bbb C}^n$.
We use this basis to obtain a parameterization
for Schubert cells and their intersections.
Flags are represented by $n\times n$
matrices $M$: 
Let $(g_1,\ldots,g_n):=M\cdot e^T$ be the ordered basis given
by the `change of basis' matrix $M$.
The $i$th row of $M$ gives the components of $g_i$.
Then $M$ represents the flag $\SPan{g_1,\ldots,g_n}$.
We adopt some conventions for the entries of $M$:
a dot (${\setlength{\unitlength}{1pt}\begin{picture}(5,2)
\put(2.5,3){\circle*{2}}\end{picture}}$) will denote an entry of zero
and an asterix
(\raisebox{-4pt}{\epsfxsize=.1in \epsfbox{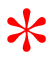}}) an
entry which may assume any value in ${\Bbb C}$.
One last convention is that the flags $\Edot, \Fdot$, 
etc.~will always be defined to be $\Edot:=\SPan{e_1,\ldots,e_n}$
and the `primed' flags $\Epdot, \Fpdot$, etc.~which are opposite to
their unprimed cousins will be defined by
$\Epdot:=\SPan{e_n,e_{n-1},\ldots,e_2,e_1}$.
We refer to these as the {\em standard flags}.

\subsection{Theorem~\ref{thm:B}~({\em ii})}

In Theorem~\ref{thm:B}~({\em ii}), we 
had $u\leq_k w$, $x\leq_k z$, and $w u^{-1}=z x^{-1}$ and we
studied $X_{w_0 w}\Edot\bigcap X_u\Epdot$
and $X_{w_0z}\Edot\bigcap X_x\Epdot$.
The main result was that, in $\Gr_k {\Bbb C}^n$,
$$
\pi_k\left(X_{w_0 w}\Edot\bigcap X_u\Epdot\right)\ 
=\ \pi_k\left(X_{w_0z}\Edot\bigcap X_x\Epdot\right).
$$

The general case of Theorem~~\ref{thm:B}~({\em ii}) was reduced to 
Lemma~\ref{lem:skew_critical}, where $w$ was Grassmannian of descent
$k$, and $k<i\Longrightarrow u(i)=x(i)$ (and hence also $w(i)=z(i)$).
The first example illustrates this case.

Let $n=7$, $k=4$, and 
$$
\begin{array}{r}u\ =\ 1436257\\w\ =\ 4567123\end{array}\qquad
\begin{array}{r}x\ =\ 4631257\\z\ =\ 5764123\end{array}
$$
The following matrices respectively represent general flags in the
Schubert cells  $X^\circ_{w_0 w}\Edot$, $X^\circ_u\Epdot$,
$X^\circ_{w_0z}\Edot$,  and 
$X^\circ_x\Epdot$:
$$
\epsfxsize=4.8in \epsfbox{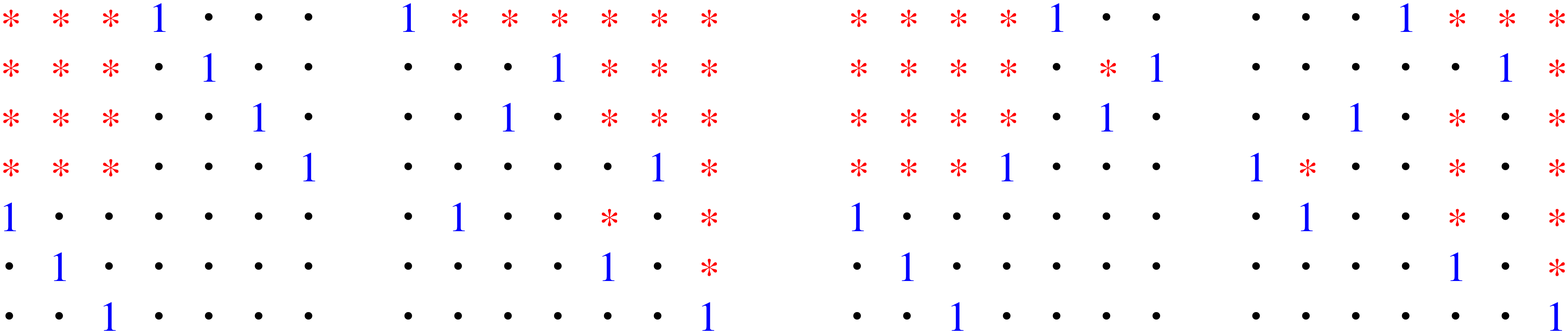}
$$
$w$ is chosen to be Grassmannian so that the cell
$X^\circ_{w_0 w}\Edot$ has a particularly simple form,
which gave an easy parameterization for the intersection of
the two cells, $X^\circ_{w_0 w}\Edot\bigcap X^\circ_u\Epdot$.
In the proof of Lemma~\ref{lem:skew_calculation} we describe how to find
bases parameterized by $A:=\{M\in M(w)\;|\; M\in X^\circ_u\Edot\}$.
In practice, this method may be used to determine the subvariety $A$ of
$M(w)$.

First, let $g_1,\ldots,g_7$ be the rows of the following matrix,
where $\alpha,\beta,\gamma,\delta,x,\rho,\sigma$, and $\tau$ are
arbitrary elements of ${\Bbb C}$ with $\alpha\delta x\tau\neq 0$:
$$
\epsfxsize=.9in \epsfbox{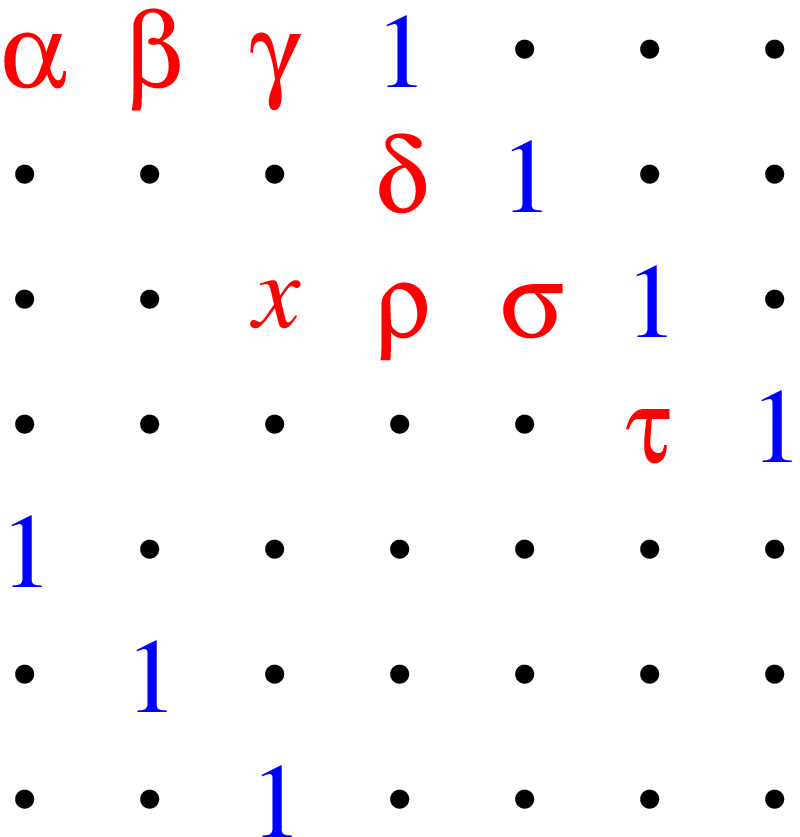}
$$
These parameters were chosen so that  for each $j=1,2,3,4$, 
$g_j\in E_{w(j)}\bigcap E'_{n+1-u(j)}$ and 
does not lie in either of $E_{w(j)-1}$ or
$E'_{n-u(j)}$, hence the 1's, the condition on 
$\alpha,\delta,x,\tau$, and the 0's
(${\setlength{\unitlength}{1pt}\begin{picture}(5,2)
\put(2.5,3){\circle*{2}}\end{picture}}$) in their initial columns.

This matrix determines a flag 
$\Gdot:=\SPan{g_1,\ldots,g_n} \in X^\circ_{w_0w}\Edot$, since it
is in $M(w)$.
Also, since $g_j\in E'_{n+1-u(j)}-E'_{n-u(j)}$ for $j\leq k$, at least
$G_1,\ldots,G_k$ satisfy the conditions for the flag $\Gdot$ to be in
$X^\circ_u\Epdot$. 
The remaining conditions for $\Gdot\in X^\circ_u\Epdot$, 
$$
G_{j-1}\bigcap E'_{n+1-u(j)}\subsetneq G_j\bigcap E'_{n+1-u(j)}
\qquad \mbox{for}\qquad k<j, 
$$ 
impose additional restrictions on the parameters.
In practice this means we seek conditions to ensure that 
$({\Bbb C}^\times g_j + G_{j-1})\bigcap E'_{n+1-u(j)}$ 
is non-empty.
For instance, for $j=6$, since 
$\Span{g_5,g_5}=
(\raisebox{-4pt}{\epsfxsize=.1in \epsfbox{figures/asterix.eps}}, 
\raisebox{-4pt}{\epsfxsize=.1in
\epsfbox{figures/asterix.eps}},0,0,0,0,0)$
and $E_{n+1-u(6)} = (0,0,0,0,
\raisebox{-4pt}{\epsfxsize=.1in\epsfbox{figures/asterix.eps}}, 
\raisebox{-4pt}{\epsfxsize=.1in\epsfbox{figures/asterix.eps}},
\raisebox{-4pt}{\epsfxsize=.1in\epsfbox{figures/asterix.eps}})$, 
some cancellation must occur.
Indeed, since
\begin{eqnarray*}
-\alpha g_5 -\beta g_6 + g+1 &=&
(0,0,\gamma,1,0,0,0)\\
g_3&=& (0,0,x,\rho,\sigma,1,0),
\end{eqnarray*}
we must have $\gamma\rho-x =0$ in order that 
$({\Bbb C}^\times g_6 + G_5)\bigcap E'_{n+1-u(6)} \neq \emptyset$ .
In the general situation, more complicated determinantal conditions
may arise.
\bigskip

From these considerations, we arrive at a parameterization
for $X_{w_0w}\Edot\bigcap X_u\Epdot$.
$$
\epsfxsize=2.2in \epsfbox{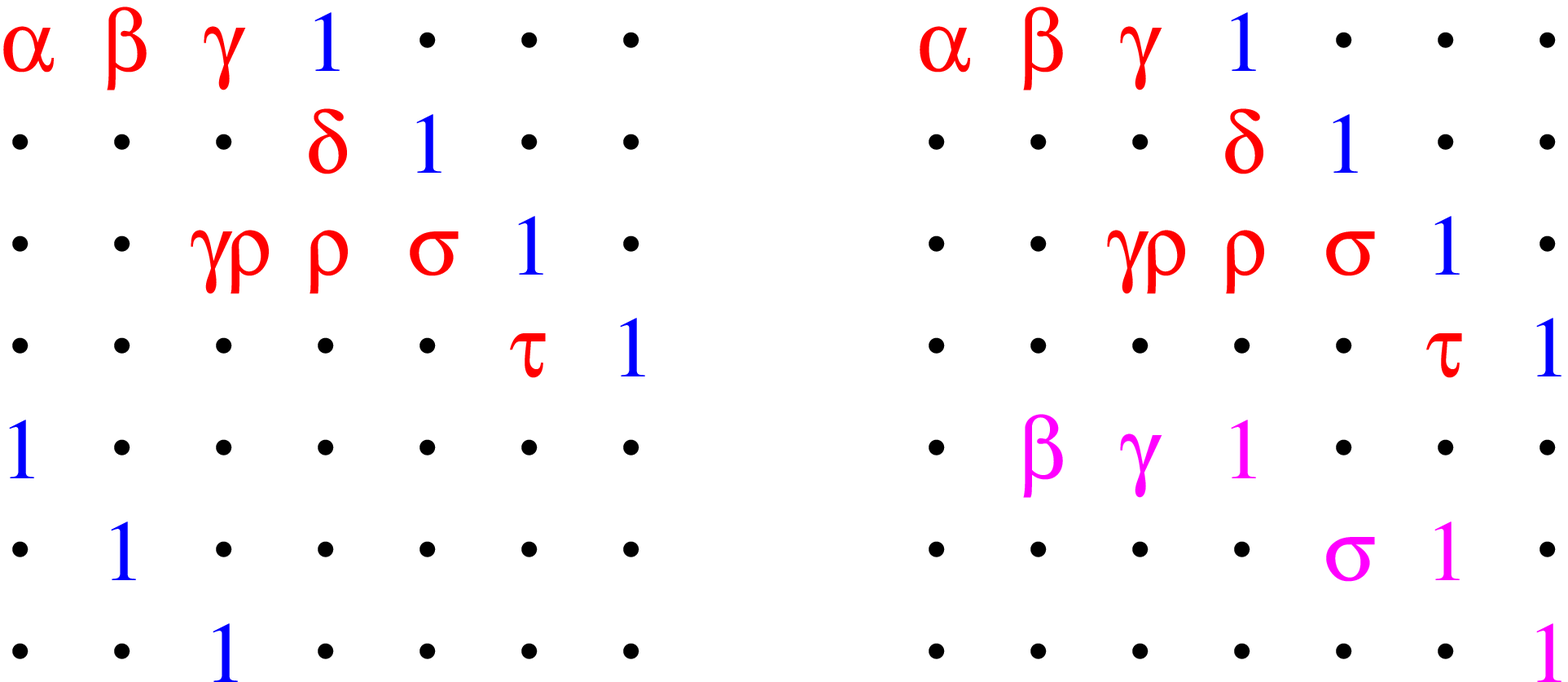}
$$
For $\alpha,\ldots,\tau \in {\Bbb C}^\times$,
both matrices represent the same flag in the intersection.
To see this,  
let $g_1,\ldots,g_7$ be the basis determined by the first matrix, and
$g'_1,\ldots,g'_7$ the basis determined by the second matrix.
Then, by the definition~(\ref{eq:Schubert_cell_definition}) of Schubert
cells in \S\ref{sec:flag}, the flags 
$\SPan{g_1,\ldots,g_7}\in X_{w_0w}^\circ\Edot$ and
$\SPan{g'_1,\ldots,g'_7}\in X^\circ_u\Epdot$.
Since $g_i=g'_i$ for $i=1,2,3,4$,
\begin{eqnarray*}
g'_5 &=& g_1 - \alpha g_5\\
g'_6 &=& g_3 - \rho(g_1-\alpha g_5-\beta g_6)\\
g'_7 &=& g_4 - \tau\left[g_3- \rho(g_1-\alpha g_5-\beta g_6)
-\sigma\left(g_2-\delta(g_1-\alpha g_5-\beta g_6-\gamma
g_7)\right)\right],
\end{eqnarray*}
showing $\SPan{g_1,\ldots,g_7}=\SPan{g'_1,\ldots,g'_7}$.
Lastly, since $\ell(w)-\ell(u)=12 -5=7$, and $\Edot,\Epdot$ are opposite
flags, we see that $X^\circ_u\Epdot\bigcap X^\circ_{w_0w}\Edot$ is
irreducible of dimension 7.
This identifies a 7-parameter family of flags in this
intersection, which must be dense.

Similarly, (with the same restrictions on parameters), the two
matrices below both represent the same flag in 
$X^\circ_{w_0z}\Edot\bigcap X^\circ_x\Epdot$:
$$
\epsfxsize=2.2in \epsfbox{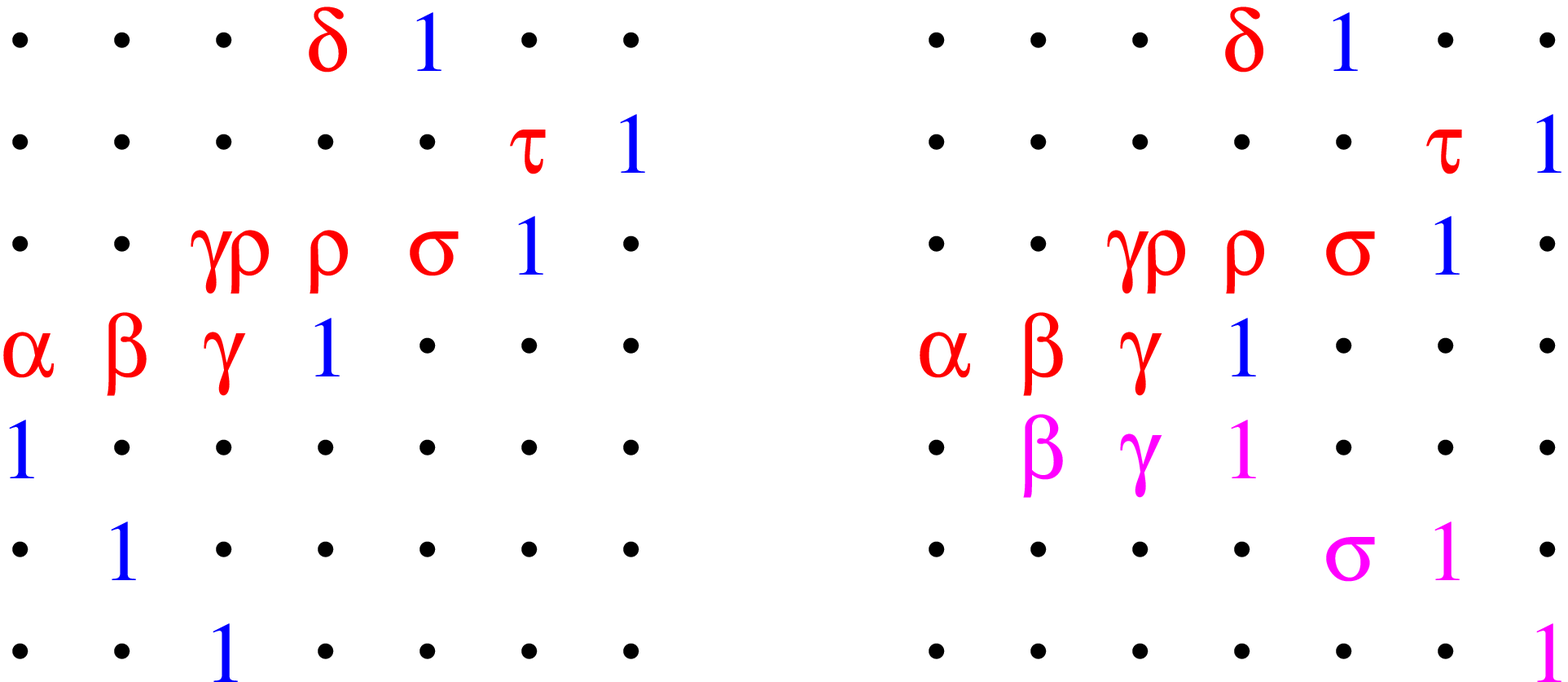}
$$
If $h_1,\ldots,h_7$ is the basis determined by the first matrix, then 
$h_1=g_2$, $h_2=g_4$, $h_3=g_3$, and $h_4=g_1$.
Thus $\Span{g_1,g_2,g_3,g_4}= \Span{h_1,h_2,h_3,h_4}$, which proves
$$
\pi_k\left(X_u\Epdot\bigcap X_{w_0 w}\Edot\right)\ 
=\ \pi_k\left(X_x\Epdot\bigcap X_{w_0y}\Edot\right).
$$

This is true even when $u,w,x,z$ do not satisfy the
extra hypotheses of Lemma~\ref{lem:skew_critical} (one may construct a
proof using the geometric analogs of the arguments that reduce
Theorem~\ref{thm:B}~({\em ii}) to Lemma~\ref{lem:skew_critical}).
We illustrate this on another example.

Here, let $n=7$, $k=3$, and 
$$
\begin{array}{r}u\ =\ 2134765\\w\ =\ 3571624\end{array}\qquad
\begin{array}{r}x\ =\ 2316475\\z\ =\ 3752164\end{array}
$$
Then the following four matrices  represent, respectively, the Schubert
cells  
$X^\circ_{w_0 w}\Edot$, $X^\circ_u\Epdot$,  $X^\circ_{w_0z}\Edot$, and 
$X^\circ_x\Epdot$:
$$
\epsfxsize=4.8in \epsfbox{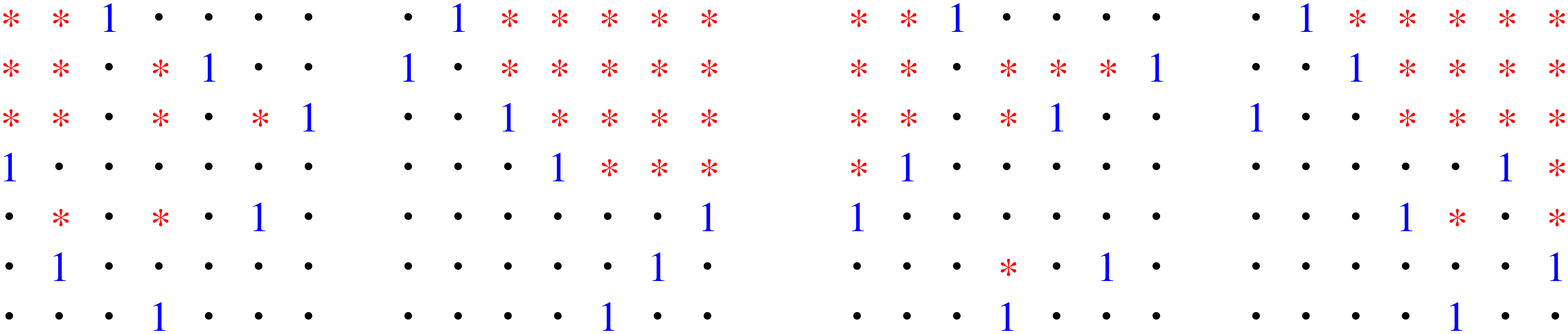}
$$
Consider the  (equivalent pairs of) parameterizations for flags in the
intersections of the cells,
$X^\circ_{w_0 w}\Edot\bigcap X^\circ_u\Epdot$ (the left-hand pair),
and  $X^\circ_{w_0z}\Edot\bigcap X^\circ_x\Epdot$ (the right-hand pair):
$$
\epsfxsize=4.8in \epsfbox{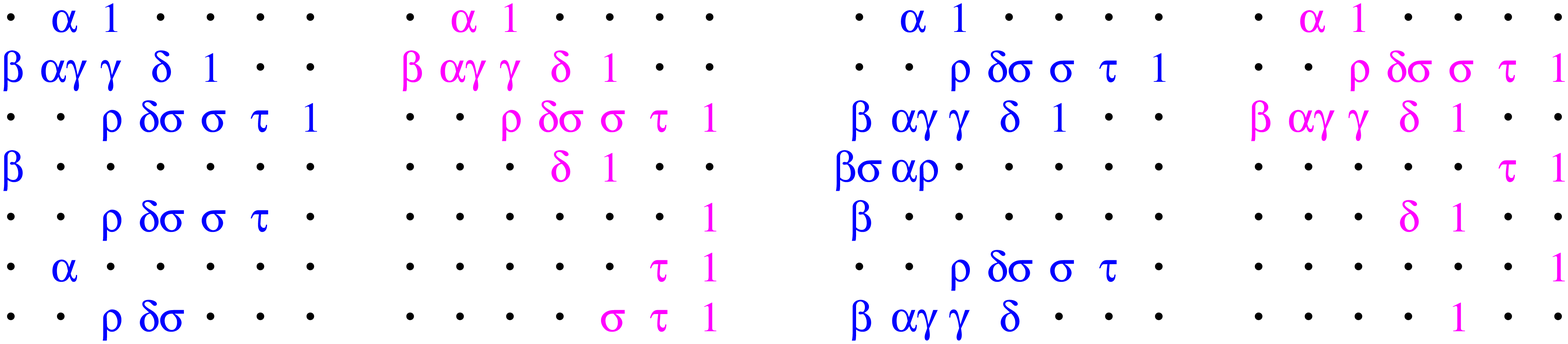}
$$
To see that each pair of matrices does indeed give the same flag, let 
$g_1,\ldots,g_7$, $g'_1,\ldots,g'_7$, $h_1,\ldots,h_7$, and 
$h'_1,\ldots,h'_7$ be the bases given by the four matrices (read
left-to-right).
Then, for $i=1,2,3$, $g_i=g'_i$ and $h_i=h'_i$.
Also,
\begin{eqnarray*}
g'_4&=& - \alpha g_1 + g_2 -g_4\\
g'_5 &=& g_3 - g_5\\
g'_6&=& g_3-\sigma g'_4 -\rho( g_1  - g_6)\\
g'_7&=& g_3-g_7
\end{eqnarray*}
and 
\begin{eqnarray*}
\qquad
h'_4&=& h_2 -\sigma h_3 +h_4 + (\gamma -\rho)h_1\\
h'_5 &=& h_3-\gamma h_1-h_5\\
h'_6&=& h_2-h_6\\
h'_7&=& h_3-h_7,
\end{eqnarray*}
thus, $\SPan{g_1,\ldots,g_7}=\SPan{g'_1,\ldots,g'_7}$
and  $\SPan{h_1,\ldots,h_7}=\SPan{h'_1,\ldots,h'_7}$.
As before, these parameterized bases
give dense subsets of flags in each of
$X_{w_0w}\Edot\bigcap X_u\Epdot$ and $X_{w_0z}\Edot\bigcap X_x\Epdot$.
Moreover, since $\Span{g_1,g_2,g_3}=\Span{h_1,h_2,h_3}$,
we see that
$$
\pi_k\left(X_{w_0w}\Edot\bigcap X_u\Epdot\right)\ 
=\ \pi_k\left(X_{w_0z}\Edot\bigcap X_x\Epdot\right).
$$

\subsection{Theorem~\ref{thm:C}~({\em ii})}
We complete Example~\ref{example:disjoint}, giving the geometric side of
the story.
The permutation $(1978)(26354)$ is the disjoint product of 
$\zeta=(1978)$ and $\eta=(26354)$.
Note that $u=372186945\leq_4   586913724=(\zeta\eta)u =:w$.
Let $\Gdot$ and $\Gpdot$ be the standard flags in ${\Bbb C}^9$.
The following matrices parameterize the Schubert cells 
$X^\circ_{w_0586913724}\Gdot$ and $X^\circ_{372186945}\Gpdot$:
$$
\epsfxsize=3.1in \epsfbox{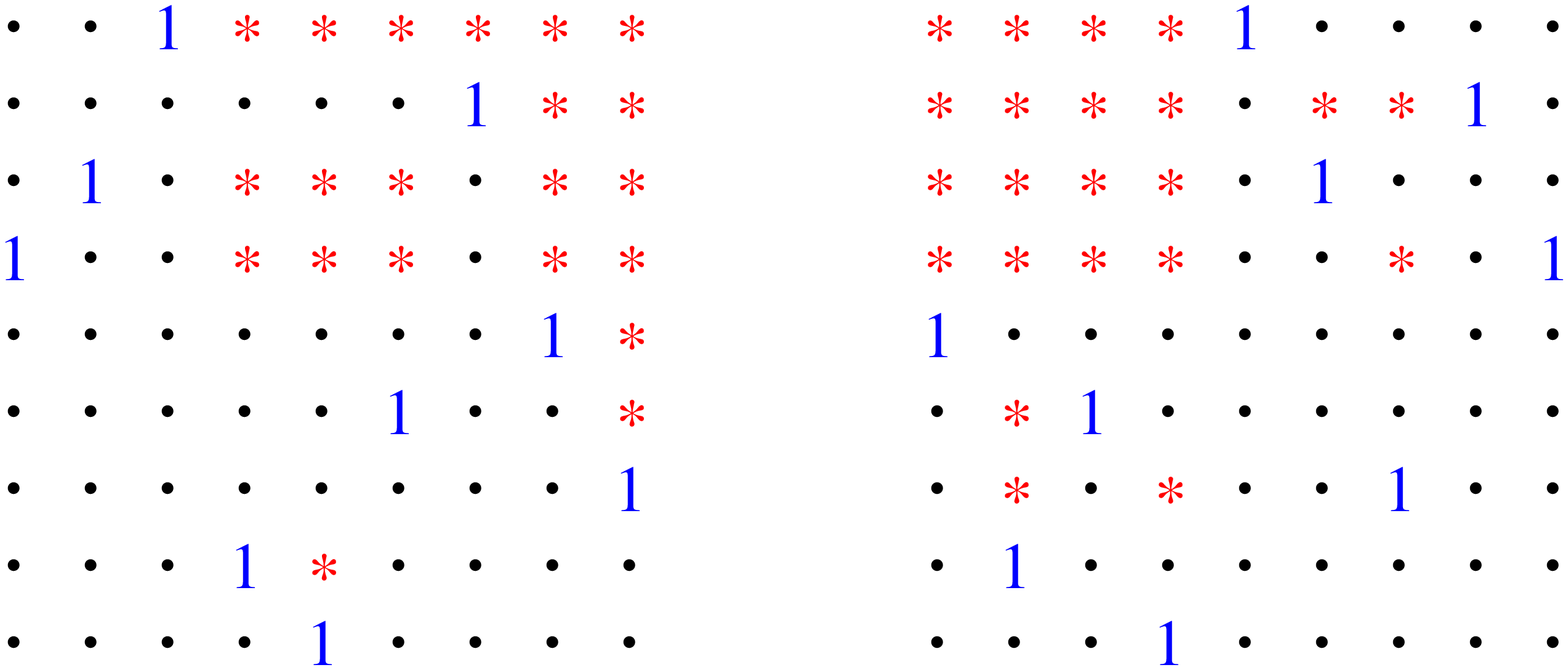}
$$
As before, here are two parameterized matrices, each of which give bases
defining the same flag in the intersection of the two cells:
$$
\epsfxsize=3.1in \epsfbox{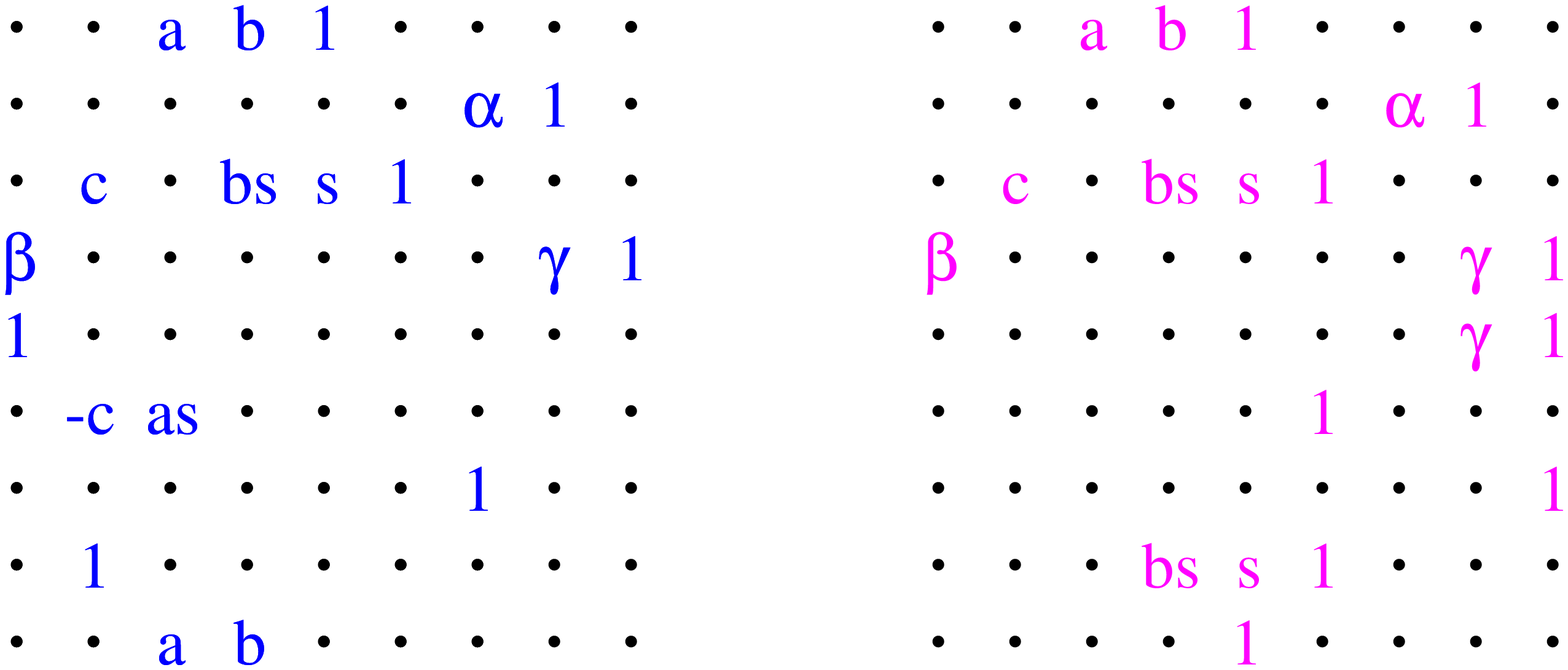}
$$
It is clearer to display two matrices `on top of each other':
$$
\epsfxsize=1.4in \epsfbox{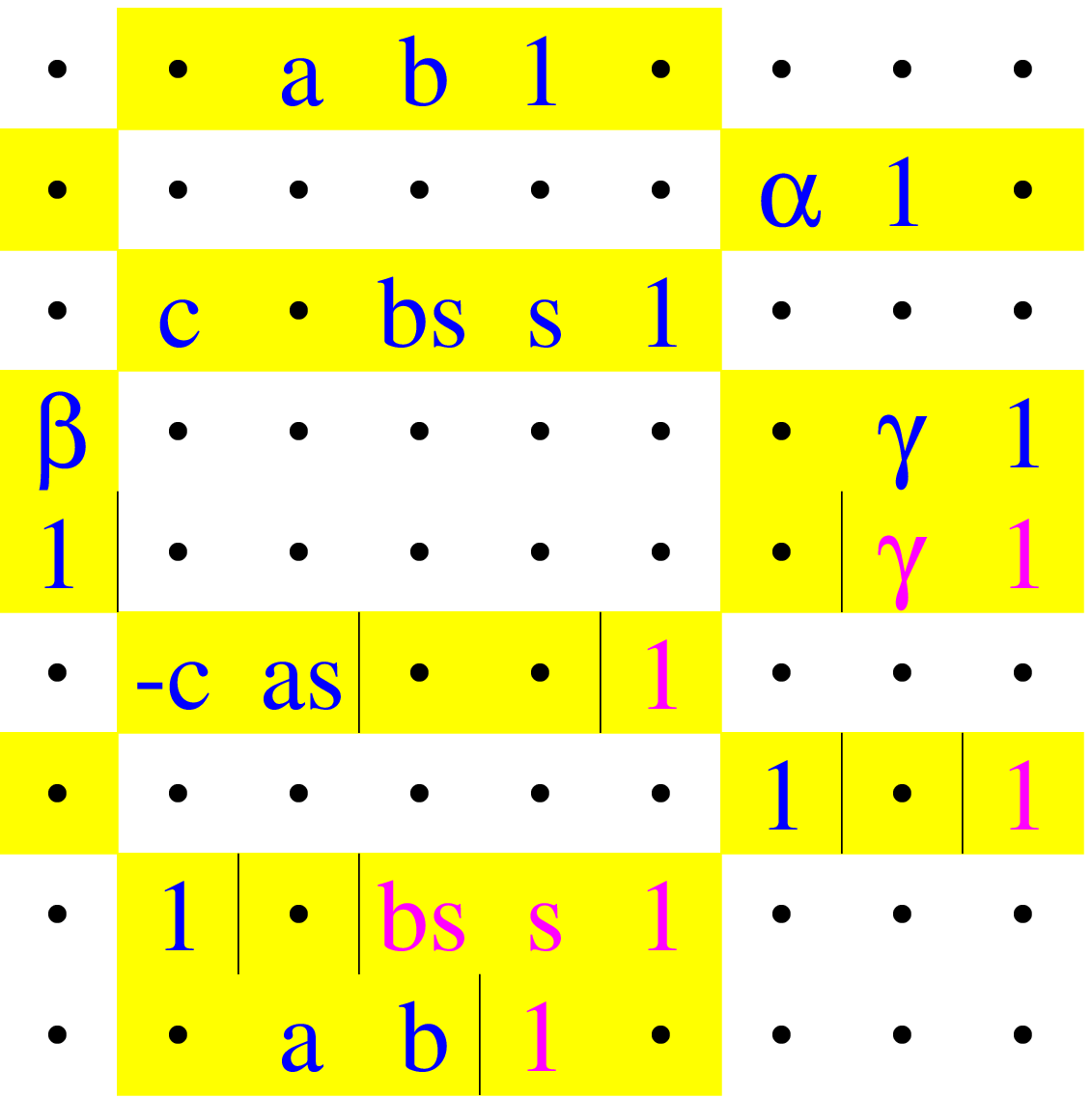}
$$
The vertical lines in the last 5 rows illustrate that
the left- and right-sides of those rows come from different (equivalent)
bases. 
The shading accentuates its `block form':
Let $Q=\{1,7,8,9\}$ and $P=\{2,4,5,7\}=u^{-1}(Q)$.
The $P^c=\{1,3,6,8,9\}=u^{-1}(Q^c)$, where $Q^c=\{2,3,4,5,6\}$.
Then the shaded regions are $(P\times Q) \bigcup (P^c\times Q^c)$.
We see that $\zeta':= (1423)$ and $\eta':=(15243)$ are uniquely defined
by $\phi_Q\zeta' = \zeta$ and $\phi_{Q^c}\eta'=\eta$.
Moreover, we may define permutations $v$ and $w$ as in
Lemma~\ref{lem:long}; let $v=2134$ and $w=21534$.
Then
\begin{enumerate}
\item $v\leq_2 \zeta'v = 3412$ and $w\leq_2 \eta' w = 45213$.
\item $u = \varepsilon_{P,Q}(v,w)$ and 
$(\zeta\eta)u= \varepsilon_{P,Q}(\zeta'v,\eta'w)$.
\end{enumerate}
For the last part of Lemma~\ref{lem:long}, 
let $\Fdot, \Fpdot$ be the standard flags in ${\Bbb C}^4$, 
and $\Edot, \Epdot$ the standard flags in ${\Bbb C}^5$.
Then the following four matrices parameterize the Schubert cells
 $X^\circ_{w_0\zeta'v}\Fdot$, $X^\circ_v\Fpdot$,
 $X^\circ_{w_0\eta'w}\Edot$, and $X^\circ_w\Epdot$: 
$$
\epsfxsize=4.in \epsfbox{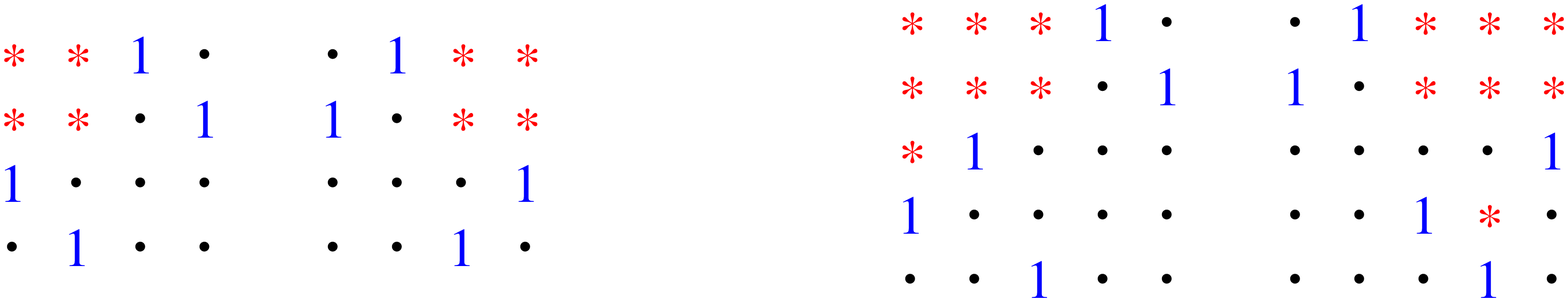}
$$
Then the following two matrices parameterize the two intersections.
Again, we have drawn two matrices on top of each other.
$$
\epsfxsize=1.9in \epsfbox{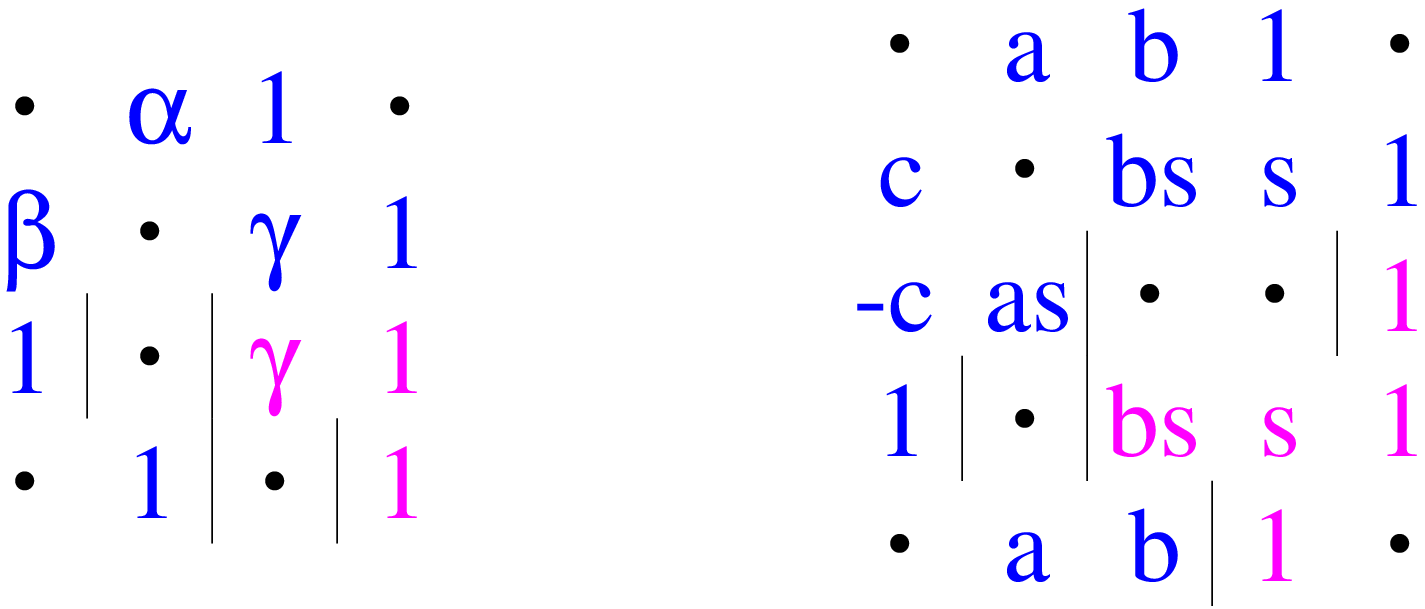}
$$
Next, note that $\Gdot=\psi_Q(\Edot,\Fdot)$ and 
$\Gpdot=\psi_{w_0^{(9)}Q}(\Epdot,\Fpdot)$.
Finally, verifying that 
$$
\psi_P\left[
\left(X_{w_0^{(4)}\zeta' v}\Edot \bigcap X_v\Epdot\right) \times
\left(X_{w_0^{(5)}\eta' w}\Fdot \bigcap X_w\Fpdot\right)\right]
$$
is equal to 
$$
X_{w_0^{(9)}(\zeta\eta)u}\Gdot \bigcap
X_u\Gpdot.
$$
may be done by comparing these parameterized matrices.

In the final part of the proof of Theorem~\ref{thm:C}~({\em ii}),
we compare images of 
these intersections under projections to Grassmannians.
The row spans of the next three parameterized matrices represent  
$\pi_2\left(X^\circ_{w_0^{(4)}\zeta' v}\Edot 
\bigcap X^\circ_v\Epdot\right)$, 
$\pi_2\left(X^\circ_{w_0^{(5)}\eta' w}\Fdot 
\bigcap X^\circ_w\Fpdot\right)$,  and 
$\pi_4\left(X^\circ_{w_0^{(9)}(\zeta\eta)u}\Gdot \bigcap
X^\circ_u\Gpdot\right)$ in each of $\Gr_2{\Bbb C}^4$,  $\Gr_2{\Bbb C}^5$,  
and $\Gr_4{\Bbb C}^9$,  respectively.
$$
\epsfxsize=3.5in \epsfbox{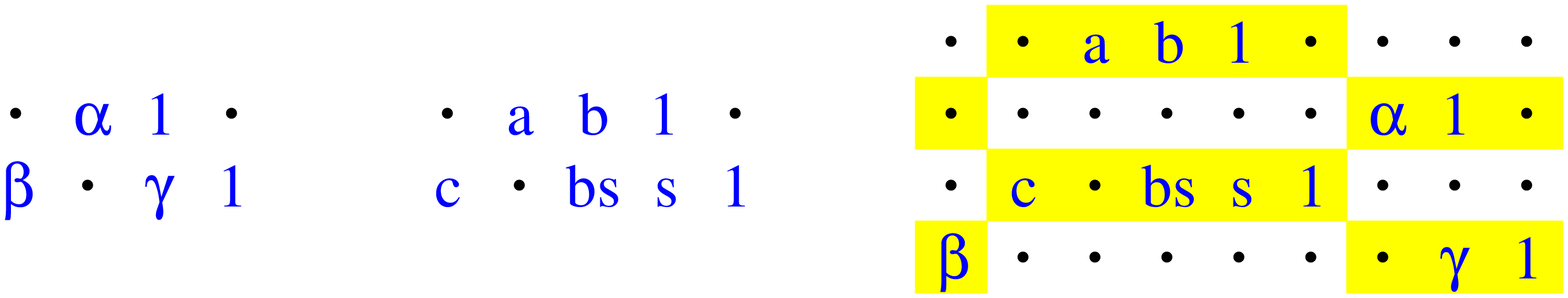}
$$
Thus 
$$
\pi_4\left(X_{w_0^{(9)}(\zeta\eta)u}\psi_Q(\Edot,\Fdot)
\bigcap X_u\psi_{w_0^{(9)}Q}(\Epdot,\Fpdot)\right)
$$
is equal to 
$$
\varphi_{2,2}\left(
\pi_2\left(X_{w_0^{(4)}\zeta'v}\Edot\bigcap X_v\Epdot\right) \times
\pi_2\left(X_{w_0^{(5)}\eta'w}\Fdot\bigcap X_w\Fpdot\right) \right).
$$

Consider the images of 
$X_{w_0^{(4)}\zeta'v}\Edot \times X_{w_0^{(5)}\eta'w}\Fdot$ 
and $X_v\Epdot\times X_w\Fpdot$ under $\psi_P$ in ${\Bbb F}\ell_9$:
$$
\epsfxsize=3.1in \epsfbox{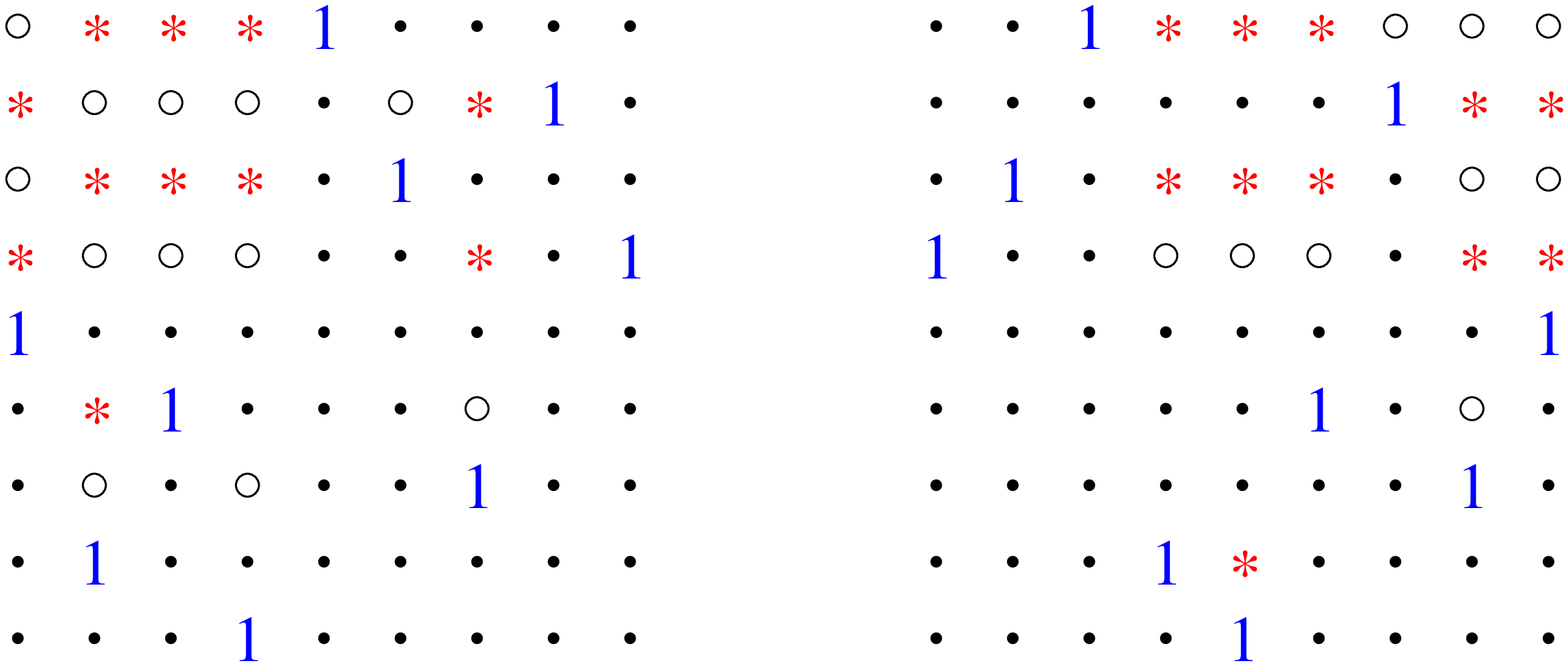}
$$
This should be compared with the first figure of this section, which
shows the cells  $X^\circ_{w_0586913724}\Gdot$ and 
$X^\circ_{372186945}\Gpdot$.
Here, the circles  (\begin{picture}(6,6)\thinlines
\put(3,3){\circle{5}}\end{picture}) indicate the `surprise' 
entries of 0; those which are not zero in that first figure.
This illustrates the two inclusions, and serves to illustrate
Lemma~\ref{lem:product_subset}: 
\begin{eqnarray*}
\psi_P \left(X_{w_0^{(4)}\zeta'v}\Edot \times 
X_{w_0^{(5)}\eta'w}\Fdot\right) &\subset &
 X^\circ_{w_0586913724}\Gdot\\
\psi_P \left(X_v\Epdot\times X_w\Fpdot\right)
&\subset &X^\circ_{372186945}\Gpdot
\end{eqnarray*}

\subsection{Theorem~\ref{thm:D}}
We illustrate  `cyclic shift'.
Let $u=21354$ and $w=45123$ so that $wu^{-1}=\zeta=(15243)$.
Define $x,z\in{\cal S}_5$ as in the proof of Theorem~\ref{thm:D}$'$
(\S~\ref{sec:thmd}) to be $x=31245$ and $z=53124$.
Then $zx^{-1}=(13542) = \zeta^{(12345)}$.
Here, $p=4$, $m=1$, and $l=2$.
Let $\Fdot, \Fpdot$ be the standard flags for ${\Bbb C}^5$.
We define $\Gdot=\SPan{e_5,e_1,e_2,e_3,e_4}$ and 
$\Gpdot = \SPan{e_4,e_3,e_2,e_1,e_5}$.
Then, with respect to these flags,  the Schubert cells
$X^\circ_{w_0 w}\Fdot$, $X^\circ_u\Fpdot$, $X^\circ_{w_0 z}\Gdot$,  and
$X^\circ_x\Gpdot$ are:
$$
\epsfxsize=3.7in \epsfbox{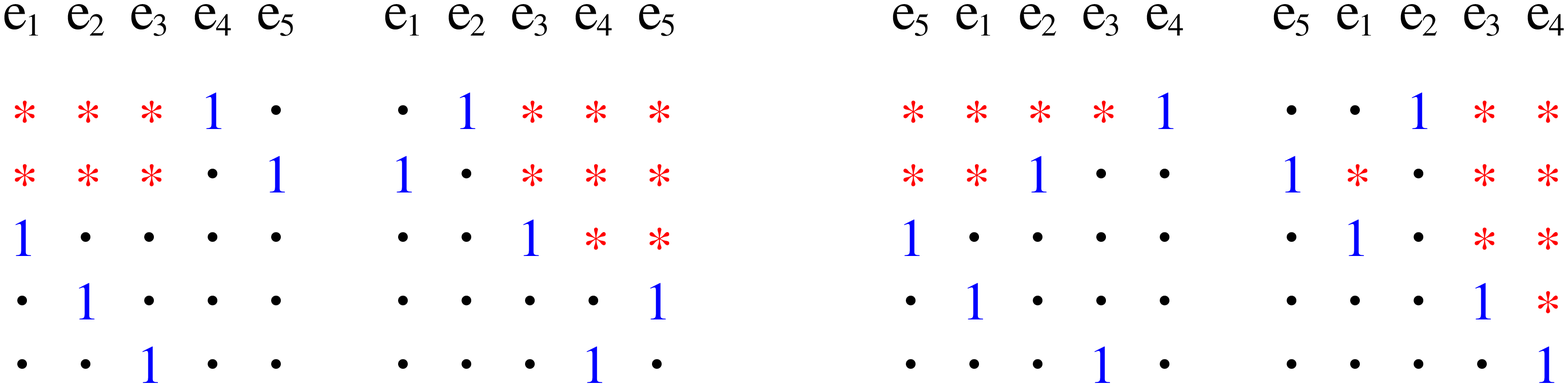}
$$
Here, since the flags are different, the columns of the matrices on the
right correspond to different elements of our fixed basis,
$e_1,e_2,e_3,e_4,e_5$, as indicated.

Here are two matrices giving (equivalent) parameterized bases for
flags in the intersection of the cells 
$X^\circ_{w_0 w}\Fdot\bigcap X^\circ_u\Fpdot$:
$$
\epsfxsize=1.8in \epsfbox{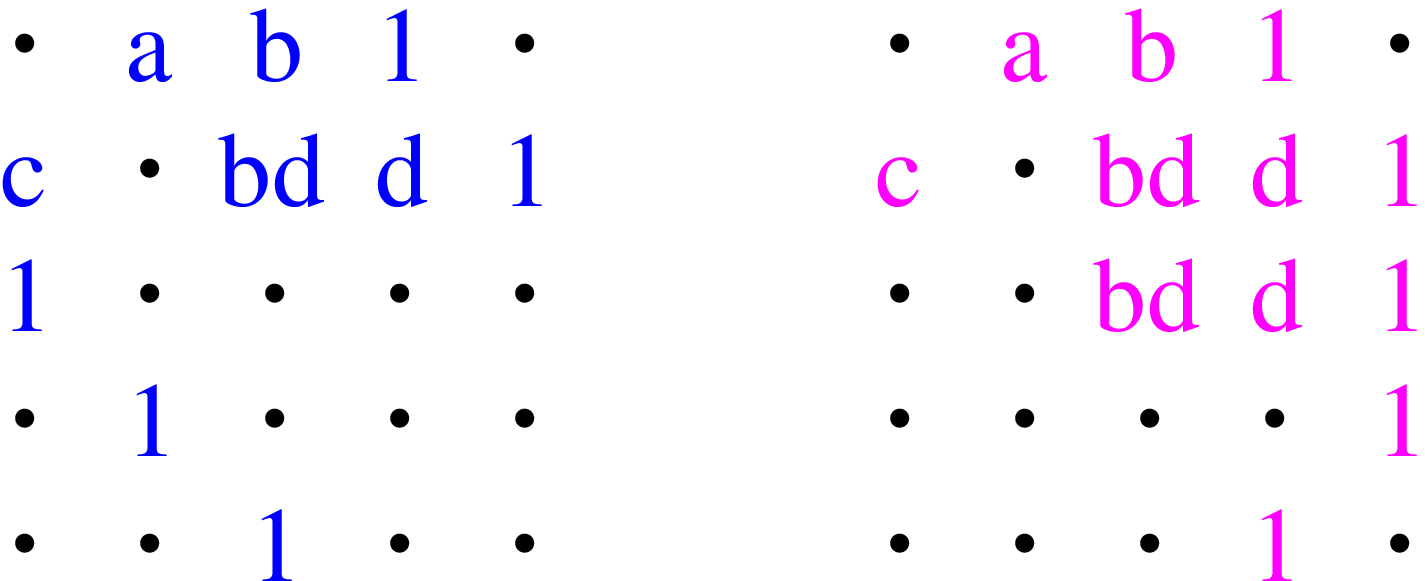}
$$
Here $a,b,c,d\in {\Bbb C}^\times$, showing that $({\Bbb C}^\times)^4$
parameterizes the set $A$ of the intersection of cells.
Let $g_1,\ldots,g_5$ be the basis given by the left matrix and 
$g'_1,\ldots,g'_5$ the basis given by the right matrix.
Since
$$
g_2(a,b,c,d)\ =\ e_5\; + c\,e_1+bd\,e_3+ d\,e_4,
$$ 
$(\beta_1,\beta_2,\beta_3,\beta_4)=(c, 0, bd, d)$
are regular functions on $A$.
Also, since 
$$
e_5 = - d\,g_1  + g_2 -c\,g_3 + da\,g_4,
$$
$\delta_1= -d$, $\delta_3 = -c$, and $\delta_4=da$ are regular
functions on $A$ with $\delta_4$ nowhere vanishing.
The bases $h_1,\ldots,h_5$ and
$h'_1,\ldots,h'_5$ defined in the proof of Theorem~\ref{thm:D}$'$, are
parameterized by the following two matrices:
$$
\epsfxsize=1.8in \epsfbox{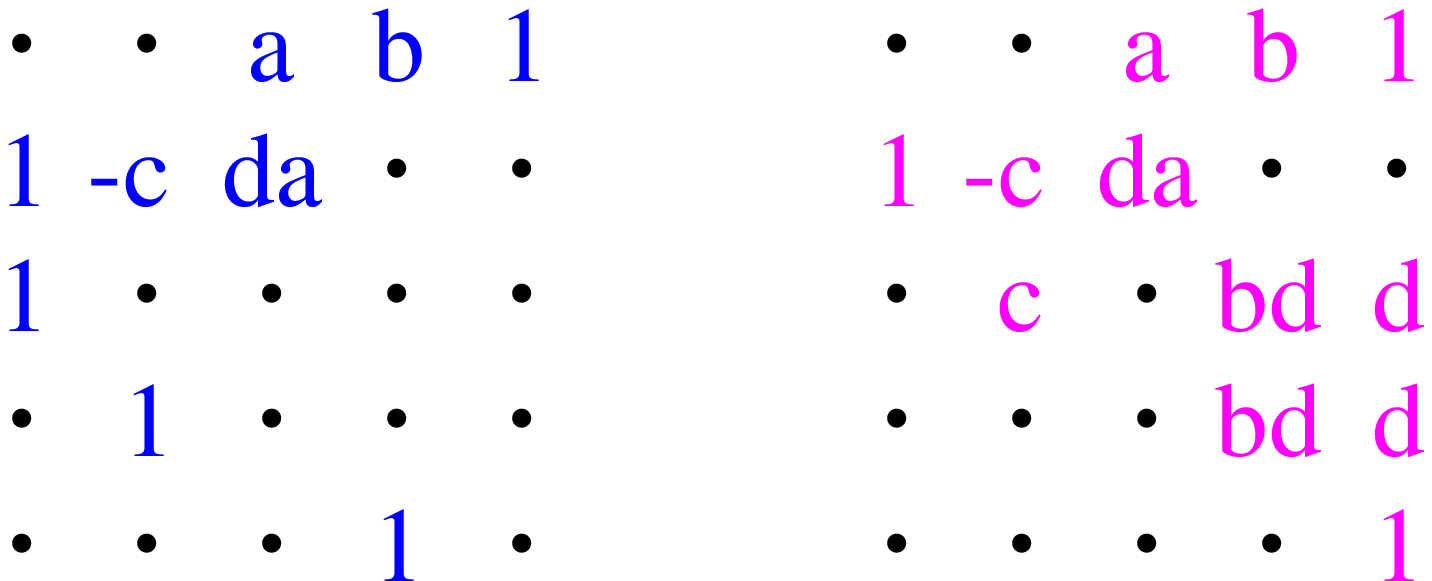}
$$
These matrices give equivalent bases, and hence
define the same flag.
Also, comparing them to the rightmost two matrices in the first figure
of this subsection, shows this flag is in the intersection
$X^\circ_{w_0 z}\Gdot\bigcap X^\circ_x\Gpdot$.
Since the first two rows of each matrix have the same span, 
$$
\pi_2\left(X^\circ_{w_0 w}\Fdot\bigcap X^\circ_u\Fpdot\right)\ =\ 
\pi_2\left(X^\circ_{w_0 z}\Gdot\bigcap X^\circ_x\Gpdot\right),
$$
which is the main geometric result needed to deduce Theorem~\ref{thm:D}$'$.

\section{Combinatorial and algebraic examples}

\subsection{Suborders of ${\cal S}_4$}
The Bruhat order is one of our main objects of study in this paper.
Here is a picture of the (full) Bruhat order and  the 
2-Bruhat order on ${\cal S}_4$.
$$
\epsfxsize=2.5in \epsfbox{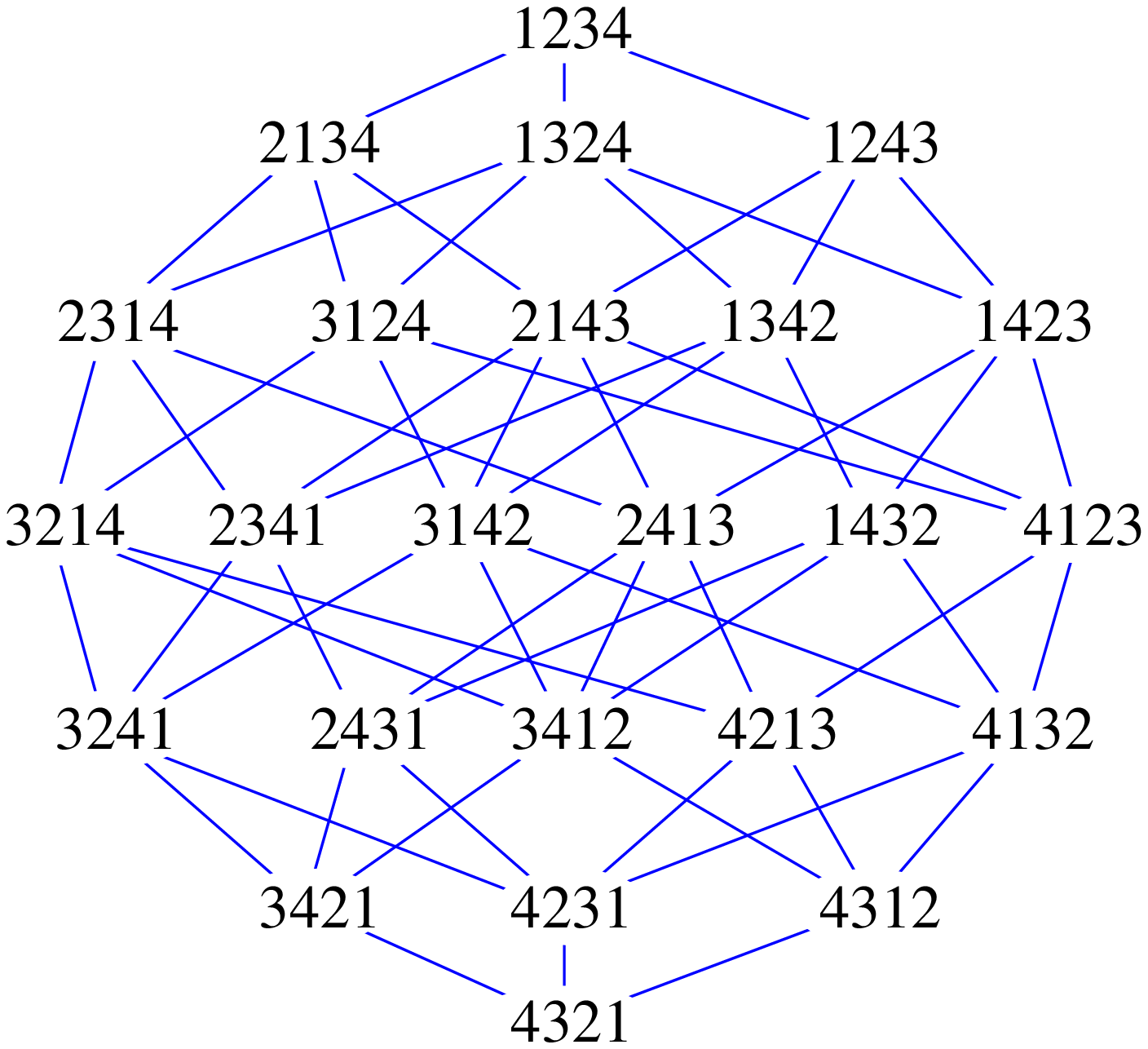}\qquad
\epsfxsize=2.5in \epsfbox{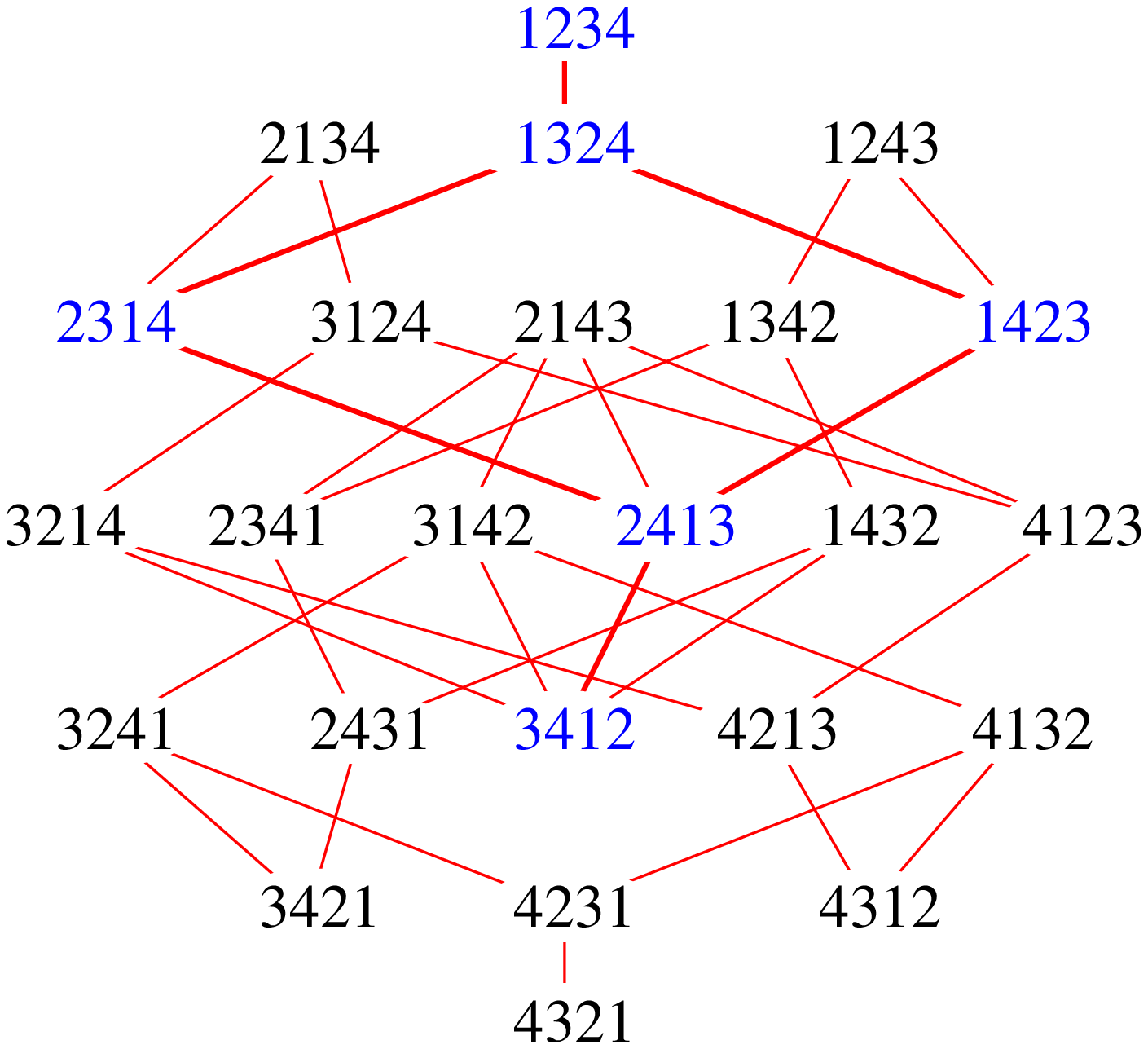}
$$
For comparison, here is the $\preceq$-order on ${\cal S}_4$ (reproduced
from \S\ref{sec:new:order}).
$$\epsfxsize=3.3in \epsfbox{figures/new_order.S_4.eps}$$

\subsection{Chains in the $P$-Bruhat order}
Theorem~\ref{thm:chains} describes the relation between chains in the
$P$-Bruhat order and the structure constants $c^w_{u\,v}$, when $v$ is
a minimal coset representative in $vP$.
We consider an instance of this.
Let  $P:=\Span{(1,2),(4,5)}\subset{\cal S}_5$.
Then $32154\leq_P 45312$
and this is the interval $[32154,\;45312]_P$:
$$
\epsfxsize=2.4in \epsfbox{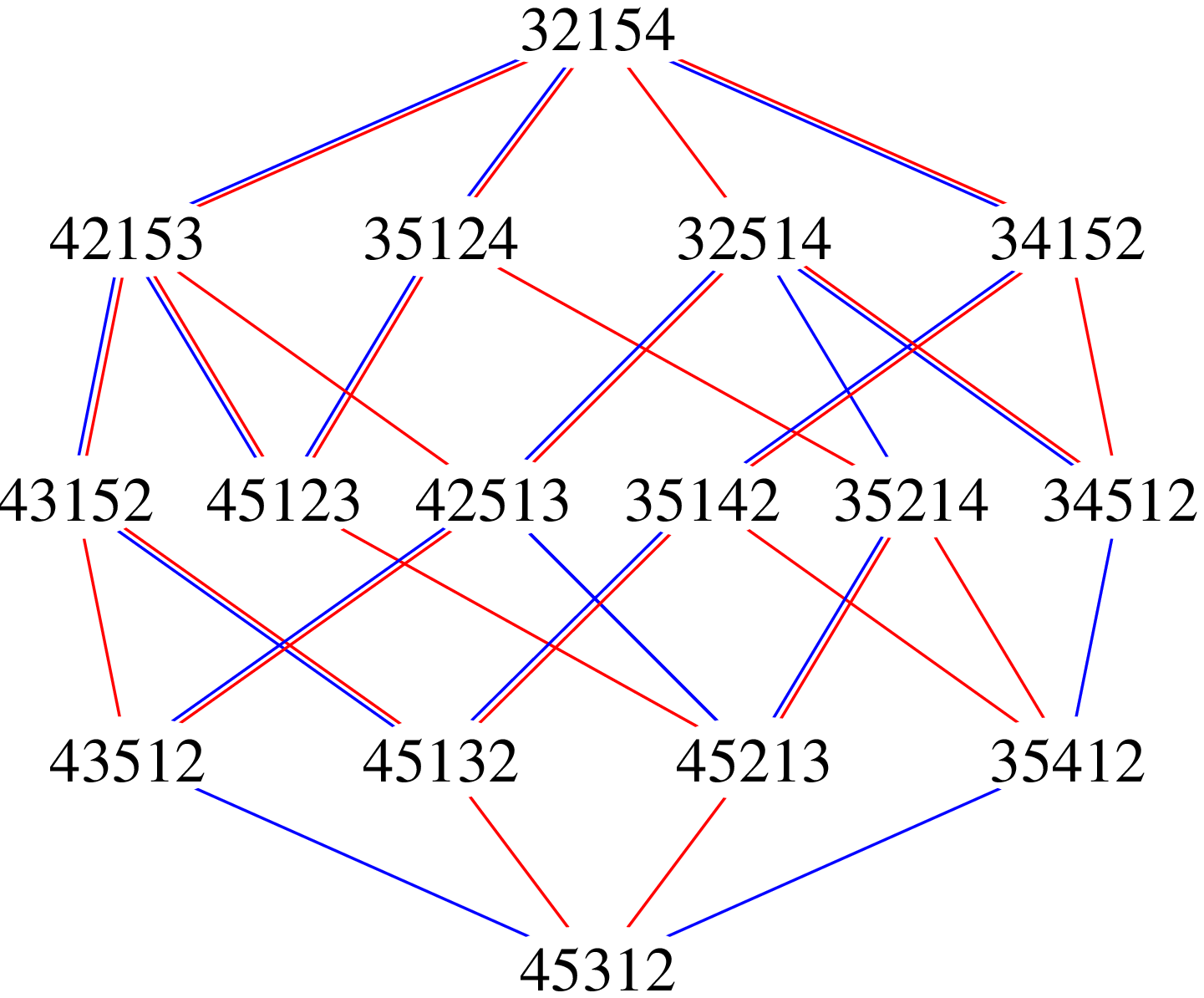}
$$
The multiple edges are those with two possible
colourings.
One may verify that $f^{45312}_{32154}(P) = 57$.
To check Theorem~\ref{thm:chains}, we first compute
$c^{45312}_{32154\: v}$ for those $v\in {\cal S}_5$ of length 4 which
are minimal in their $P$-coset.
$$
25134, \quad 34125, \quad 24315, \quad 15324,  \quad 14523,
 \quad \mbox{and} \quad  23514.
$$
The first two are Grassmannian of descent 2, and the last two are
Grassmannian of descent 3.
Since $32154\not\leq_2 45312$,
we have 
$$
c^{45312}_{32154\ 25134} \quad =\quad c^{45312}_{32154\ 34125}
\quad =\quad 0,
$$
Let $\zeta=(13425)$.
Then $45312 = \zeta\cdot 32145$ and $(13425)^{(12345)} = (12435)$.
Since $(12435)=v(\,%
\begin{picture}(6,9)    \put(0,0){\line(1,0){6}}
\put(0,0){\line(0,1){9}}\put(0,3){\line(1,0){6}}
\put(3,0){\line(0,1){9}}\put(0,6){\line(1,0){6}}
\put(6,0){\line(0,1){6}}\put(0,9){\line(1,0){3}}
\end{picture}\,,\, 3)\cdot
v(%
\begin{picture}(3,3)    
\put(0,0){\line(0,1){3}}\put(0,0){\line(1,0){3}}
\put(3,0){\line(0,1){3}}\put(0,3){\line(1,0){3}}
\end{picture}\,,\, 3)^{-1}$
and $32154\leq_3 45312$, Theorem~\ref{thm:D} implies 
$$
c^{45312}_{32154\ 14523} \quad = \quad
c^{\,\begin{picture}(4,6) \put(2,0){\line(1,0){2}}
\put(0,2){\line(0,1){4}}\put(0,2){\line(1,0){4}}
\put(2,0){\line(0,1){6}}\put(0,4){\line(1,0){4}}
\put(4,0){\line(0,1){4}}\put(0,6){\line(1,0){2}}
\end{picture}}_{\,\begin{picture}(4,4) 
\put(0,0){\line(0,1){4}}\put(0,0){\line(1,0){4}}
\put(2,0){\line(0,1){4}}\put(0,2){\line(1,0){4}}
\put(4,0){\line(0,1){4}}\put(0,4){\line(1,0){4}}
\end{picture}} \quad = \quad 1
\qquad \mbox{and}\qquad
c^{45312}_{32154\ 23514} \quad = \quad
c^{\,\begin{picture}(4,6) \put(2,0){\line(1,0){2}}
\put(0,2){\line(0,1){4}}\put(0,2){\line(1,0){4}}
\put(2,0){\line(0,1){6}}\put(0,4){\line(1,0){4}}
\put(4,0){\line(0,1){4}}\put(0,6){\line(1,0){2}}
\end{picture}}_{\,\begin{picture}(4,6) 
\put(0,0){\line(0,1){6}}\put(0,0){\line(1,0){4}}
\put(2,0){\line(0,1){6}}\put(0,2){\line(1,0){4}}
\put(4,0){\line(0,1){2}}\put(0,4){\line(1,0){2}}
			\put(0,6){\line(1,0){2}}
\end{picture}} \quad = \quad 1.
$$
Next, let $\Fdot,\Fpdot, \Fppdot$ be in general position.
If  $\Edot \in X_{15324}\Fdot\bigcap X_{32154}\Fpdot$, then 
$E_2\subset F'_4$ and $E_2 \supset F_1$, contradicting 
$\Fdot$ and $\Fpdot$ in general position.
Thus 
$$
c^{45312}_{32154\ 15324}\quad=\quad
\#\left(X_{15324}\Fdot\bigcap X_{32154}\Fpdot
\bigcap X_{w_045312}\Fppdot\right) \quad=\quad 0.
$$
To compute $c^{45312}_{32154\ 24315}= \deg\left(
{\frak S}_{w_045312}\cdot {\frak S}_{32154}\cdot 
{\frak S}_{24315}\right)$, note that 
${\frak S}_{w_0\;45312}={\frak S}_{21354}=
{\frak S}_{(1,2)}\cdot{\frak S}_{(4,5)}$.
Two applications of Monk's formula show 
$c^{45312}_{32154\ 24315}= 1$.
(The other computations could also have proceeded via Monk's formula.)

To compute $f^v_e(P)$ for these minimal coset representatives, consider
the part of the $P$-Bruhat order rooted at $e$ and restricted to
permutations of length at most 4:
$$
\epsfxsize=2.4in \epsfbox{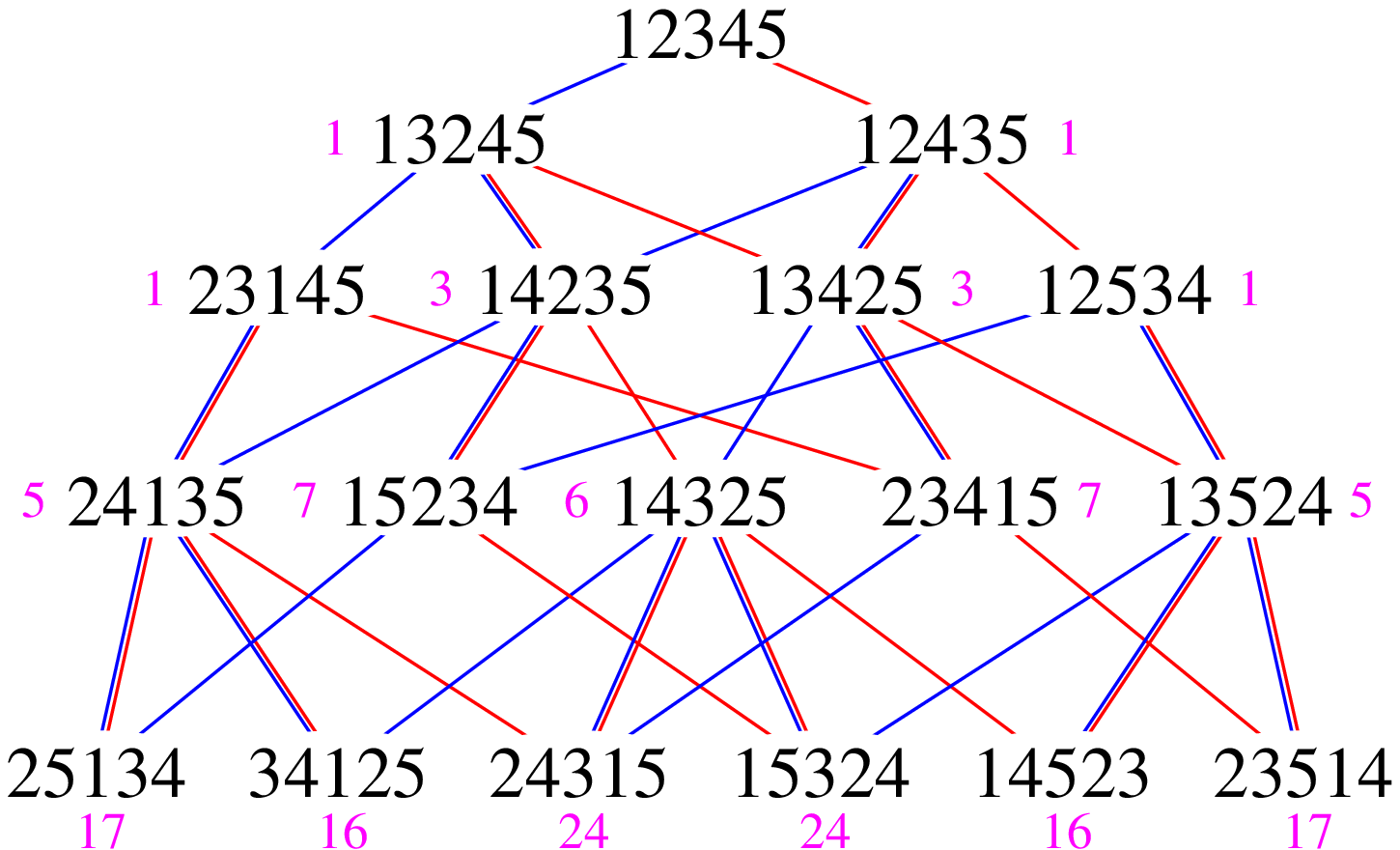}
$$
The small numbers adjacent to each permutation $v$ are $f^v_e(P)$.
Thus
$$
\sum_v f^v_e(P)\; c^{45312}_{32154\: v}\quad =\quad
17\cdot 0\ +\ 16\cdot 0\ +\ 24\cdot 1\ +\ 24\cdot 0
\ +\ 16\cdot 1\ +\ 17\cdot 1 \quad =\quad 57,
$$
which equals $f^{45312}_{32154}(P)$.

\subsection{Instance of Theorem~\ref{thm:substitution}}
We consider $\Psi_{\{1,3,5,\ldots\}}({\frak S}_{516432})$.
\begin{eqnarray*}
{\frak S}_{516432}  &=& \ 
x_1^4x_2^2x_3^3x_5 + x_1^4x_2x_3^3x_4x_5 + x_1^4x_3^3x_4^2x_5\\
&& + x_1^4x_2^3x_3^2x_5 + x_1^4x_2^2x_3^3x_4 + 
x_1^4x_2^2x_3^2x_4x_5 + 
x_1^4x_2x_3^3x_4^2 + x_1^4x_2x_3^2x_4^2x_5\\
&& + 
x_1^4x_2^3x_3^2x_4 + x_1^4x_2^3x_3x_4x_5 + x_1^4x_2^2x_3^2x_4^2 + 
x_1^4x_2^2x_3x_4^2x_5 \\
&& + 
x_1^4x_2^3x_3x_4^2 + x_1^4x_2^3x_4^2x_5.
\end{eqnarray*}
$\Psi_{\{1,3,5,\ldots\}}({\frak S}_{516432})
= {\frak S}_{516432}(y_1,z_1,y_2,z_2,y_3,z_3,\ldots)$,
which is
$$
y_1^4y_2^3y_3(z_1^2 + z_1z_2 + z_2^2)\  + \ 
y_1^4y_2^2y_3(z_1^3 + z_1z_2^2 + z_1^2z_2)\  +\ 
y_1^4y_2^3(z_1^2z_2 + z_1z_2^2)$$
$$+\ 
(y_1^4y_2^2 + y_1^4y_2y_3)(z_1^3z_2 + z_1^2z_2^2)\ +\ 
(y_1^4y_2+y_1^4y_3)z_1^3z_2^2.
$$

\begin{figure}[htb]
    $$\epsfxsize=4.in \epsfbox{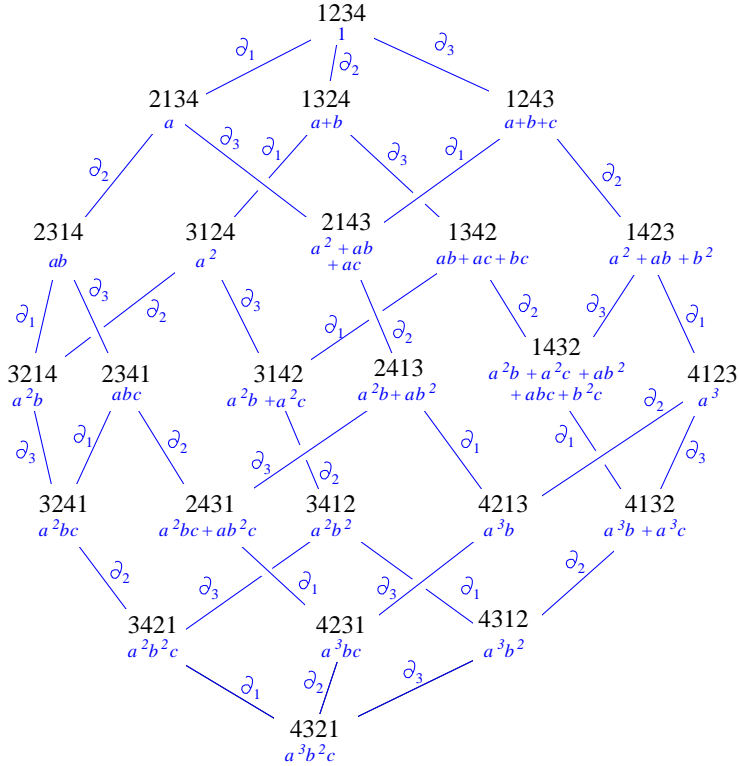}$$
\caption{Schubert polynomials in ${\cal S}_4$\label{fig:Schubert_S_4}}
  \end{figure}

Using the definition of Schubert polynomials in \S\ref{sec:schubert},
one may check
$$
\begin{array}{lcl}
{\frak S}_{54213} \ =\  x_1^4x_2^3x_3  &\qquad&
{\frak S}_{53214}  \ =\  x_1^4x_2^2x_3\\
{\frak S}_{54213}  \ =\  x_1^4x_2^3  &\qquad&
{\frak S}_{53124}  \ =\ x_1^4x_2^2\\
{\frak S}_{52314} \ =\ x_1^4x_2x_3&\qquad&
{\frak S}_{51324}  \ =\ x_1^4x_2 + x_1^4x_3
\end{array}
$$
The Schubert polynomials ${\frak S}_w$ for 
$w\in {\cal S}_4$ are indicated in Figure~\ref{fig:Schubert_S_4}.
The Schubert polynomial ${\frak S}_w$ is written below the
permutation $w$, and these data are 
displayed at the vertices of the permutahedron 
(Cayley graph of ${\cal S}_4$).
The divided difference operators are displayed on the edges of
this figure.

We see that 
$\Psi_{\{1,3,5,\ldots\}}{\frak S}_{516432} = 
{\frak S}_{516432}(y_1,z_1,y_2,z_2,y_3,z_3,\ldots)$
is equal to 
$$
{\frak S}_{54213}(y)
{\frak S}_{1423}(z)\  +\  
{\frak S}_{53214}(y)\left[
{\frak S}_{4123}(z) + 
{\frak S}_{2413}(z)\right]\   + \  
{\frak S}_{54123}(y){\frak S}_{2413}(z)$$
$$ +\  
\left[{\frak S}_{53124}(y)+
{\frak S}_{52314}(y)\right]\left[
{\frak S}_{4213}(z)+
{\frak S}_{3412}(z)\right]\   +\  
{\frak S}_{51324}(y)
{\frak S}_{4312}(z).
$$

\subsection{Automorphisms of $({\cal S}_\infty,\preceq)$}
The definition of the $k$-Bruhat orders imply that if $u,w\in{\cal
S}_n$, and $k<n$, then the following are equivalent:
$$
u\leq_k w \qquad
w_0 w\leq_k w_0 u\qquad
w w_0 \leq_{n-k} u w_0 \qquad
w_0w w_0 \leq_{n-k} w_0u w_0.
$$
These induce the following isomorphisms (which were stated in
Theorem~\ref{thm:new_order}) of intervals in the $\preceq$-order on 
$\cal S_\infty$.
Suppose $\zeta\in {\cal S}_n$ and $\overline{\zeta}=w_0\zeta w_0$.
Then
$$
[e,\zeta]_\preceq\quad \simeq\quad 
[e,\overline{\zeta}]_\preceq\quad \simeq\quad 
[e,\zeta^{-1}]^{\mbox{\scriptsize op}}_\preceq\quad \simeq\quad 
[e,\overline{\zeta}^{-1}]^{\mbox{\scriptsize op}}_\preceq.
$$
These are illustrated in the posets displayed in
%the next section.
\S\ref{sec:schensted_example}. 

\subsection{Canonical algorithms?}
Besides Algorithm~\ref{alg:chain}, there are three other `canonical'
algorithms for finding a chain between $u$ and $w$ when
$u\leq_k w$, each induced from Algorithm~\ref{alg:chain} by one of the
automorphisms of the previous section.  
For example, here is one.

\begin{alg}[Produces a chain in the $k$-Bruhat order]
\mbox{ }

\noindent{\tt input: }Permutations $u,w\in {\cal S}_\infty$ with 
$u\leq_k w$.

\noindent{\tt output: }A chain in the $k$-Bruhat order from $w$ to $u$.

Output $w$.
While $u\neq w$, do
\begin{enumerate}
\item[1] Choose $a\leq k$ with $w(a)$ maximal subject to $u(a)< w(a)$.
\item[2] Choose $k< b$ with $w(b)$ minimal subject to 
$w(b)\leq u(a)< u(b)$.
\item[3] $u:=u(a,b)$, output $u$.   
\end{enumerate}
\end{alg}
In general, these algorithms produce different chains.
In $S_7$, consider the two permutations
$2317546 <_3 4671235$.
Here are chains produced by the four algorithms:
$$
\begin{array}{lclclcl}
2317546 &\ & 2317546 &\ & 2317546 &\ & 2317546 \\
2417536 &\ & 2417536 &\ & 2371546 &\ & 2371546 \\
2517436 &\ & 2517436 &\ & 2571346 &\ & 2571346 \\
2617435 &\ & 4517236 &\ & 2671345 &\ & 3571246 \\
4617235 &\ & 4617235 &\ & 3671245 &\ & 4571236 \\
4671235 &\ & 4671235 &\ & 4671235 &\ & 4671235 
\end{array}
$$

\subsection{Schensted insertion and the
$c^w_{u\,v(\lambda,k)}$}\label{sec:schensted_example} 
In \S\ref{sec:further}, we discussed how the conclusion of
Theorem~\ref{thm:skew_shape_prime} holds for many permutations in 
${\cal S}_6$, even most which are not skew permutations.
We illustrate that here.

\begin{figure}[htb]
$$\epsfxsize=5.4in \epsfbox{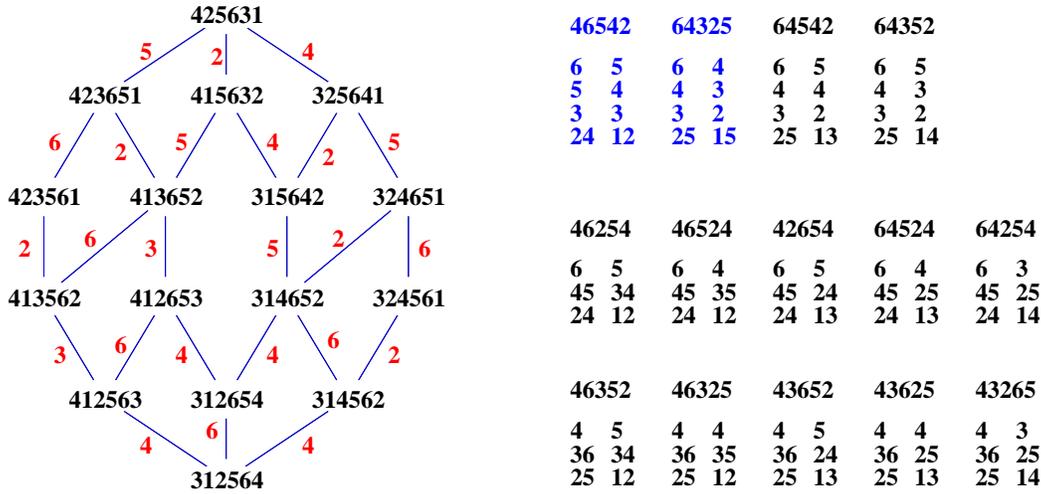}$$
\caption{Labeled Hasse diagram of $[214365,\;345612]_4$ and 
Schensted insertion\label{fig:125634}}
\end{figure}

Let $\zeta = (143562)$.
Then $214365 \leq_4 \zeta \cdot 214365 = 345612$.
In Figure~\ref{fig:125634}, we display the labeled Hasse diagram of 
$[214365,\;345612]_4$
and beside it a table of the words of the 14 chains in 
this interval, each displayed above its insertion and recording tableau.

\begin{figure}[htb]
$$\epsfxsize=5.4in \epsfbox{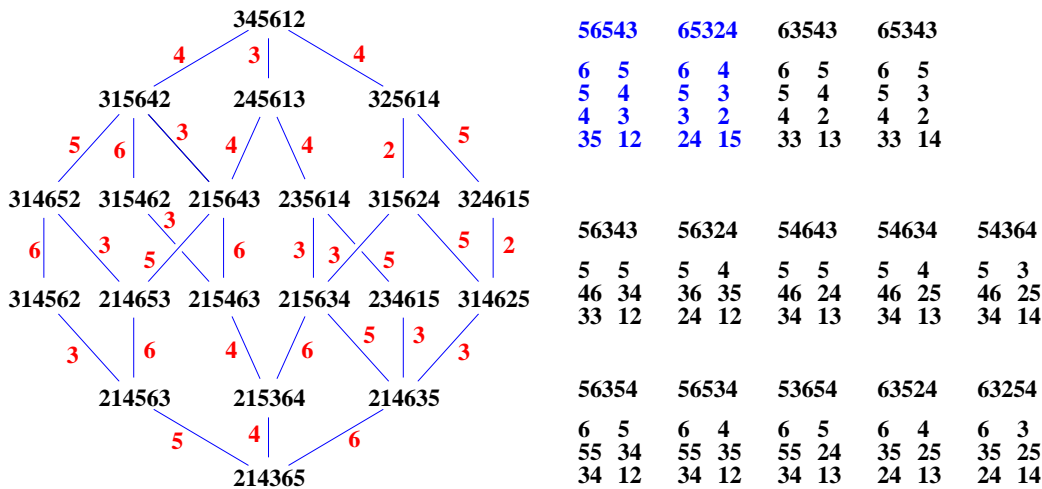}$$
\caption{Labeled Hasse diagram of $[312564,\,425635]_4$ and 
Schensted insertion\label{fig:145236}}
\end{figure}

Note that $ \eta:=(125634)=\zeta^{(123456)}$
and $312564 \leq_4 \eta\cdot 312564 = 425635$.
We continue this example, and illustrate Theorem~\ref{thm:D}.
In Figure~\ref{fig:145236} are the labeled Hasse diagram of 
$[312564,\,425635]_4$, and the insertion and recording tableaux for all 
14 chains in this interval.

\begin{figure}[htb]
$$\epsfxsize=4.6in \epsfbox{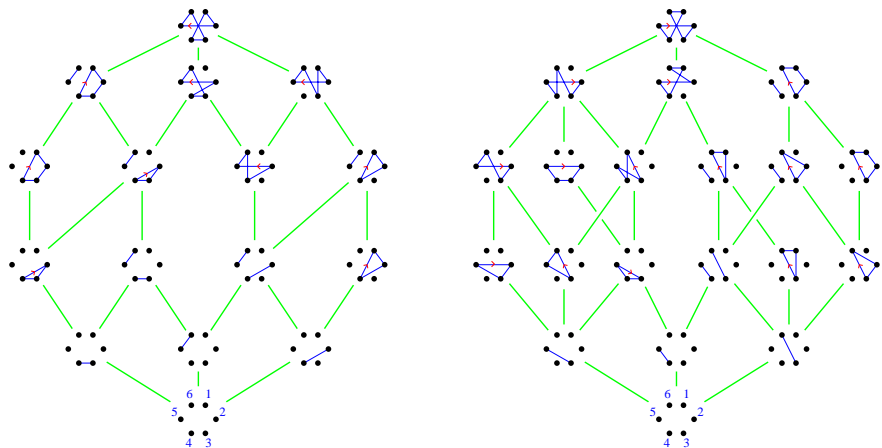}$$
\caption{Geometric graphs of permutations in $[e,\zeta]_\preceq$ and 
$[e,\eta]_\preceq$\label{fig:cyclic_geometric}}
\end{figure}

For these last two intervals, it is interesting to view them with 
the permutation $v\in [u,w]_k$ replaced by the geometric graph of
$vu^{-1}$, as illustrated in Figure~\ref{fig:cyclic_geometric}.
This gives an idea of the effect of a `cyclic shift' on the 
$\preceq$-order.

\subsection{Simplicial complexes and  $\leq_k$}
%The $[u,w]_k$ are {\em not} shellable}
In the theory of partially ordered sets, one often constructs a
simplicial complex $\Delta(P)$ from a poset, $P$.
We compute one such for an interval in the $k$-Bruhat order,
which shows these intervals are not in general shellable.
We illustrate this with one example drawn from this paper.
In Example~\ref{ex:neworder}, we considered the interval
$[21342,\,45123]_2$.
We display that interval below, together with the Hasse diagram of an
isomorphic poset:
$$
\epsfxsize=3.1in \epsfbox{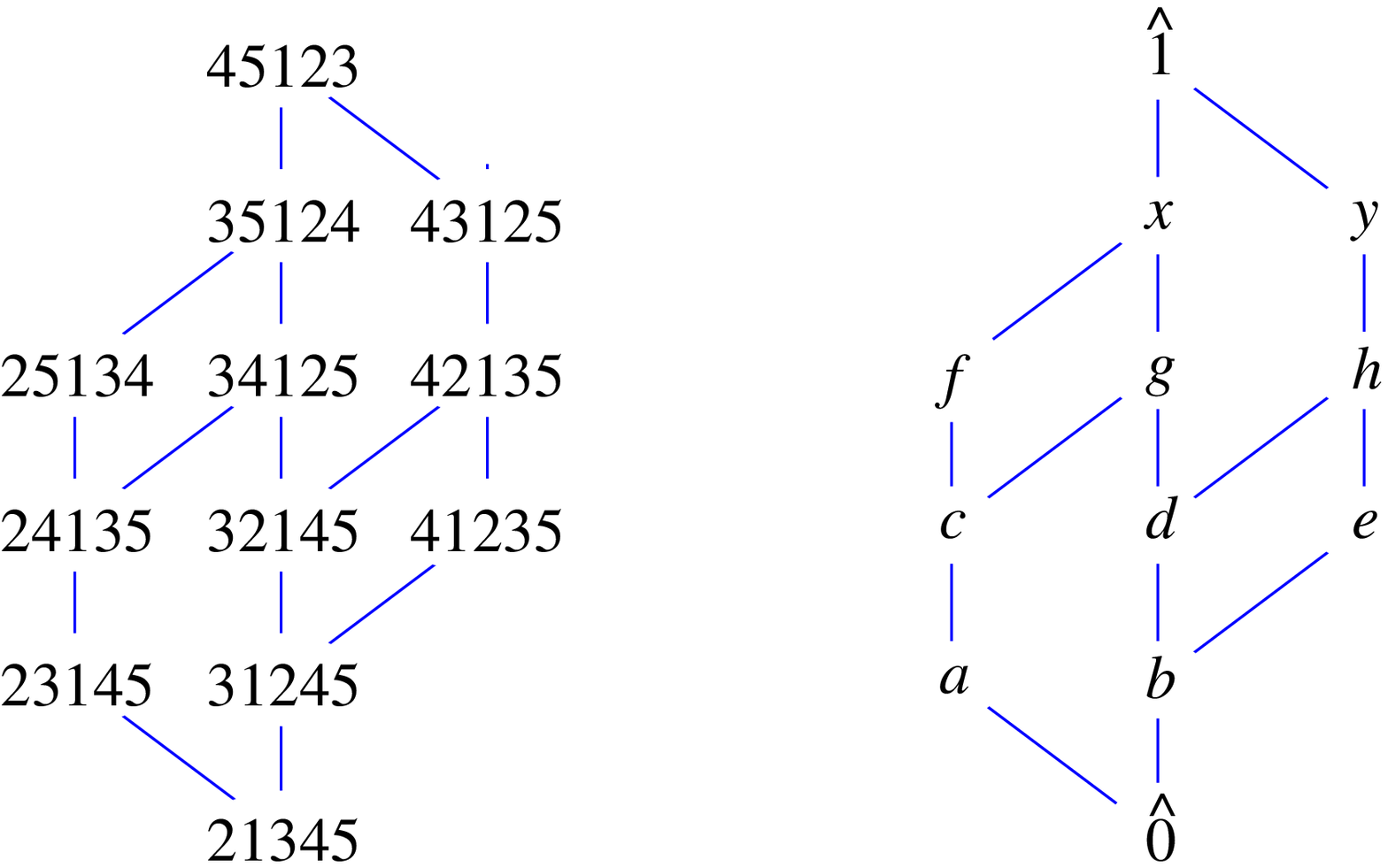}
$$
The simplicial complex $\Delta(P)$ associated to a poset $P$ has as
simplices all chains, including the non-maximal ones.
In our case above, the maximal simplices are
$$
\{a,c,f,x\},\ \{a,c,g,x\},\ \{b,d,g,x\},\ 
\{b,d,h,y\},\ \{b,e,h,y\}.
$$

While $(\{a,c,f,x\}, \{a,c,g,x\})$  and 
$(\{b,d,h,y\},\{b,e,h,y\})$ are attached along facets
($\{a,c,x\}$ and $\{b,h,y\}$, respectively), the pairs
$(\{a,c,g,x\},\{b,d,g,x\})$ and $(\{b,d,g,x\},\{b,d,h,y\})$
are not.
They are attached along codimension 2 faces, 
$\{g,x\}$ and $\{b,d\}$, respectively.
Thus this simplicial complex is not shelable.
Below, we display a geometric realization of this simplicial complex:
$$
\epsfxsize=2.5in \epsfbox{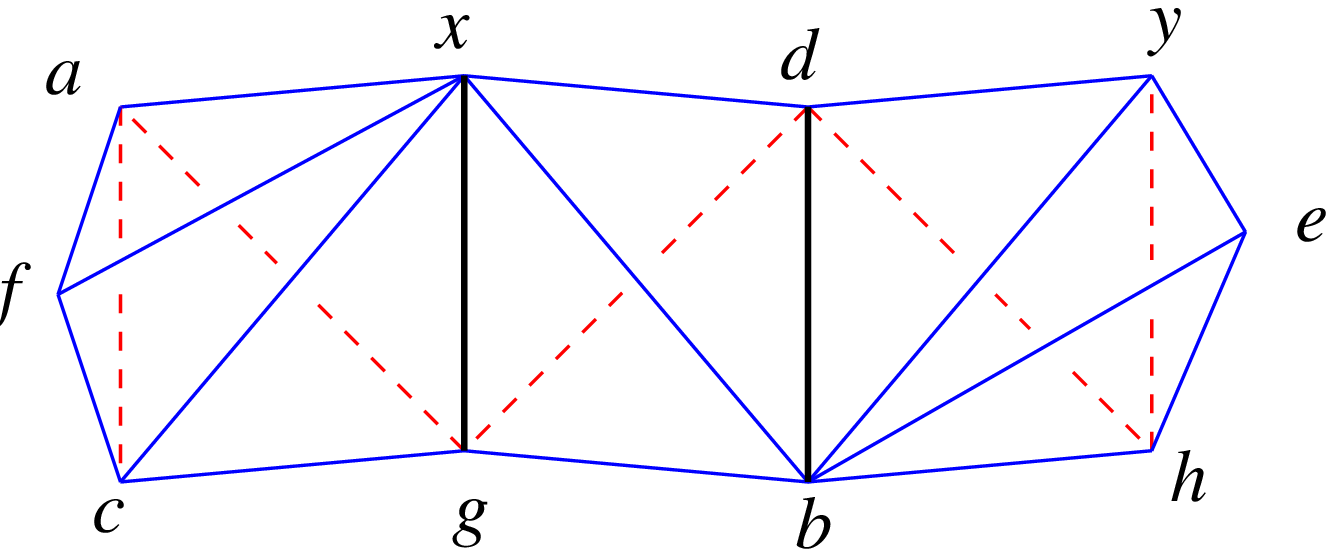}
$$

\end{document}